\documentclass[journal ]{new-aiaa}
\usepackage[utf8]{inputenc}
\usepackage{textcomp}

\usepackage{graphicx}
\usepackage{amsmath}
\usepackage{amsmath,empheq}
\usepackage{hyperref}
\usepackage[version=4]{mhchem}
\usepackage{siunitx}

\usepackage{bm}
\usepackage{amsfonts}

\usepackage{caption}
\usepackage{subcaption}

\usepackage[verbose]{placeins}
\usepackage{soul}
\usepackage{adjustbox}
\usepackage{multirow}
\usepackage{array}
\newcolumntype{L}[1]{>{\raggedright\arraybackslash}p{#1}}
\newcolumntype{C}[1]{>{\centering\arraybackslash}p{#1}}
\usepackage{ulem}
\usepackage{nicefrac, xfrac}

\usepackage{graphicx,import}
\usepackage{pgfplots}
\usepackage{tikz}
\usetikzlibrary{positioning,external}
\pgfplotsset{compat=newest,
every plot/.append style={color=black, mark=none},
every axis/.append style={
label style={font=\fontsize{9pt}{1em}\color{white!15!black}\selectfont},
tick style={font=\fontsize{9pt}{1em}\selectfont, color=black, line cap=round},
ticklabel style= {font=\fontsize{9pt}{1em}\selectfont},
legend style={legend cell align=left, font=\fontsize{9pt}{1em}\selectfont, align=left, draw=white!15!black},
/pgf/number format/1000 sep={},
yticklabel style={
        /pgf/number format/fixed,
        /pgf/number format/precision=5
},
scaled y ticks=false
}}

\usepackage[capitalise]{cleveref}
\usepackage{algorithm}
\usepackage{algpseudocode}

\usepackage[nopostdot,acronym,nogroupskip,nonumberlist]{glossaries}
\newacronym{cr3bp}{CR3BP}{Circular Restricted Three-Body Problem}
\newacronym{bcr4bp}{BCR4BP}{Bi-Circular Restricted Four-Body Problem}
\newacronym{bc}{BC}{Ballistic Capture}
\newacronym{abc}{ABC}{Asymptotic Ballistic Capture}
\newacronym{etd}{ETD}{Energy Transition Domain}
\newacronym{zvc}{ZVC}{Zero Velocity Curves}
\newacronym{zvs}{ZVS}{Zero Velocity Surfaces}
\newacronym{dro}{DRO}{Distant Retrograde Orbit}
\newacronym{soi}{SOI}{Sphere of Influence}
\newacronym{ltb}{LTB}{Lunar Trailblazer}
\newacronym{loi}{LOI}{Lunar Orbit Insertion}
\newacronym{ic}{IC}{Initial Condition}
\newacronym{raan}{RAAN}{right ascension of the ascending node}

\usepackage{enumitem}

\usepackage{lipsum}

\usepackage[version=4]{mhchem}
\usepackage{siunitx}
\usepackage{longtable,tabularx}
\newlength{\minuslength}
\usepackage[toc,page]{appendix} 

\title{Extensive Database of Spatial Ballistic Captures with Application to Lunar Trailblazer}

\author{Lorenzo Anoè \footnote{PhD student, Te P\=unaha \=Atea -  Space Institute, University of Auckland, 20 Symonds Street, Auckland 1010, New Zealand. \\ Corresponding author. Email: \textit{lorenzo.anoe@gmail.com}} }

\author{Roberto Armellin \footnote{Professor, Te P\=unaha \=Atea -  Space Institute, University of Auckland, 20 Symonds Street, Auckland 1010, New Zealand.}}
\affil{The University of Auckland, Auckland, 1010, NZ}

\author{Gregory Lantoine \footnote{Mission Design Engineer, Mission Design and Navigation Section, Jet Propulsion Laboratory, 4800 Oak Grove Drive, Pasadena, CA 91011, USA.}}
\affil{Jet Propulsion Laboratory, California Institute of Technology, Pasadena, California 91011, USA}

\author{Claudio Bombardelli \footnote{Associate Professor, Space Dynamics Group, Technical University of Madrid, Plaza Cardenal Cisneros 4, Madrid 28040, Spain.}}
\affil{Technical University of Madrid, Madrid,  28040, Spain}

\begin{document}

\maketitle

\begin{abstract}
For low-energy missions to the Moon and beyond, Ballistic Capture has proven to be a valuable technique for enabling orbital insertion while alleviating propulsion system requirements. This approach offers two key advantages. First, it extends the insertion window, allowing multiple maneuver opportunities to mitigate potential failures at the nominal insertion point. Second, it enables the required insertion maneuver to be distributed across multiple revolutions, reducing propulsion system constraints in terms of single-burn thrust.
Prior research introduced the concept of Energy Transition Domain to support the creation of a comprehensive database of Ballistic Captures in the planar Circular Restricted Three-Body Problem. However, to apply these trajectories to a real mission scenario, a three-dimensional, spatial analysis and transition to an ephemeris model are necessary. This paper first extends the Energy Transition Domain framework to the spatial case, constructing an extensive database of spatial Ballistic Captures. Then, using Lunar Trailblazer as a case study, a subset of the trajectories is filtered using a mission-specific distance metric, and transitioned into an ephemeris model. Finally, interesting features of this subset are analyzed, and sample high-fidelity trajectories are selected as potential backup options for Lunar Trailblazer.
\end{abstract}

\section{Introduction}

\gls{bc} is a key dynamical process in celestial mechanics and astrodynamics, wherein a spacecraft or small celestial body transitions from an unbound trajectory to a temporary bound orbit around a primary body without requiring an impulsive maneuver. This phenomenon has been widely studied for its implications in both planetary science and space mission design. 

The study of \gls{bc} is particularly relevant for space missions seeking low-energy transfer options. Previous research has shown that \glspl{bc} can significantly reduce the propellant needed for orbit insertion, making them attractive for both interplanetary and lunar missions. A wide range of methods have been proposed to generate such trajectories, including weak stability boundary theory~\cite{Belbruno-Topputo-3, Hyeraci-Topputo2010MethodToDesignBCinER3BP}, the use of periodic orbits~\cite{DeiTos-Russell-Topputo}, and techniques based on invariant manifolds~\cite{koonRoss2001low-energyTransfers, Topputo-Vasile-Bernelli}.
These approaches have been applied in both theoretical studies using various three-body models~\cite{Winter-Neto, Astakhov, LuoTopputo2014constructingBCinRealModel} and in practical contexts, such as low-thrust transfer design~\cite{Low-thrust2BC_Yash, CARLETTA_Mars_low-thr_low-en_mission}, and multi-body dynamics across different planetary systems~\cite{luo2020role}. Notably, they have also found application in real mission scenarios, such as BepiColombo’s transfer to Mercury~\cite{jehn2004low}, and in the study of temporarily captured asteroids~\cite{Urrutxua-Scheeres, Urrutxua-Bombardelli, granvik2012population, fedorets2017AsteoidsDistribution, fedorets2020CD3}.
The application of \gls{bc} techniques to lunar missions is also well established, with several works proposing efficient capture strategies to the Moon~\cite{Griesemer_TargetingLunarBCusingPO, belbruno-lunarBC, SousaTerraCeriotti_FastTransfersBC}.

Despite these advancements, identifying and classifying \gls{bc} trajectories remains a non-systematic task, particularly in higher-fidelity models. A previous work~\cite{BC-journal} introduced the concept of the \gls{etd}, which provides a geometric framework for identifying \gls{bc} initial conditions within the \gls{cr3bp}. The \gls{etd} is defined as the set of points in configuration space, at a fixed Jacobi constant, where the two-body energy with respect to the secondary body is exactly zero. The latter represents a necessary condition that a trajectory must satisfy to be temporarily captured. By focusing on this energy transition surface, the \gls{etd} offers a structured and computationally efficient approach to mapping \gls{bc} trajectories. However, this prior study was limited to a planar formulation of the problem.

The present work extends the \gls{etd} framework to the full spatial \gls{cr3bp}, providing a more comprehensive characterization of \gls{bc} dynamics. The spatial extension allows for the inclusion of out-of-plane motion, which is essential for modeling inclined orbit insertions—particularly relevant for many lunar and planetary exploration scenarios—and leads to a more complete representation of possible capture trajectories.
Additionally, this methodology is applied to a real mission scenario by considering \gls{ltb} as a test case. The goal is to assess whether the database of \gls{bc} trajectories derived from the \gls{etd} can provide practical alternatives for lunar insertion strategies.

To bridge the gap between the simplified \gls{cr3bp} model and operational mission design, the study introduces a transition to an ephemeris model (based on the Spice toolkit~\cite{spice}). This transition is simplified thanks to the \gls{etd} definition, which allows for a straightforward transformation of initial conditions into the ephemeris model, avoiding the need for a tool that adapts entire trajectories. In addition, the transition is facilitated by a tailored distance metric, which identifies \gls{bc} candidates with dynamical properties similar to the nominal \gls{ltb} trajectory. This filtering process significantly reduces the number of candidate trajectories that need to be transitioned, thereby improving computational efficiency. Once transferred to the ephemeris model, suitable trajectories for \gls{ltb} are identified based on specific characteristics. Extended capture durations (i.e., multiple lunar revolutions) can be leveraged to optimize fuel efficiency by exploiting their inherent stability. Moreover, such trajectories could offer increased robustness to miss-thrust events, providing alternative opportunities for lunar insertion. To this end, trajectories featuring repeated close approaches within a short time span are specifically sought. Furthermore, stable \glspl{bc} can support gradual insertion maneuvers across successive revolutions, reducing the need for a large impulsive burn and increasing flexibility in maneuver execution.

In summary, this work builds upon previous work by extending the \gls{etd} to the spatial \gls{cr3bp}, demonstrating its applicability to a real mission, and developing a method for integrating \gls{bc} trajectories into an ephemeris model. The results provide both theoretical insights and practical tools for low-energy mission design, particularly for lunar exploration.

\subsection{Scope and structure of this study}

The first objective of this work is to compute a comprehensive database of \glspl{bc} in the spatial \gls{cr3bp}, building on the planar approach introduced in~\cite{BC-journal}. The second objective is to develop a method for transitioning mission-relevant trajectories into a full ephemeris model. Finally, the third objective focuses on analyzing the resulting trajectories. Key features of the capture set are examined, and representative trajectories are selected based on specific mission constraints, demonstrating the applicability of \glspl{bc} to low-energy, real-mission scenarios.

\begin{figure}[tbp]
    \centering
    \includegraphics[width=0.33\linewidth]{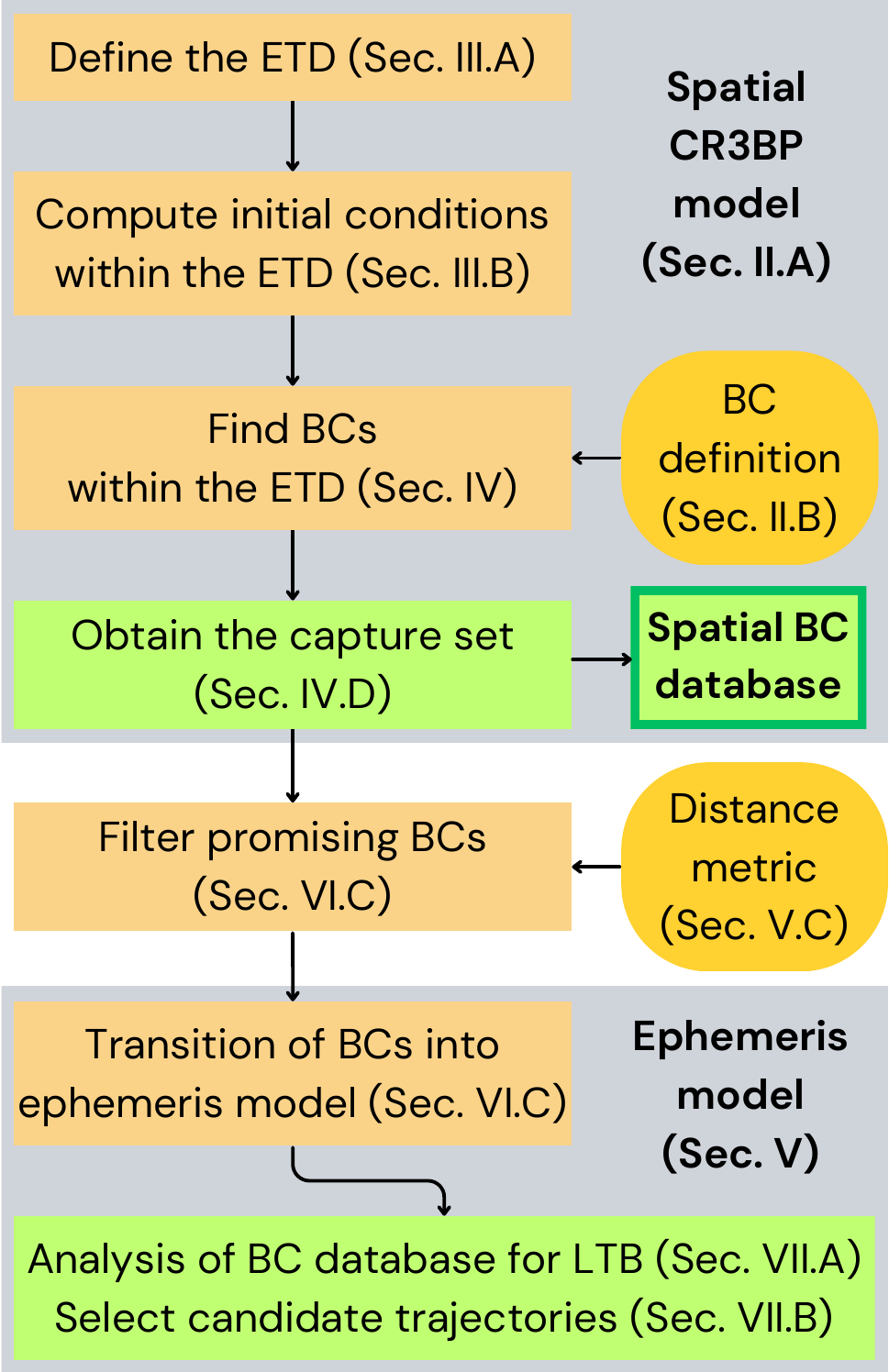}
    \caption{Structure of the paper.}
    \label{fig: flowchart structure}
\end{figure}

The structure of this paper is summarized in the flowchart of \cref{fig: flowchart structure}, and it is discussed here in more detail.
\cref{sec: background} introduces the \gls{cr3bp} model and the definition of \gls{bc} used in this work.
This model will be used to define the \gls{etd}. After addressing the definition of \gls{etd} in the position space (for a fixed three-body energy) in \cref{sec: ETD definition and ICs}, the discussion shifts to determining the possible initial conditions defined by the \gls{etd} itself in \cref{sec: ics in ETD}.
Provided this, it is possible to compute an extensive database of \glspl{bc} in the spatial \gls{cr3bp}. For this purpose, a method (similar to the ones proposed in~\cite{BC-journal}) to provide comprehensive coverage with a reasonable computational cost is introduced in \cref{sec: capture set computation}. Sample results are also provided in \cref{sec: sample spatial capture sets} to highlight the structure of the resulting capture sets. 
In \cref{sec: transition to ephemeris}, a full-ephemeris model and a method to transform initial conditions from the \gls{cr3bp} to this new model are introduced. In addition, a distance metric is also introduced in \cref{sec: distance metric} to preliminary filter all the most promising \glspl{bc} for a specific application.
Using this new high-fidelity model, a test case study is introduced in \cref{sec: test case}. Some generic and mission design related features of the specific missions (i.e. Lunar Trailblazer) are provided, before introducing a method for the transition of the promising \glspl{bc} into the high-fidelity model in \cref{sec: transition to ephemeris algorithm}.
Finally, in \cref{sec: capture set analysis}, the most relevant features of the mission-specific capture set are analyzed. Then, \cref{sec: sample BCs for LTB} presents representative \glspl{bc} that illustrate potential alternatives to Lunar Trailblazer’s nominal insertion strategy. They present a ballistic lunar insertion at the Moon that features promising characteristics such as long time in capture and therefore potentially higher stability, polar revolutions and/or robustness to miss-thrust events.
All the initial conditions for the sample trajectories displayed in this work are provided in Appendix \ref{appendix: sample ICs}.

\section{Background} \label{sec: background}
This section introduces the \gls{cr3bp} and the definition of \gls{bc} as essential elements for establishing the framework of this study.

\subsection{Circular restricted three-body problem} \label{sec: CR3BP}

The \gls{cr3bp} is a fundamental model in celestial mechanics that describes the motion of a body $M_3$ under the gravitational influence of two celestial bodies $M_1$ and $M_2$ (the primaries), with masses $m_1$ and $m_2$, respectively. The mass of the satellite $M_3$ is assumed to be negligible compared to the primaries: $m_3 \ll m_1,m_2$. The primaries are assumed to be point masses located at fixed positions in the synodic rotating frame. In the present analysis, the synodic frame is centered at the barycenter $B$ with the $x-$axis pointing towards $M_2$, providing a convenient framework to study the system~\cite{BattinZVC}. 

While the developed methodology is applicable to any generic system composed of primaries $M_1$ and $M_2$, this study adopts the Earth–Moon system as a case study for all numerical applications and analyses.

\begin{figure}[tbh]
    \centering
    \def\svgwidth{0.6\textwidth}
    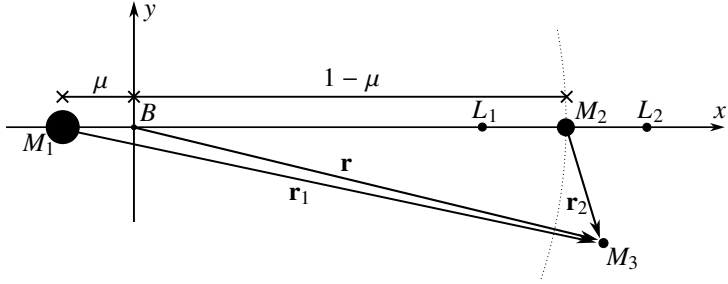
    \caption{\label{fig: sketch CR3BP} \gls{cr3bp} parametrization in a planar $x$–$y$ view.}
\end{figure}

The Cartesian coordinates of the spacecraft are $\mathbf{x}=(\mathbf{r}, \mathbf{v})$, where $\mathbf{r}=(x, y, z)$ and $\mathbf{v}=( \dot{x}, \dot{y}, \dot{z})$. The frame is represented in \cref{fig: sketch CR3BP} for the planar case, with the $z$-axis completing the orthogonal coordinate system.

The $M_1$-$M_2$ distance and inverse mean motion $\tau$ of $M_2$ orbit around the $M_1$ are employed as units of distance (LU) and time (TU), respectively. The mass ratio of the two primaries is denoted as $\mu=m_{2}/(m_{1}+m_{2})$. All the scaling units for the Earth-Moon system are summarized in \cref{tab: scaling units}.

\begin{table}
\caption{\label{tab: scaling units} Approximate scaling units used in this work for the Earth-Moon system.}
\centering
\begin{tabular}{L{2cm} L{3.5cm} L{4cm} L{4.8cm}}
\hline
\noalign{\vskip\doublerulesep}
Unit & Symbol & Value & Note \\ \hline
\noalign{\vskip\doublerulesep}
\noalign{\vskip\doublerulesep}
Mass & $MU=G(m_1+m_2)$ & $4.035032\cdot10^{5} \; km^3 \, s^{-2}$ & System gravitational constant \\
\noalign{\vskip\doublerulesep}
\noalign{\vskip\doublerulesep}
Length & $LU$ & $384399 \; km$ & Mean Earth-Moon distance \\
\noalign{\vskip\doublerulesep}
\noalign{\vskip\doublerulesep}
Time & $TU=\left( LU^3/MU \right)^{0.5}$ & $2.357381 \cdot 10^6 \; s$ & Moon's mean revolution period \\
\noalign{\vskip\doublerulesep}
\noalign{\vskip\doublerulesep}
Velocity & $VU=2\pi LU/TU$ & $1.024548 \; km \, s^{-1}$ & Mean orbital velocity of the Moon \\
\noalign{\vskip\doublerulesep}
\noalign{\vskip\doublerulesep}
Energy & $EU=VU^2=MU/LU$ & $1.049699 \; km^2 \, s^{-2}$ & 
Moon's keplerian energy \\
\noalign{\vskip\doublerulesep}
\end{tabular}
\end{table}

The equations of motion for the spatial \gls{cr3bp} are expressed in terms of the six state variables $(x,y,z,\dot{x},\dot{y},\dot{z})$:
\begin{equation}
\left\{ \begin{array}{l}
\ddot{x} = 2\dot{y} + x - (1-\mu)\cfrac{x+\mu}{(r_1)^3} - \mu \cfrac{(x-(1-\mu))}{(r_2)^3}
\\
\ddot{y} = -2\dot{x} + y - (1-\mu)\cfrac{y}{(r_1)^3} - \mu \cfrac{y}{(r_2)^3}
\\
\ddot{z} = - (1-\mu)\cfrac{z}{(r_1)^3} - \mu \cfrac{z}{(r_2)^3}
\end{array}\right. \, ,
\label{eq: equations of motion}
\end{equation}
where $r_1$ and $r_2$ denote the distance of $M_3$ to $M_1$ and $M_2$ respectively:
\begin{equation}
r_1 = \sqrt{(x+\mu)^2+y^2+z^2} \, ,
\end{equation}
\begin{equation}
r_2 = \sqrt{[x-(1-\mu)]^2+y^2+z^2} \, .
\label{eq: r_2 norm}
\end{equation}

The \gls{cr3bp} is a Hamiltonian system, meaning that a three-body (i.e. total) energy of $M_3$ is conserved in the natural dynamics. This quantity is also known as the Jacobi constant:
\begin{equation}
C_{J} = -v^{2}+x^2+y^2+2\frac{\left(1-\mu\right)}{r_{1}}+2\frac{\mu}{r_{2}} \, ,
\label{eq: CJ}
\end{equation}
where $v = \sqrt{\dot{x}^2+\dot{y}^2+\dot{z}^2}$ represents the norm of $M_3$ velocity.


Five equilibrium points are defined in the \gls{cr3bp} and are commonly referred to as Lagrange points: the three collinear equilibrium points L1, L2 (see \cref{fig: sketch CR3BP}), and L3 lie among the $M_1-M_2$ line. Instead, L4 and L5 are known as triangular points.
The present analysis quantifies the total energy of $M_{3}$ in the synodic frame by employing its own Jacobi constant and those of the Lagrange points L1 and L4. This is achieved through the definition of a three-body energy parameter~\cite{BC-journal}
\begin{equation}
    \Gamma=\frac{C_{J}-C_{J}^{L1}}{C_{J}^{L4}-C_{J}^{L1}} \, .
    \label{eq: gamma}
\end{equation}
This parameter serves as a convenient measure for characterizing the current energy level. It is obtained by comparing and normalizing the energy of $M_{3}$ with respect to two critical values. The first is $C_{J}=C_{J}^{L1}$, corresponding to $\Gamma=0$, which marks the opening of the \gls{cr3bp} forbidden regions at L1. As the three-body energy increases, the parameter $\Gamma$ grows, eventually reaching the second value at $\Gamma=1$, when the forbidden regions vanish at $C_{J}=C_{J}^{L4}$. Moreover, the three-body energy parameter $\Gamma$ facilitates the comparison between different systems and requires fewer digits than the Jacobi constant when reporting the results.

A \gls{bc} is only possible for $\Gamma > 0$, as first introduced by Conley~\cite{conley1968}. In addition, numerical results previously presented~\cite{BC-journal} show that \glspl{bc} (as defined in the next section) are impossible when $\Gamma$ exceeds a limiting value of $\sim 1.36$.


\subsection{Ballistic Capture definition} \label{sec: definition BC}

A clear definition of \gls{bc} is essential for understanding its mechanism~\cite{BC-journal}. Let $\mathbf{x}_0$ be an initial condition of the mass $M_3$ belonging to the \gls{etd} at $\tau_0 =0$, which is a state that satisfies $\varepsilon_2(\tau_0)=0$. When propagated forward (or backward) with \cref{eq: equations of motion}, this initial condition produces the trajectory $\mathbf{x}(\tau)$ for $\tau>0$ (or $\tau<0$).

\underline{Capture state}: The mass $M_3$ is in a capture state at time $\tau$ when its two-body energy is negative, i.e. $\varepsilon_2(\tau)<0$.

\underline{(Temporary) capture phase}: The mass $M_3$ is in a (temporary) capture phase between a start time $\tau_s$ and a final time $\tau_f$ if $\varepsilon_2(\tau_s)=\varepsilon_2(\tau_f)=0$ and $\varepsilon_2(\tau)<0$, for $\tau_s<\tau<\tau_f$. 

\underline{Backwards escape}: Backwards escape occurs when $M_3$ moves away from $M_2$ beyond a certain threshold: $r_2(\tau_e) \geq r_{2,lim} = 0.9$, at time $\tau_e<\tau_0$. Moreover, the trajectory must satisfy $\varepsilon_2(\tau)>0$ for $\tau_{e} \leq \tau<\tau_0$. A similar definition can be applied in the case of forward escape.

\underline{\gls{bc}}: The mass $M_3$ undergoes a \gls{bc} if it remains in a temporary capture phase for a sufficient duration, completing at least one full revolution around the second primary before escaping.
In addition, it must result in an escape backward in time.

\section{The spatial Energy Transition Domain} \label{sec: ETD definition and ICs}

As previously discussed, the \gls{etd} serves as a filter to identify potential \glspl{bc}. In this section, we formalize its mathematical definition and extend its analysis to the spatial case in Cartesian coordinates. This derivation builds upon the equations presented in~\cite{BC-journal}, maintaining consistency with the planar case while revealing new geometric properties. While a detailed discussion of the planar \gls{etd} can be found in~\cite{BC-journal}, here we focus on the spatial extension, briefly referencing established \gls{etd} features — such as three-body energy thresholds and sign-constrained regions — where relevant.

\subsection{ETD definition} \label{sec: ETD}
The first step for the identification of \gls{bc} trajectories is to characterize the \gls{etd}.
The latter is defined as the locus of points in the configuration space (i.e. position space) with zero two-body energy $\varepsilon_{2}$ for a fixed value of the Jacobi constant $C_{J}$.
This section describes how the \gls{etd} can be obtained analytically, and aims to answer the following question: does a point $(x,y,z)$ in the configuration space belong to the \gls{etd}? We show that - given a value of $C_J$ - it is possible to answer the question.

\begin{figure}[tbp]
    \centering
    \def\svgwidth{0.7\textwidth}
    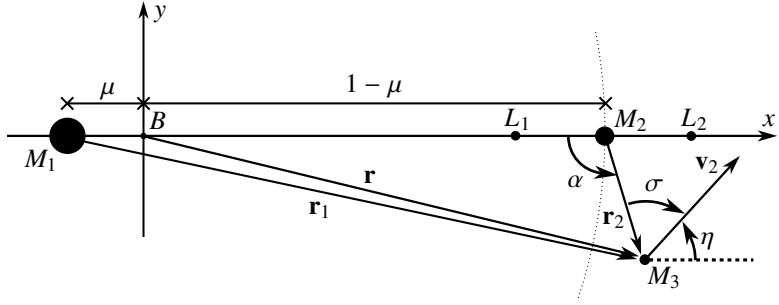
    \caption{\label{fig: sketch params ETD} \gls{etd} parameters in a planar $x$–$y$ view.}
\end{figure}

Building on \cref{fig: sketch CR3BP}, \cref{fig: sketch params ETD} includes the parameters used in this section to define the \gls{etd}. Note that the $z$-axis is again completing the orthogonal coordinate system, and projections of the spatial vectors into the $x$–$y$ plane will be indicated with an (additional) subscript $xy$. Throughout the entire work, we consider an initial time $\tau_0=0$ when the synodic frame is instantaneously aligned with the inertial frames centered at both $M_1$ and $M_2$.
Among other implications, this yields the following transformation:
\begin{equation}
\left\{ \begin{array}{c}
x_2=x-(1-\mu) \\
y_2=y \\
z_2=z
\end{array} \right. .
\label{eq: coord wrt M2}
\end{equation}
To describe the projection of $\mathbf{r}_{2}$ on the $x$–$y$ plane, the variable $\mathbf{r}_{2,xy} = (x_2,y_2,0)$ is introduced.
The angle $\alpha$ is also introduced to define the orientation of $\mathbf{r}_{2,xy}$. Specifically, $\alpha$ is measured counterclockwise from the negative $x-$axis:
\begin{equation}
\cos\alpha=-\frac{x_2}{r_{2,xy}}\qquad\sin\alpha=-\frac{y_2}{r_{2,xy}}
\label{eq: alpha definition}
\end{equation}
To characterize the motion of $M_3$ with respect to an inertial frame centered at $M_2$ (and initially aligned with the synodic frame), the velocity vector $\mathbf{v}_2$ is introduced (with norm $v_2$):
\begin{equation}
\mathbf{v}_{2}=\dot{x}_2\,\mathbf{u}_{x}+\dot{y}_2\,\mathbf{u}_{y}+\dot{z}_2\,\mathbf{u}_{z} \, ,
\label{eq: v2_vec}
\end{equation}
where $\mathbf{u}_{x}$, $\mathbf{u}_{y}$, and $\mathbf{u}_{z}$ are the unit vectors along the three Cartesian axes.
To describe the projection of $\mathbf{v}_{2}$ on the $x$–$y$ plane, the variable $\mathbf{v}_{2,xy} = (\dot{x}_2,\dot{y}_2,0)$ is introduced.
The orientation of the planar $\mathbf{v}_{2,xy}$ is described by the angle $\eta$, which is measured counterclockwise from the positive $x$-axis (see \cref{fig: sketch params ETD}).
Instead, the tilt angle of $\mathbf{v}_2$ along the $z$-axis is represented by $\zeta$.
In other words, $\eta$ is the right ascension and $\zeta$ is the declination of the velocity vector in the velocity space $\dot{x}_2$-$\dot{y}_2$-$\dot{z}_2$.
In this way, the velocity vector components can be parameterized using $\eta$ and $\zeta$:
\begin{equation}
\left\{ \begin{array}{c}
\dot{x}_2 = v_2 \cos{\eta} \cos{\zeta} \\
\dot{y}_2 = v_2 \sin{\eta} \cos{\zeta} \\
\dot{z}_2 = v_2 \sin{\zeta}
\end{array} \right. .
\label{eq: vel eta zeta}
\end{equation}

Finally, the \textit{injection angle} $\sigma$ can be retrieved using the relation: 
\begin{equation}
\alpha=\eta-\sigma \, .
\label{eq: eta sigma}
\end{equation}
Using this configuration, the orientation of the velocity projection $\mathbf{v}_{2,xy}$ is described either by $\eta$ or $\sigma$. Given that the choice is completely arbitrary, $\sigma$ is used while describing the \gls{etd} in conformity with~\cite{BC-journal}. Instead, $\eta$ is used to retrieve the initial conditions in the computational process (see following subsection).

The motion of $M_3$ relative to $M_2$ can be analyzed using the Keplerian two-body energy in an inertial frame centered at $M_{2}$:
\begin{equation}
\varepsilon_{2}=\cfrac{v_{2}^{\,2}}{2}-\cfrac{\mu}{r_2} \, ,
\label{eq: vis viva}
\end{equation}
and imposing zero two-body energy $\varepsilon_{2}=0$, the value of $v_2$ for \cref{eq: vel eta zeta} can be obtained:
\begin{equation}
v_{2}^{\,2}=\cfrac{2\mu}{r_2} \, .
\label{eq: v2_energy_constraint}
\end{equation}
This constraint ensures that the velocity magnitude is exactly providing a transitioning condition that can potentially lead to a capture state with $\varepsilon_{2}<0$ when propagated forward in time.

The velocity transformation from the inertial to the synodic frame is given by:
\begin{equation}
\mathbf{v} = \mathbf{v}_{2} - \mathbf{k}\times\mathbf{r}_{2} = \left(\dot{x}_2+y_2\right)\,\mathbf{u}_{x} + \left(\dot{y}_2-x_2\right)\,\mathbf{u}_{y} + \dot{z}_2\,\mathbf{u}_{z} \, ,
\label{eq: v_vec}
\end{equation}
where $\mathbf{k}$ is the dimensionless angular velocity vector of the synodic frame with respect to an inertial frame, and the $M_{1}-$centered and $M_{2}-$centered position vectors of $M_{3}$ read:
\begin{equation}
\mathbf{r}_{1}=(x+\mu)\,\mathbf{u}_{x}+y\,\mathbf{u}_{y}+z\,\mathbf{u}_{z} \qquad \mathbf{r}_{2}=x_2\,\mathbf{u}_{x}+y_2\,\mathbf{u}_{y}+z_2\,\mathbf{u}_{z}
\label{eq: r_vecs}
\end{equation}
Note that the relationships for $x_2$, $y_2$, and $z_2$ were introduced in \cref{eq: coord wrt M2}.
The inverse relationship of \cref{eq: v_vec} states instead the transformation from the synodic to the inertial frame centered at $M_2$, and it reads:
\begin{equation}
\mathbf{v_2} = \mathbf{v} + \mathbf{k}\times\mathbf{r}_{2} = \left(\dot{x}_2-y_2\right)\,\mathbf{u}_{x} + \left(\dot{y}_2+x_2\right)\,\mathbf{u}_{y} + \dot{z}_2\,\mathbf{u}_{z} \, ,
\label{eq: v_vec syn2inert2}
\end{equation}

\begin{figure}[tbp]
    \centering
    \def\svgwidth{0.49\textwidth}
    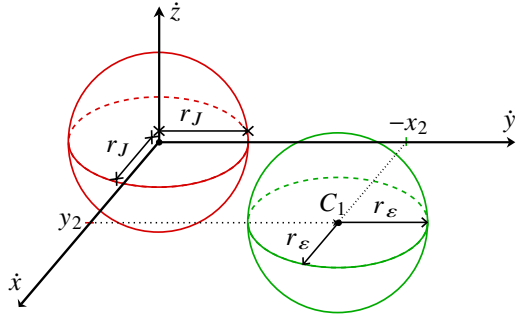
    \caption{Geometric representation of the constraints defining the \gls{etd}.}
    \label{fig: spheres}
\end{figure}

The following equations provide a convenient perspective to analyze the two constraints defining the \gls{etd}, and are obtained from the norm of \cref{eq: v_vec syn2inert2} and \cref{eq: CJ}, respectively. For a given position $(x,y,z)$, the first sets the two-body energy $\varepsilon_2=0$, while the second fixes the value of $C_J$:
\begin{subequations}
  \begin{empheq}[left=\empheqlbrace]{align}
  & (\dot{x}-y_2)^2 + (\dot{y}+x_2)^2 + \dot{z}^2 = 2\mu/r_2 = v_2 = f_1(x,y,z, \varepsilon_2=0) 
  \\
  & \dot{x}^2+\dot{y}^2+\dot{z}^2=f_2(x,y,z,C_J)
  \end{empheq}
\end{subequations}
These can also be written in the form
\begin{subequations}
  \begin{empheq}[left=\empheqlbrace]{align}
  & (\dot{x}-y_2)^2 + (\dot{y}+x_2)^2 + \dot{z}^2 = r_\varepsilon^2 \label{eq: eps2=0 sphere}
  \\
  & \dot{x}^2+\dot{y}^2+\dot{z}^2=r_J^2 \label{eq: CJ=const sphere}
  \end{empheq}
\end{subequations}
where $r_J$ can be obtained from \cref{eq: CJ}.
The two equations represent two spheres in the $(\dot{x}, \dot{y}, \dot{z})$ space with radii $r_\varepsilon$ and $r_J$, as represented in \cref{fig: spheres}. Therefore, this problem can be solved by finding the intersection of two spheres, which is usually a circumference defined in the three-dimensional space. In particular, \cref{eq: CJ=const sphere} represents a sphere centered in the origin, while \cref{eq: eps2=0 sphere} is centered in the point $C_1=\{y_2,\,-x_2,0\}$. The distance of this point from the origin is
\begin{equation}
r_{C1}=\sqrt{x_2^2+y_2^2} \, .
\label{eq: r_C1}
\end{equation}
Note that both spheres lie in the $\dot{x}-\dot{y}$ plane, slightly simplifying the solution process.

When an intersection exists, then the point belongs to the \gls{etd}. Vice versa, a solution for the \gls{etd} cannot be found when the spheres are not intersecting, which can occur only in two cases. First, if the separation between the spheres is bigger than the sum of their radii, hence $r_{C1}>r_J+r_\varepsilon$ (as represented in \cref{fig: spheres}). Or second, if one sphere is completely contained by the other one, hence $r_{C1}<\lVert r_J-r_\varepsilon \rVert$. As a consequence, checking if a point in the physical space belongs to the  \gls{etd} reduces to analyzing these two inequalities.

\begin{figure}[tb]
    \centering
    \includegraphics[width=0.7\linewidth]{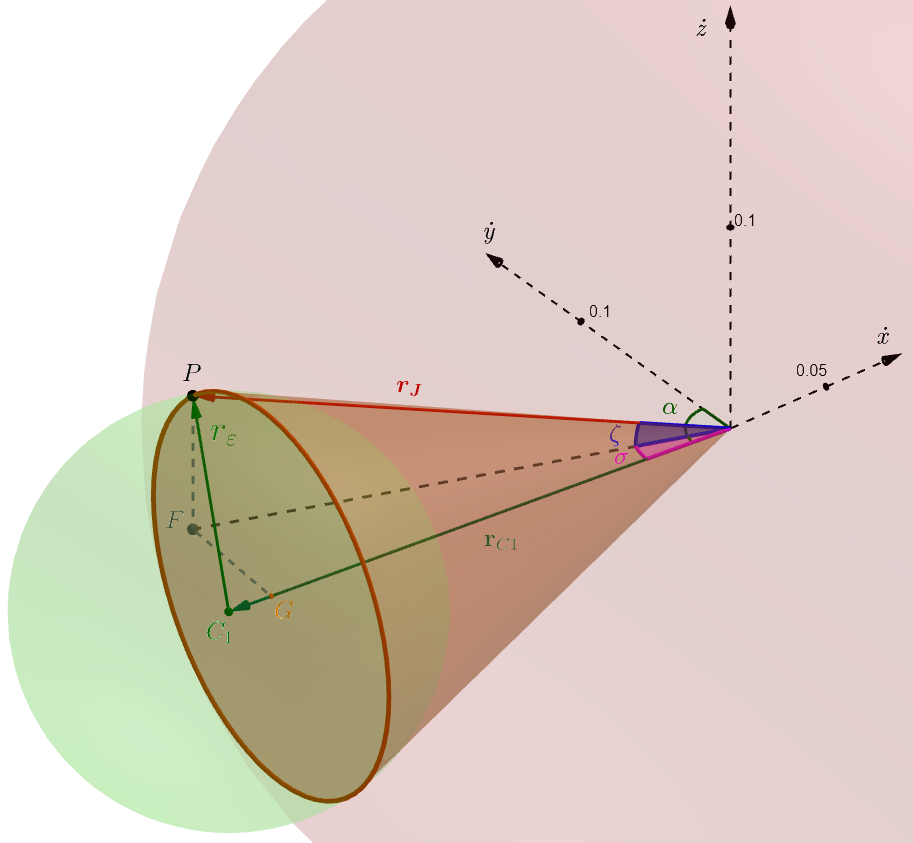}
    \caption{Definition of \gls{etd} through constraints and parameters.}
    \label{fig: spheres Geogebra}
\end{figure}

A three-dimensional representation of this problem is given in \cref{fig: spheres Geogebra}\footnote{An interactive version of this figure
is available online: https://www.geogebra.org/m/jvaxqweb}. Here, the spheres of \cref{eq: eps2=0 sphere} and \cref{eq: CJ=const sphere} are represented in green and red, respectively, for the values contained in \cref{tab: parameters spheres geogebra}. The radii of the two spheres $r_\varepsilon$ and $r_J$ are colored accordingly. The center of the green sphere is indicated with the point $C_1$, whose position taken from the $\dot{y}$ axis is defined by the angle $\alpha$ [see \cref{eq: alpha definition}]. The intersection of the two spheres is represented by the orange circumference centered in the point $G$, which is used to define a cone with vertex in the origin. Therefore, the axis of the cone connects $G$ to the origin. The point $C_1$ can be found by prolonging this axis. Therefore, the vector $\mathbf{r}_{C1}$ connecting the origin with $C_1$ is highlighted in green. The angles $\sigma$ [see \cref{eq: eta sigma}] and $\zeta$ define the lateral surface of the cone starting from the vector $\mathbf{r}_{C1}$, i.e. the orange circumference.
For example, the black point $P$ in \cref{fig: spheres Geogebra} represents a particular velocity that belongs to the set of admissible solutions, i.e. the circumference. The red vector $r_J$ defines this point in the $(\dot{x}, \dot{y}, \dot{z})$ space. Projecting the vector into the $\dot{x}-\dot{y}$ plane, the $F$ point can be obtained. The angle $\sigma$ is defined as the angle spanning between $G$ and $F$, while the angle $\zeta$ spans between $F$ and $P$. A relation between $\sigma$ and $\zeta$ will be introduced in \cref{sec: ics in ETD}. This relation will allow, for a given position $(x,y,z)$, energy $\varepsilon_2 = 0$, and Jacobi constant $C_J$, to describe the full set of initial conditions as a function of $\zeta$, i.e. the circumference in \cref{fig: spheres Geogebra}. In summary, 1) determining if a point in the physical space belongs to the \gls{etd} is equivalent to checking if an intersection of two spheres in the velocity space exists; 2) the determination of the \gls{etd} initial conditions requires computing the intersection between these spheres and a convenient 1D parametrization of the intersection circumference.

Simplifying to the planar case, the intersection of the circumference with the $x$–$y$ plane yields two distinct initial conditions, as introduced in~\cite{BC-journal}.

\begin{table}[tbp]
\caption{\label{tab: parameters spheres geogebra} Approximate values of the parameters for \cref{fig: spheres Geogebra}}
\centering{}
\begin{tabular}{cccccccc}
\hline
\noalign{\vskip\doublerulesep}
$x_2$ & $y$ & $z$ & $r_2$ & $r_\varepsilon$ & $C_J$ & $\Gamma$ & $r_J$ \tabularnewline[\doublerulesep]
\hline
\noalign{\vskip\doublerulesep}
\noalign{\vskip\doublerulesep}
$-0.02$ & $-0.25$ & $0.1$ & $ 0.27$ & $ 0.1$  & $ 2.9880$ & $1.00$ & $ 0.245$ \tabularnewline[\doublerulesep]
\noalign{\vskip\doublerulesep}
\end{tabular}
\end{table}

\begin{figure}[tbp]
    \centering
    \begin{subfigure}[t!]{0.42\textwidth}
        \includegraphics[width=\textwidth]{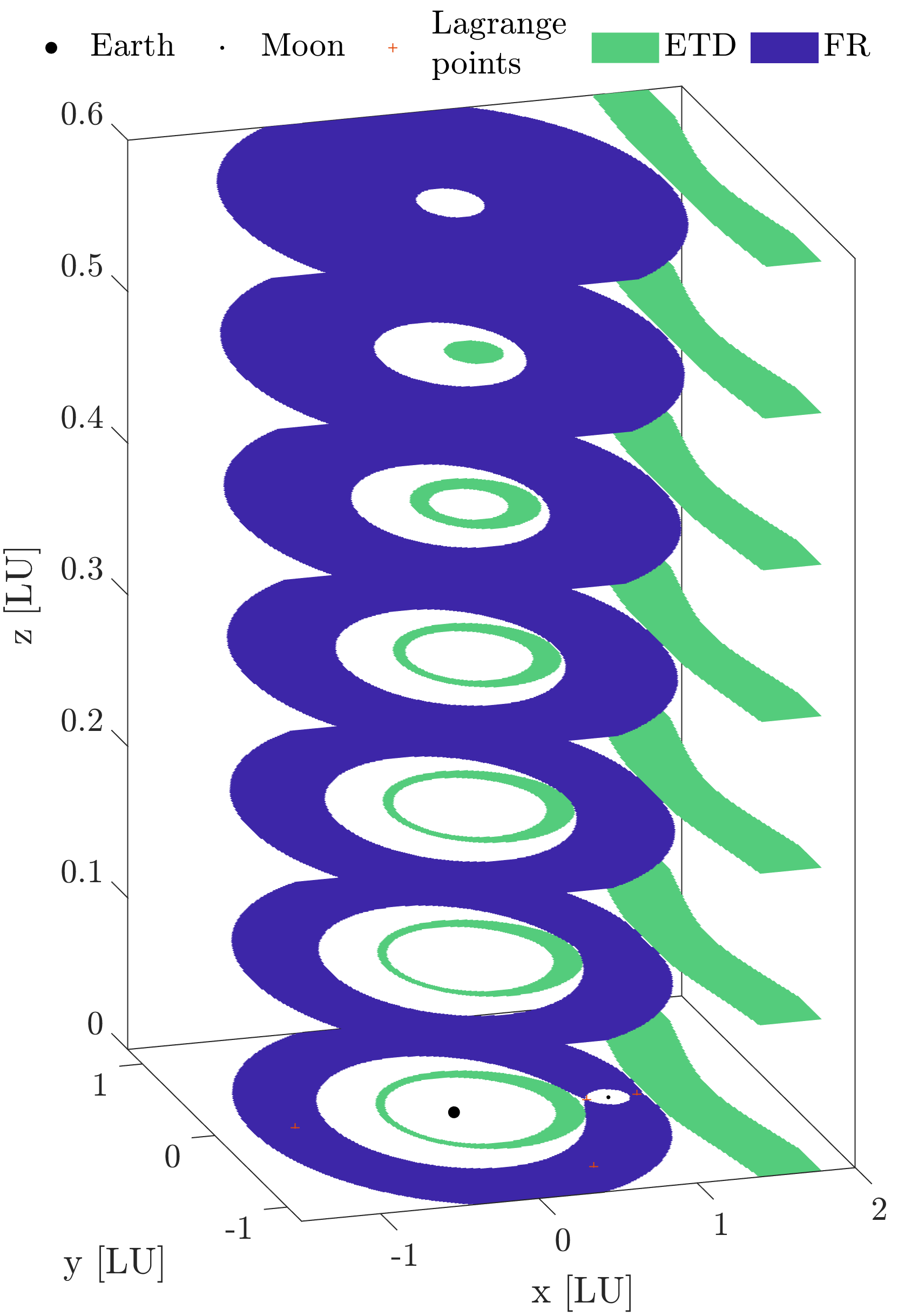}
        \subcaption{$\Gamma = 0.00$}
        \label{fig: ETD dominio Gamma 0}
    \end{subfigure}
    \hspace{5mm}
    \begin{subfigure}[t!]{0.42\textwidth}
        \includegraphics[width=\textwidth]{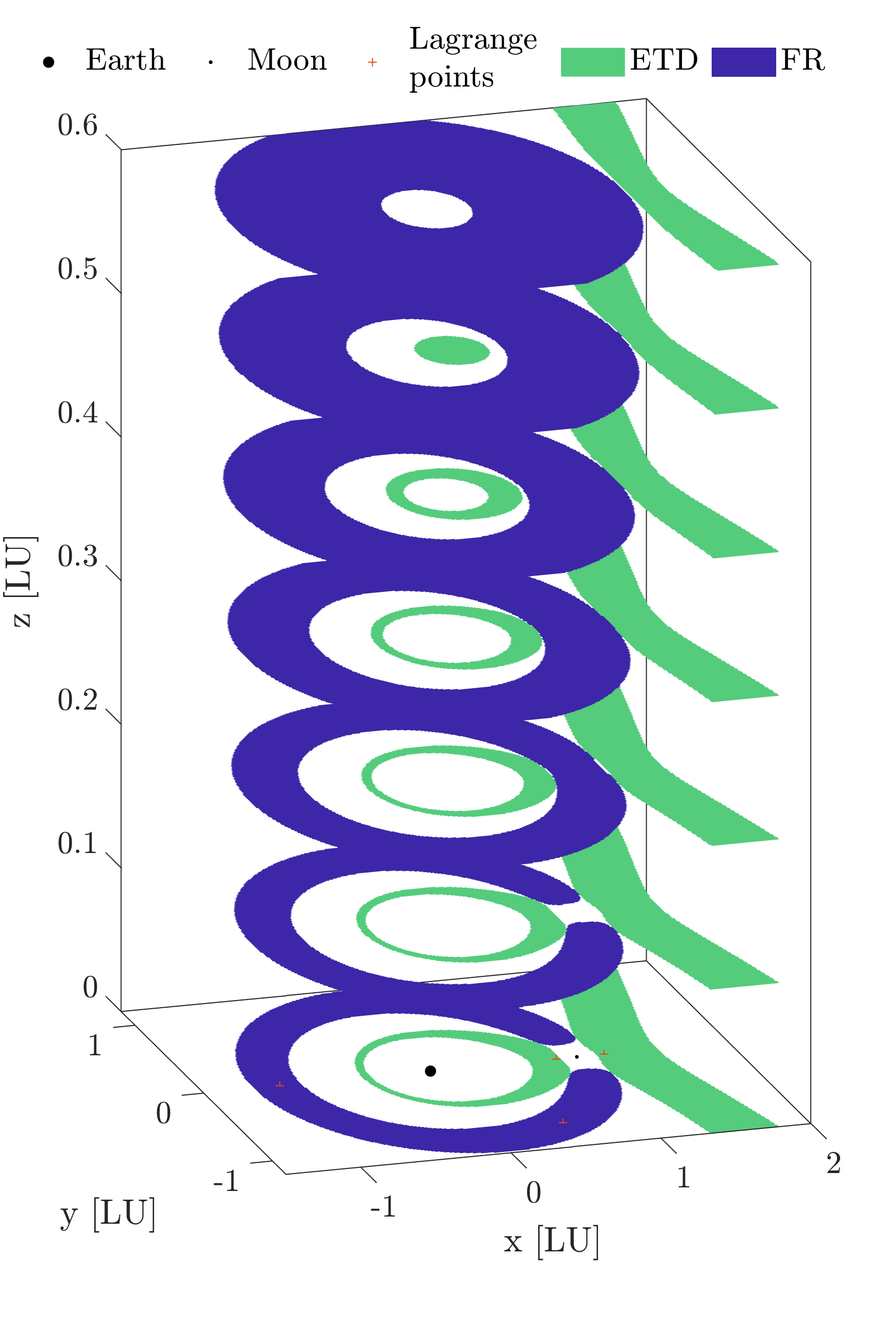}
        \subcaption{$\Gamma = 0.50$}
    \end{subfigure}
    \\
    \vspace{5mm}
    \begin{subfigure}[t!]{0.42\textwidth}
        \includegraphics[width=\textwidth]{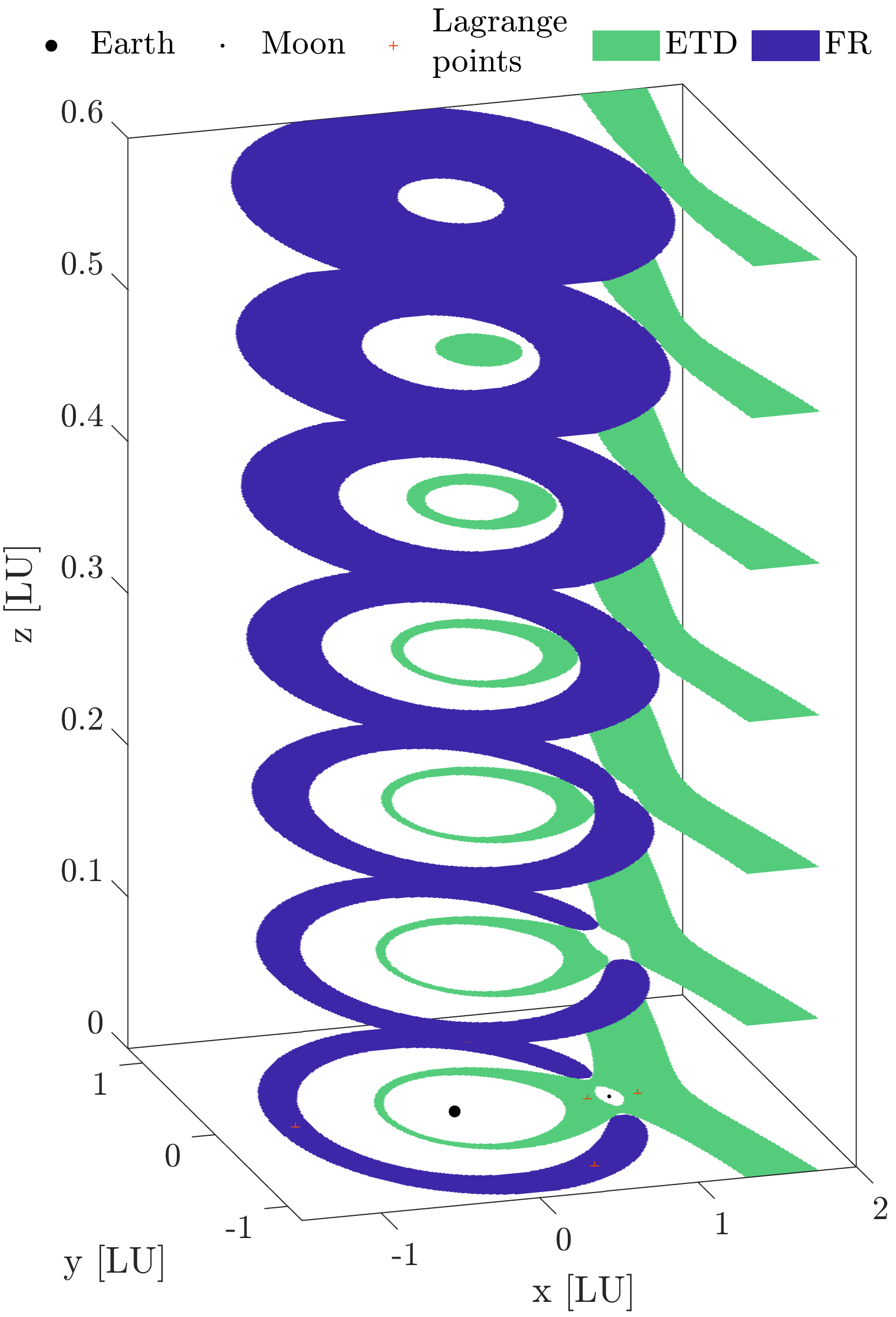}
        \subcaption{$\Gamma = 0.70$}
    \end{subfigure}
    \hspace{5mm}
    \begin{subfigure}[t!]{0.42\textwidth}
        \includegraphics[width=\textwidth]{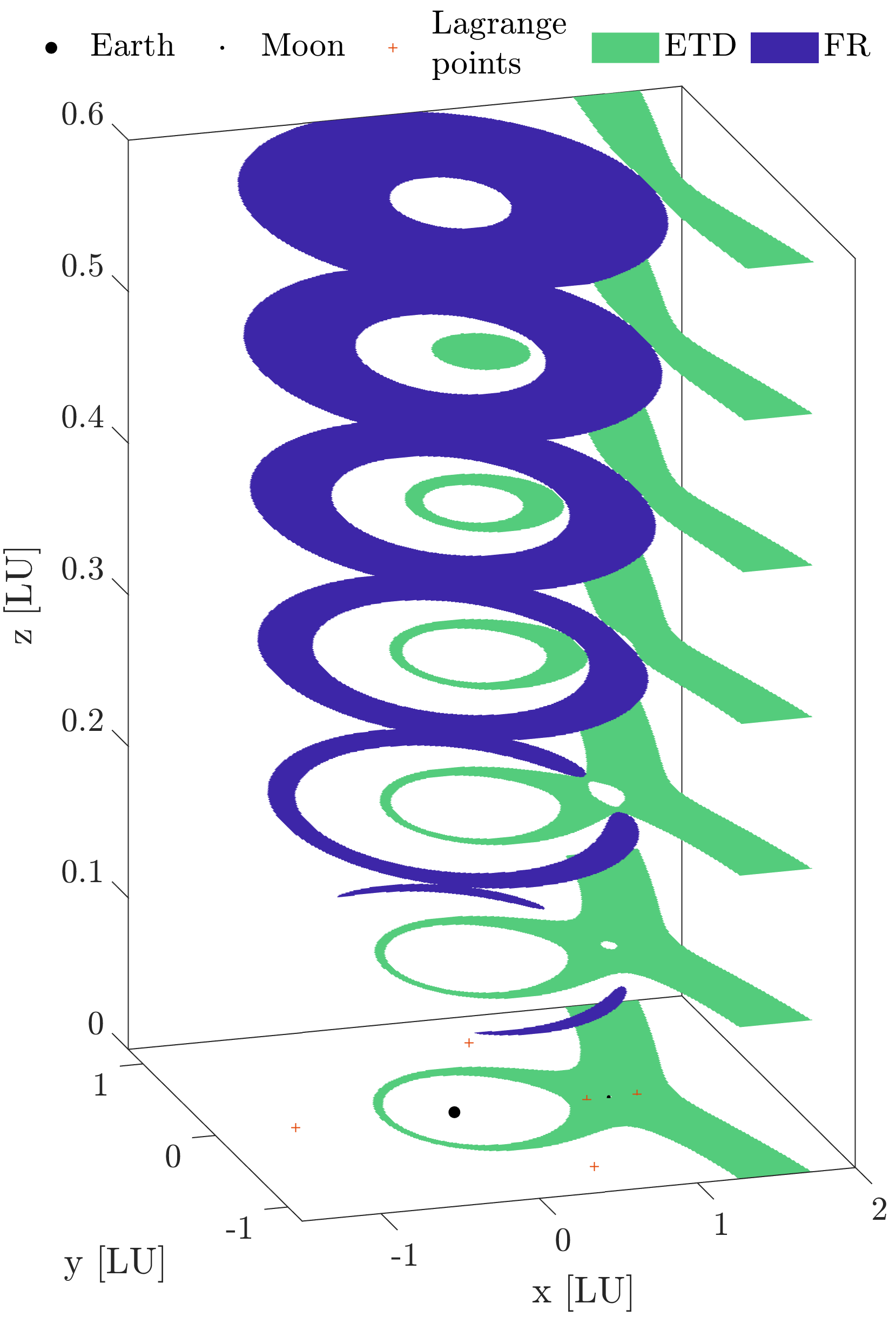}
        \subcaption{$\Gamma = 1.00$}
        \label{fig: ETD dominio Gamma 100}
    \end{subfigure}
    \caption{Evolution of the spatial \gls{etd} with $0\leq\Gamma\leq1$.}
    \label{fig: ETD dominio all 1}
\end{figure}

\begin{figure}[tbp]
    \centering
    \begin{subfigure}[t!]{0.42\textwidth}
        \includegraphics[width=\textwidth]{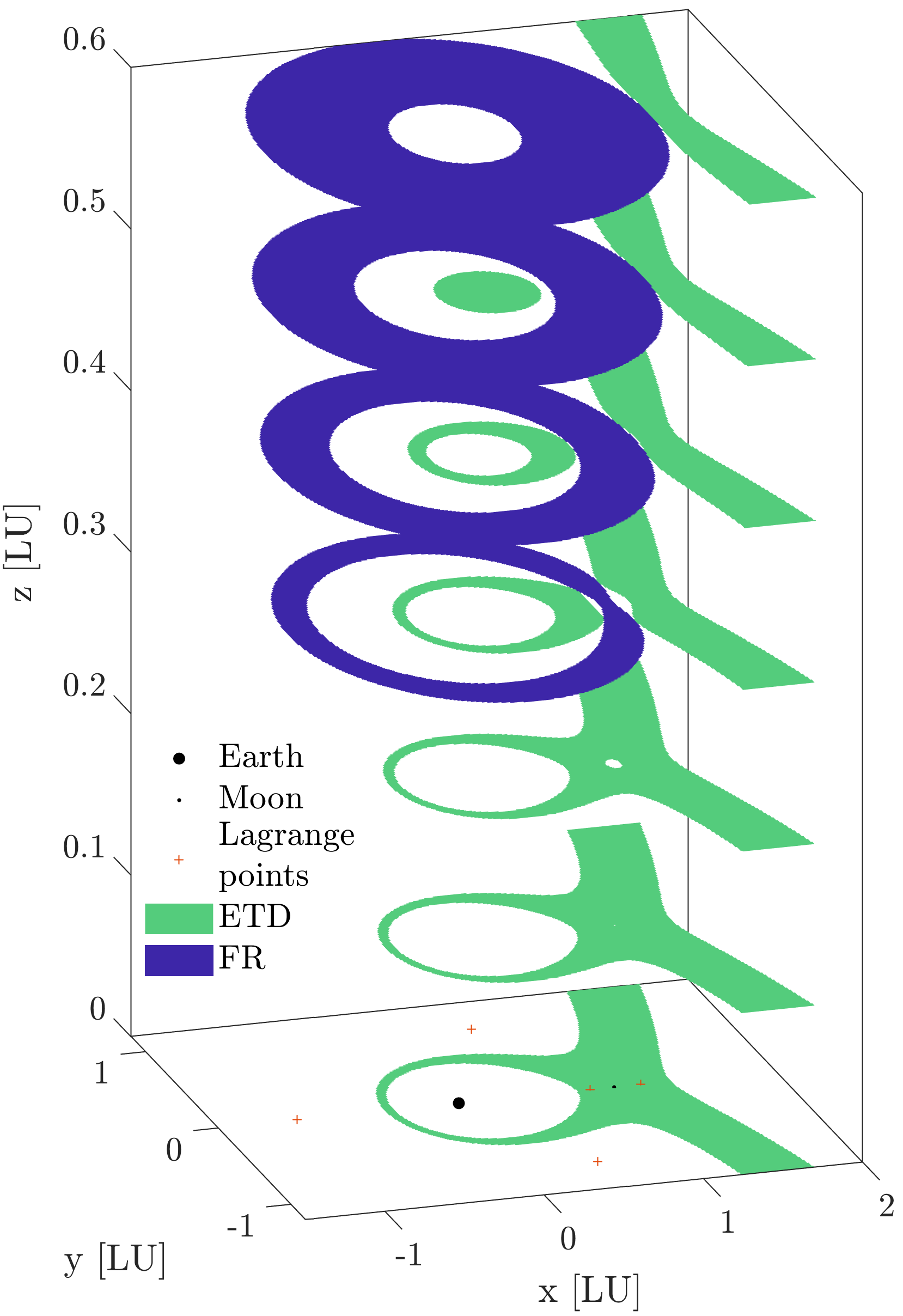}
        \subcaption{$\Gamma = 1.20$}
        \label{fig: ETD dominio Gamma 120}
    \end{subfigure}
    \hspace{5mm}
    \begin{subfigure}[t!]{0.42\textwidth}
        \includegraphics[width=\textwidth]{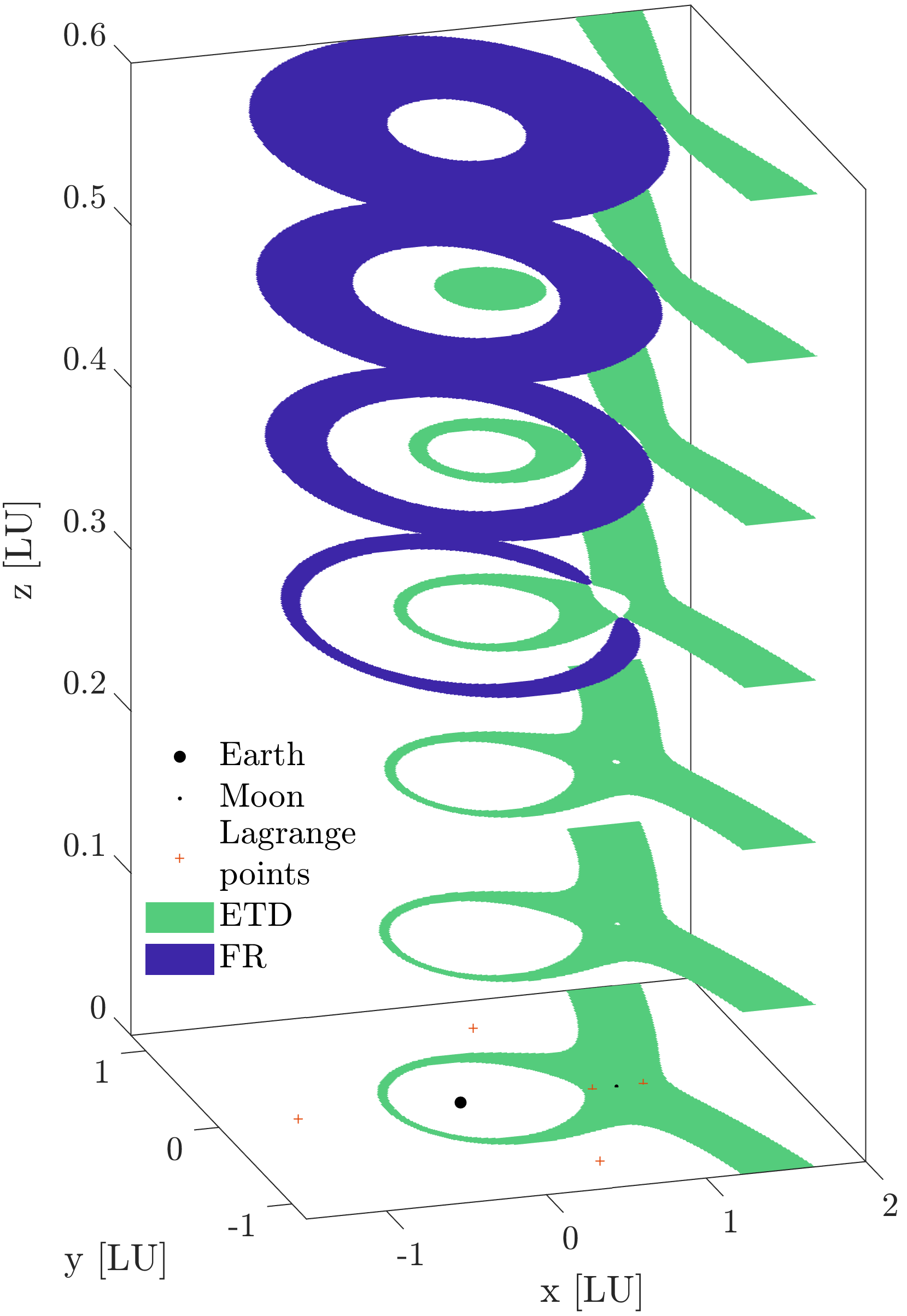}
        \subcaption{$\Gamma = 1.30$}
        \label{fig: ETD dominio Gamma 130}
    \end{subfigure}
    \\
    \vspace{5mm}
    \begin{subfigure}[t!]{0.42\textwidth}
        \includegraphics[width=\textwidth]{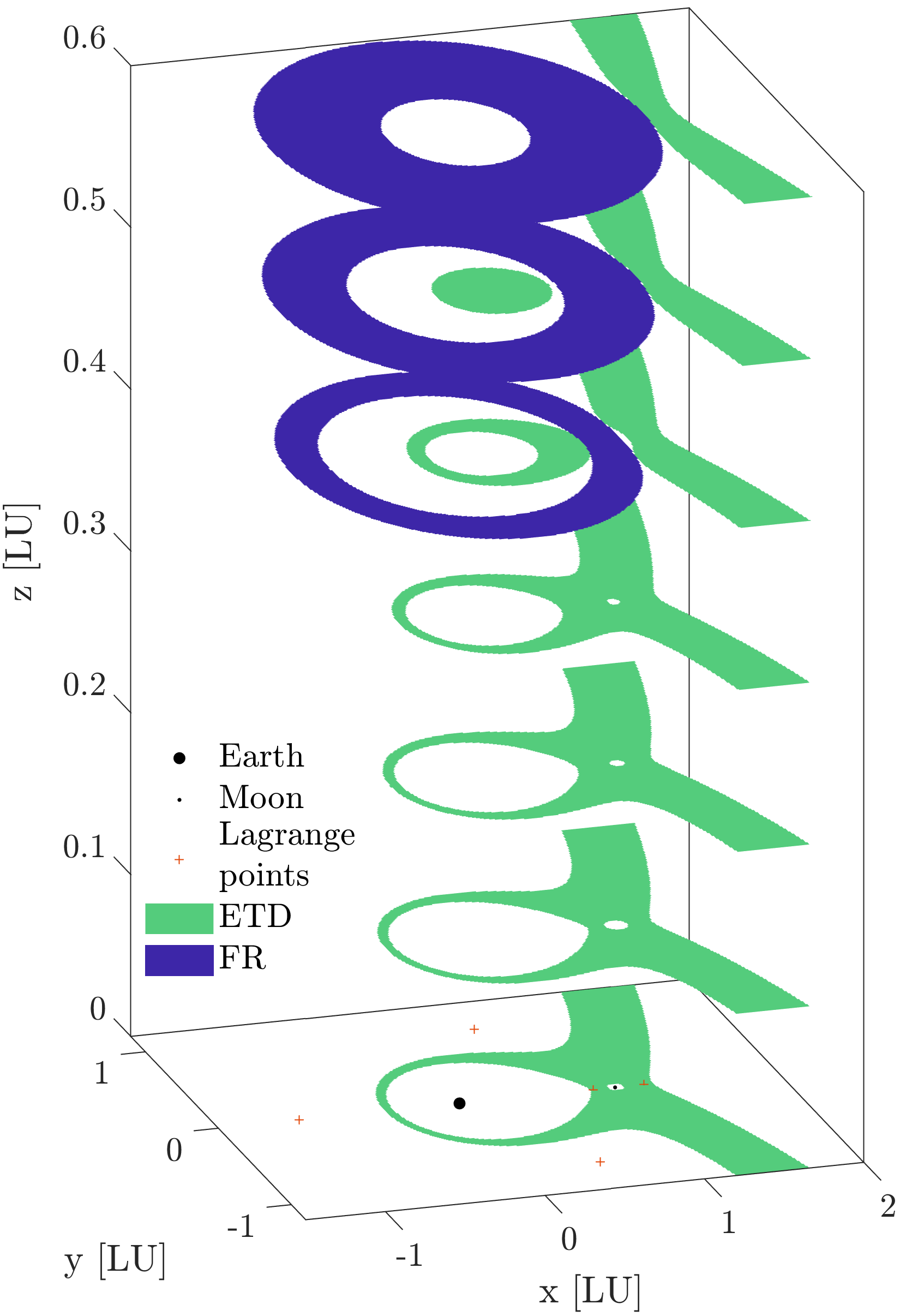}
        \subcaption{$\Gamma = 1.50$}
        \label{fig: ETD dominio Gamma 150}
    \end{subfigure}
    \hspace{5mm}
    \begin{subfigure}[t!]{0.42\textwidth}
        \includegraphics[width=\textwidth]{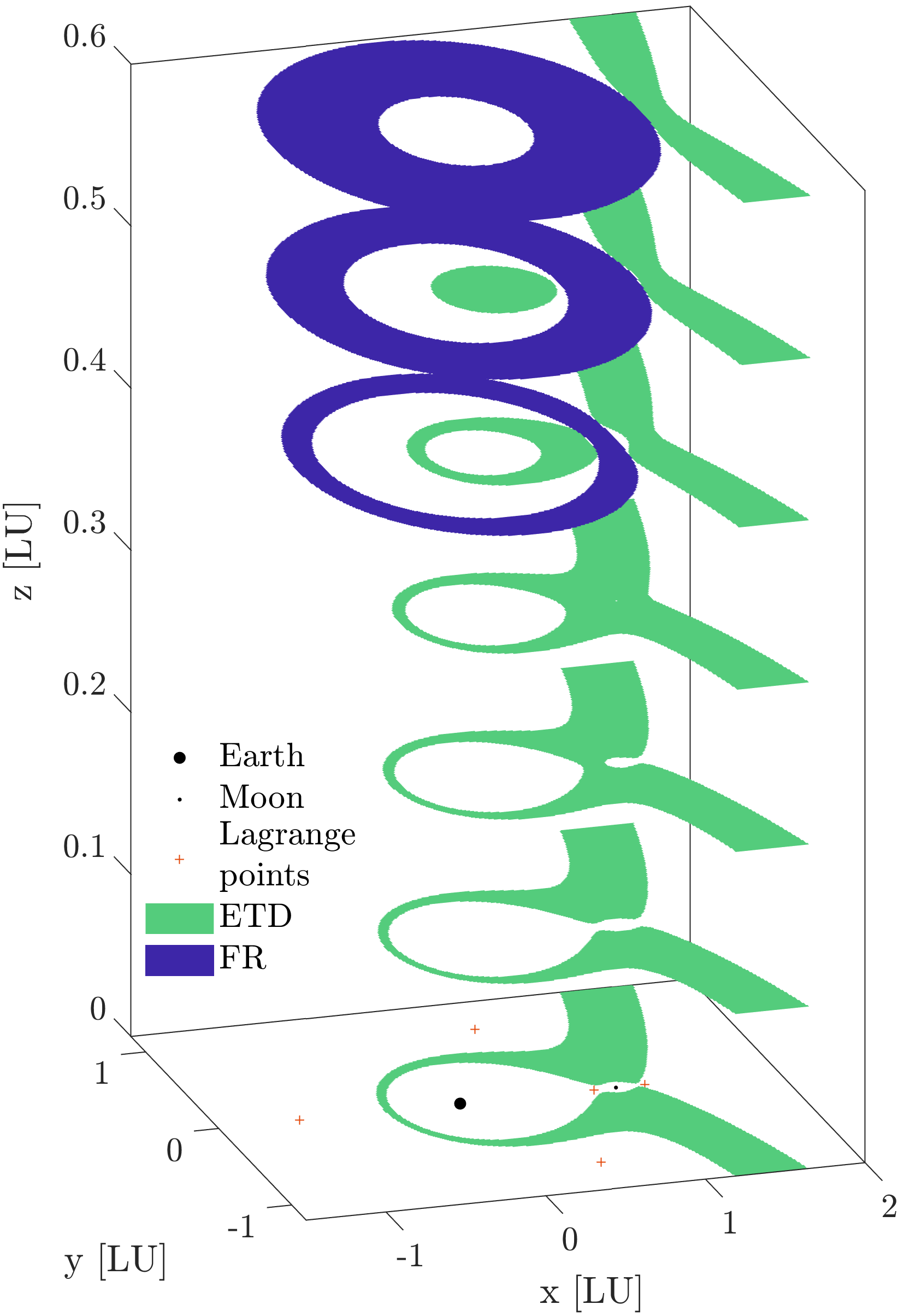}
        \subcaption{$\Gamma = 1.60$}
        \label{fig: ETD dominio Gamma 160}
    \end{subfigure}
    \caption{Evolution of the spatial \gls{etd} with $1.20\leq\Gamma\leq1.60$.}
    \label{fig: ETD dominio all 2}
\end{figure}

\begin{figure}[tb]
    \begin{minipage}[c]{0.49\linewidth}
        \centering
        \includegraphics[width=\linewidth]{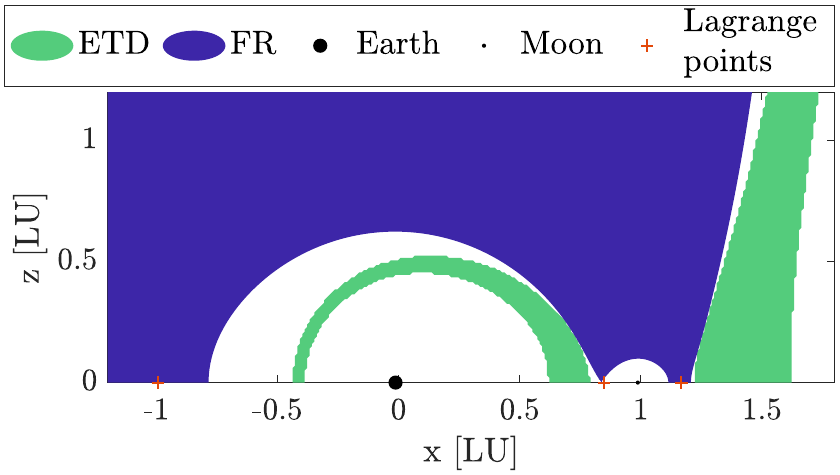}
        \subcaption{$\Gamma = 0$}
        \label{fig: ETD dominio xz section gamma 0}
    \end{minipage}
    \hfill
    \begin{minipage}[c]{0.49\linewidth}
        \includegraphics[width=\linewidth]{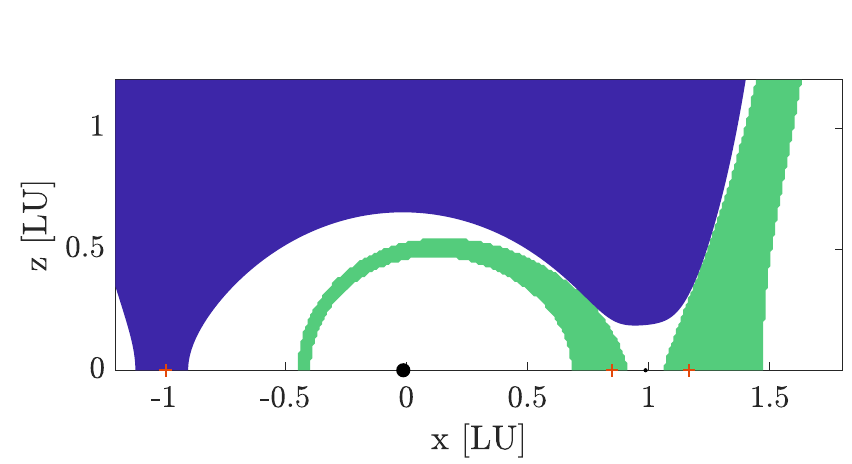}
        \subcaption{$\Gamma = 0.70$}
        \label{fig: ETD dominio xz section gamma 70}
    \end{minipage} \\
    \begin{minipage}[c]{0.49\linewidth}
        \centering
        \includegraphics[width=\linewidth]{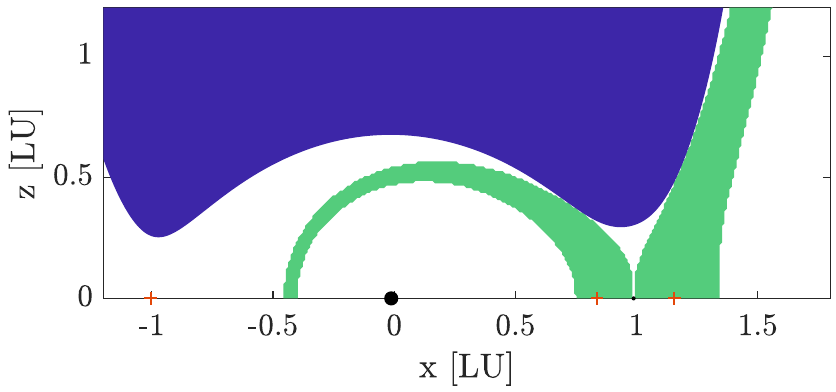}
        \subcaption{$\Gamma = 1.20$}
        \label{fig: ETD dominio xz section gamma 120}
    \end{minipage}
    \hfill
    \begin{minipage}[c]{0.49\linewidth}
        \centering
        \includegraphics[width=\linewidth]{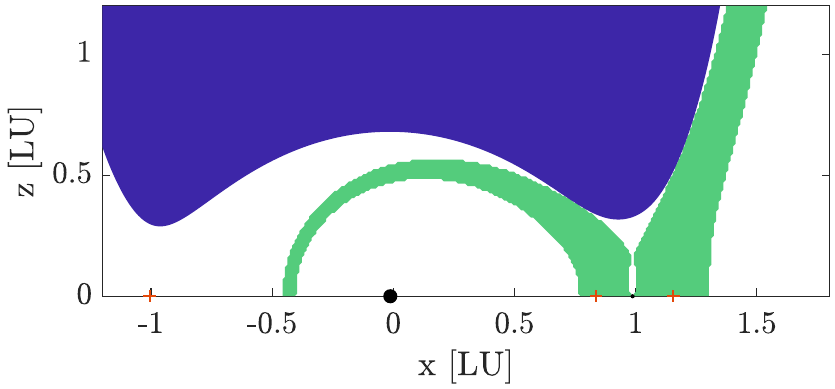}
        \subcaption{$\Gamma = 1.30$}
        \label{fig: ETD dominio xz section gamma 130}
    \end{minipage} \\
    \begin{minipage}[c]{0.49\linewidth}
        \centering
        \includegraphics[width=\linewidth]{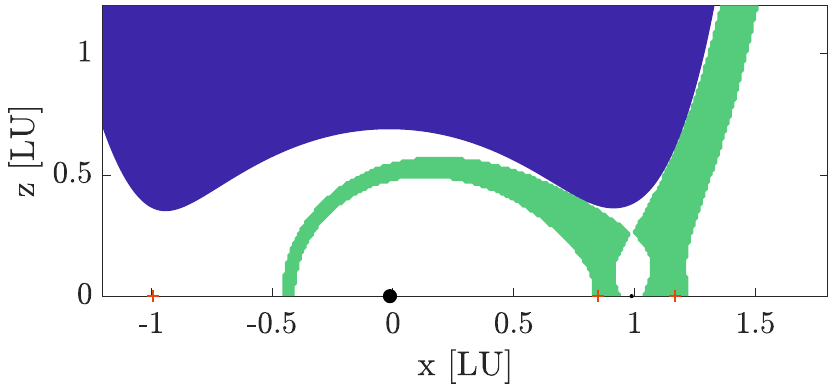}
        \subcaption{$\Gamma = 1.50$}
        \label{fig: ETD dominio xz section gamma 150}
    \end{minipage}
    \hfill
    \begin{minipage}[c]{0.49\linewidth}
        \centering
        \includegraphics[width=\linewidth]{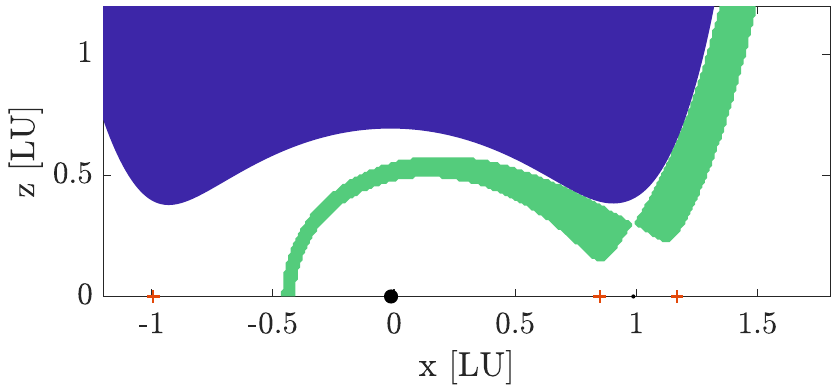}
        \subcaption{$\Gamma = 1.60$}
        \label{fig: ETD dominio xz section gamma 160}
    \end{minipage}
    \caption{\gls{etd} $x$–$z$ cross-sections at $y = 0$.} 
    \label{fig: ETD dominio xz sections}
\end{figure}

The \gls{etd}, i.e. the locus of points for which the intersection between the two spheres exists, is represented in Figs.~\ref{fig: ETD dominio all 1} and \ref{fig: ETD dominio all 2} for different values of $\Gamma$ for the Earth-Moon system. Specifically, $z-$sections of the \gls{etd} are represented, to reveal its spatial structure. Forbidden regions (FR)~\cite{BattinZVC} are also represented. In addition, \cref{fig: ETD dominio xz sections} represents cross-sections on the $x$–$z$ plane at $y=0$.
Two main features of the \gls{etd} - which are evident for low three-body energy level and/or high $z$ values - are that:
\begin{itemize}
    \item a first shell between two ellipsoid-like surfaces always encircles $M_1$;
    \item a second volume included between two planar-like surfaces always extends far out to the right of $M_2$, following the direction of the $y$–$z$ plane.
\end{itemize}
These two volumes are disconnected when $\Gamma=0$, before the opening of the forbidden regions. For increasing $\Gamma$, they expand toward $M_{2}$ and connect. Moreover, it can be demonstrated analytically that they exactly reach the coordinates of $M_2$, that is, $(x, y, z) = (1-\mu,0,0)$, for $C_J=3-4\mu+\mu^2$, which corresponds to $\Gamma$ slightly above 1 (see \cref{fig: ETD dominio Gamma 100}). Above this value of $\Gamma$, the two \gls{etd} volumes start to split again first in the $z$ direction only (e.g., Figs.~\ref{fig: ETD dominio Gamma 120} and \ref{fig: ETD dominio xz section gamma 120}) and then in the $x$–$z$ plane, as it can be seen in Figs.~\ref{fig: ETD dominio Gamma 160} and \ref{fig: ETD dominio xz section gamma 160}.

The behavior of the \gls{etd} along the $z$ axis is rather complicated to describe and can be summarized into three main features.
First, for low three-body energy values (up to $\Gamma\lesssim1.2$) the gap in the \gls{etd} around and above the position of the Moon $(x,y)=(1-\mu,0)$ increases in size when $z$ increases (see \cref{fig: ETD dominio xz section gamma 70}).
Second, for $\Gamma\approx1.2$, the gap in the \gls{etd} around the position of the Moon $(x,y)=(1-\mu,0)$ almost disappears for most of the interesting values of $z$. In fact, in Figs.~\ref{fig: ETD dominio xz section gamma 120} and \ref{fig: ETD dominio xz section gamma 130}, the gap remains stationary in size between $z=0$ and $z\approx0.2$ (see also Figs.~\ref{fig: ETD dominio Gamma 120} and \ref{fig: ETD dominio Gamma 130}).
And third, for higher values of the three-body energy $\Gamma\gtrsim1.2$, the gap in the \gls{etd} around and above the position of the Moon $(x,y)=(1-\mu,0)$ tends to shrink in size when $z$ increases. In fact, observing both Figs.~\ref{fig: ETD dominio Gamma 150} and \ref{fig: ETD dominio Gamma 160}, and Figs.~\ref{fig: ETD dominio xz section gamma 150} and \ref{fig: ETD dominio xz section gamma 160}, the gap shrinks between the values $z=0$ and $z\approx0.3$, and almost disappears for $z=0.3$. This usually occurs up to a value of $z \approx 0.4$ and reveals solutions for the \gls{etd} that are absent in the planar case, highlighting a potential region of interest for promising \glspl{bc} that are associated with high inclinations. For $z \geq 0.4$, the \gls{etd} always comprises the two disconnected main features introduced above. More details are provided in Appendix~\ref{appendix: etd characteristics}.

\subsection{Obtaining initial conditions within the ETD} \label{sec: ics in ETD}
The spatial \gls{cr3bp} has 6 degrees of freedom. Each point in the physical space that belongs to the \gls{etd} has only one degree of freedom, as the two constraints on the two-body energy and Jacobi constant are enforced. The coordinate $\zeta$ that describes the circumference solution set introduced in the previous section represents this degree of freedom. As a consequence, for a given position $(x,y,z)$ and a selected pair $(C_J, \zeta)$, an initial condition can be computed for a potential \gls{bc}. 

Following from \cref{eq: v_vec}, the magnitude of the velocity in the synodic frame hence reads:
\begin{equation}
v^{2} = \dot{x}^2+\dot{y}^2+\dot{z}^2 = \left(\dot{x}_2+y_2\right)^2 + \left(\dot{y}_2-x_2\right)^2 + \dot{z}_2^2 \, ,
\label{eq: v_mag}
\end{equation}
and using the relations of \cref{eq: vel eta zeta}, it becomes
\begin{equation}
v^{2} = \left(v_2 \cos{\eta} \cos{\zeta}+y_2\right)^2 + \left(v_2 \sin{\eta} \cos{\zeta}-x_2\right)^2 + \left(v_2 \sin{\zeta}\right)^2 \, .
\end{equation}
By substituting this equation in \cref{eq: CJ}, a relation in the form $a \sin{\eta} + b \cos{\eta} = c$ can be obtained after some algebra. In particular, this leads to
\begin{equation}
\left\{ \begin{array}{l}
a = -x_2 \\ 
b = y_2 \\ 
c = \cfrac{2(1-\mu)/r_1-2(1-\mu)x-(1-\mu)^2-C_J}{2v_2\cos{\zeta}}
\end{array} \right. .
\label{eq: added angle param}
\end{equation}
This equation can be solved analytically, for example by the added angle method, leading to the solution
\begin{equation}
\sin{(\eta-\alpha)} = c/A \, , 
\label{eq: added angle sol formula}
\end{equation}
where $A=\sqrt{a^2+b^2}$, and $\alpha$ is given by Eqs.~\eqref{eq: alpha definition}.
The variable $\eta$ can be easily obtained with
\begin{equation}
\eta = \arcsin{(c/A)}+\alpha \, .
\label{eq: eta sol}
\end{equation}
Thus, for a given $\zeta$, initial conditions belonging to the \gls{etd} can be determined by solving for $\eta$ in \cref{eq: eta sol}, then substituting into \cref{eq: vel eta zeta} and subsequently into \cref{eq: v_vec}.
Following this process it is possible to retrieve an initial condition for a potential \gls{bc}, given a position and a selected combination $(C_J, \zeta)$. Specifically, this initial condition satisfies the transition condition $\varepsilon_2=0$, ensuring the possibility of a capture phase when propagated forward in time, as $\sim50\%$ of such cases lead to a bound trajectory with $\varepsilon_2<0$.

Comparing Eqs.~\eqref{eq: eta sigma} and \eqref{eq: added angle sol formula}, the \textit{injection angle} $\sigma$ (\cref{fig: sketch params ETD}) is found as:
\begin{equation}
\sin{\sigma} = c/A = \cfrac{2(1-\mu)/r_1-2(1-\mu)x-(1-\mu)^2-C_J}{2v_2\cos{\zeta}\sqrt{x_2^2+y_2^2}} \, .
\label{eq: injection_angle}
\end{equation}
As previously introduced in~\cite{BC-journal}, the \gls{etd} can also be described by $-1 \leq \sin\sigma \leq 1$, hence:
\begin{equation}
-1 \leq \cfrac{2(1-\mu)/r_1-2(1-\mu)x-(1-\mu)^2-C_J}{2v_2\cos{\zeta}\sqrt{x_2^2+y_2^2}} \leq 1 \, .
\end{equation}
This equation for a value of $\zeta=0$ leads to the same solutions presented in \cref{sec: ETD}.

In summary, each \gls{etd} initial condition $\mathbf{x}_0 = \left( x, y, z, \dot{x}, \dot{y}, \dot{z} \right)$ can be described using the five variables $\left( x, y, z, C_J, \zeta \right)$ through Eqs.~\eqref{eq: eta sol}, \eqref{eq: vel eta zeta}, \eqref{eq: v2_energy_constraint}, and \eqref{eq: v_vec}. In the following, the parameter $\Gamma$ will replace $C_J$ for the reasons introduced in \cref{sec: CR3BP}. 
It follows that the \glspl{bc} can be represented in the $\left( x, y \right)$ section of the configuration space for a selected triplet $\left( \Gamma, z, \zeta \right)$. 
Finally, it is important to note that \cref{eq: eta sol} admits two solutions $(\eta_1, \eta_2)$. This means that two distinct initial conditions $\mathbf{x}_{0,1}=f\left( x, y, z, \eta_1, \zeta \right)$ and $\mathbf{x}_{0,2}=f\left( x, y, z, \eta_2, \zeta \right)$ can be obtained for every combination of $(x,y,z,\Gamma,\zeta)$.

\section{ETD Ballistic Capture set} \label{sec: capture set computation}
The focus shifts now to systematically compute and classify \gls{bc} trajectories. To achieve this, an extensive database of capture trajectories is constructed by propagating a large set of initial conditions and identifying those that satisfy the \gls{bc} definition. Given the high dimensionality of the problem, an efficient method is required to ensure comprehensive coverage of the entire \gls{etd} with a reasonable computational time.

Once the initial conditions are defined following the procedure of \cref{sec: ics in ETD}, the \gls{etd} Ballistic Capture set $\mathcal{C}(\Gamma,z,\zeta)$ can be computed. Using the definition in~\cref{sec: definition BC}, for each $(\Gamma,z,\zeta)$ a capture region is identified in the $(x,y)$ plane.
For brevity, we refer to this set simply as the ``capture set'' throughout the rest of this section.

To efficiently compute the capture sets, we build upon a previously introduced strategy for the planar problem~\cite{BC-journal}. This study extends that approach to the full spatial problem by incorporating the additional degrees of freedom $z$ and $\zeta$, generalizing the computation of $\mathcal{C}(\Gamma,z,\zeta)$.

\subsection{Planar method}
The previously proposed strategy for characterizing $\mathcal{C}(\Gamma,z=0,\zeta=0)$ in the planar problem~\cite{BC-journal} began with the identification of a \gls{bc} \textit{kernel} where $\mathcal{C}(\Gamma,z=0,\zeta=0)$ reduces to a point for $\Gamma\rightarrow0$. As $\Gamma$ increases, the capture set $\mathcal{C}(\Gamma,z=0,\zeta=0)$ evolves into a more complex morphology described in detail in~\cite{BC-journal}.

\subsection{Planar to spatial extension method} \label{sec: planar2spatial4symm}
The extension to the full spatial problem is facilitated by the flexibility of the previously developed polygonal-based algorithm~\cite{BC-journal}, which allows for a seamless adaptation to higher-dimensional cases. More specifically, the method relies on a polygonal description of the capture set boundary (using the \textit{polyshape} MATLAB\copyright \, function) and tracks its evolution for increasing values of $\Gamma$, $z$, or $\zeta$, with only minor modifications to the original approach. In this way, the sequence $\mathcal{C}_{i,j,k} = \mathcal{C}(\Gamma_i, z_j, \zeta_k)$ for $i = 1,\dots, N$, $j = 1,\dots, M$, $k = 1,\dots, L$ can be obtained. $\Gamma_{max} = \Gamma_{N}$ is the value of the three-body energy parameter above which \glspl{bc} are no longer found, whereas the maximum values for the remaining two parameters are $z_{max} = z_{M} = f(\Gamma_i)$ and $\zeta_{max} = \zeta_{L} = f(\Gamma_i, z_j)$.
Following the previous work, the algorithm reduces the number of candidate \gls{bc} trajectories by constraining the search space $\mathcal{S} \subset (x,y) \;\;\forall(\Gamma,z,\zeta)$ to be considered.
In the previous work, at each energy step $\Delta \Gamma$, the search space $\mathcal{S}(\Gamma_i)$ was incrementally expanded starting from $\mathcal{C}_{i-1}$ (with the \textit{polybuffer} function) until no more \gls{bc} could be found within its boundaries $\partial \mathcal{S}_i$. Here, the same procedure is repeated first for each step $\Delta z$ in \cref{alg: search algorithm z}, and then for each step $\Delta \zeta$ in \cref{alg: search algorithm zdot}. Note that \cref{alg: search algorithm z} relies on the capture sets computed by the planar algorithm~\cite{BC-journal}, i.e. $\mathcal{C}_{i,j=0,k=0} = \mathcal{C}(\Gamma_i,z=0, \zeta=0)$. In a similar way, \cref{alg: search algorithm zdot} relies on the capture sets computed by \cref{alg: search algorithm z}, i.e. $\mathcal{C}_{i,j,k=0} = \mathcal{C}(\Gamma_i,z_j, \zeta=0)$.

To reduce the computational cost, only the trajectories with $\dot{\varepsilon}_2(\tau_0)<0$ are propagated, where $\dot{\varepsilon}_2 = \dfrac{d\varepsilon_2}{d\tau}$. This ensures that only trajectories evolving toward a negative two-body energy $\varepsilon_2<0$ are considered, making them potential candidates for capture. Additionally, these trajectories lead to escape ($\varepsilon_2>0$) in backward time. Conversely, an initial condition $\mathbf{x}_0$ with $\dot{\varepsilon}_2(\tau_0)>0$ implies that its trajectory is generated by at least one other $\mathbf{x}_0 \in \mathrm{ETD}$ for $\tau<0$. Therefore, this initial condition is discarded as it is linked to another initial condition in the \gls{etd} which, in turn, can potentially satisfy the \gls{bc} definition.
The condition on $\dot{\varepsilon}_2(\tau_0)<0$ is achieved by imposing
\begin{equation}
\mathbf{a}_{3b}\cdot\mathbf{v}_{2}<0 \, ,
\label{eq: check derivative eps2}
\end{equation}
where $\mathbf{a}_{3b}$ is seen as the perturbation of $M_1$ on a orbit around $M_2$ with instantaneous velocity $v_2$, namely:
\begin{equation}
\mathbf{a}_{3b}=(1-\mu)\Bigg(-\cfrac{\mathbf{r}_{1}}{r_1^{3}}+\mathbf{r}_{1}-\mathbf{r}_{2}\Bigg) \, .
\label{eq: a_p perturbation M1}
\end{equation}

Finally, all propagations are stopped when the distance $r_2$ from the $M_2$ is lower than its physical radius $r_M$ (in which case the trajectory is cataloged as a collision). 
The same occurs - in line with what is introduced in \cref{sec: definition BC} - when $r_2 \geq r_{2,lim} = 0.9$ LU, as it is almost impossible for a \gls{bc} to be maintained at such a large distance. In line with~\cite{Belbruno-Topputo-3}, this entails that no more capture phases follow before a full revolution (or a segment of a bounded orbit) around $M_1$.

\begin{algorithm}[tbp]
\caption{\gls{etd} Ballistic Capture Set Computation: $z$-sections}
\label{alg: search algorithm z}
\begin{algorithmic}[1]
\State {Build a grid $\mathcal{G}$ around $M_2$ with width $x_{\mathcal{G}}$, height $y_{\mathcal{G}}$, and  stepsize $h$. Initialize $\Gamma=0$, $i=0$.}
\While {$\Gamma < \Gamma_{max}$}
    \State {$i \leftarrow i+1$}
    \State {$\Gamma_i \leftarrow \Gamma_{i-1}+\Delta \Gamma$}
    \State {Set $j\gets0$, $k\gets0$, and load the planar capture set $\mathcal{C}_{i,j=0,k=0} = \mathcal{C}(\Gamma_i,z=0, \zeta=0)$.}
    \While {$\mathcal{C}_{i,j,k} \neq \emptyset$} \Comment {There are still \glspl{bc} for the current value of $z_j$}
        \State {$j \leftarrow j+1$}
        \State {$z_j \leftarrow z_{j-1}+\Delta z$}
        \State {$\mathcal{S}_{i,j,k=0} \gets \mathcal{C}_{i,j-1,k=0}$} \Comment{The search space $\mathcal{S}_{i,j,k=0}$ is initialized as $\mathcal{C}_{i,j,k=0}$}
        \While {$ \mathcal{C}_{i,j,k=0} \cap \partial \mathcal{S}_{i,j,k=0} \ne \emptyset$}   \Comment {Check for \glspl{bc} on the boundary of $\mathcal{S}_{i,j,k=0}$}
            \State {$\partial \mathcal{S}_{i,j,k=0} \leftarrow \partial \mathcal{S}_{i,j,k=0} +\textrm{offset}(d_O)$}     \Comment{The polygon $\partial \mathcal{S}_{i,j,k=0}$ is enlarged with an offset of size $d_O$}
            \State {Using Eqs.~\eqref{eq: v2_energy_constraint} and \eqref{eq: injection_angle}, compute $(\dot{x},\dot{y},\dot{z})$ on every vertex $(x, y) \in \partial \mathcal{S}_{i,j,k=0}$ returned by \textit{polybuffer}, obtaining $\mathbf{x}_0$.}
            \State {Propagate each initial condition $\mathbf{x}_0$ with \cref{eq: equations of motion} for a time $\tau_{\text{sp}}$.}
            \State{Identify \glspl{bc} according to the definition in \cref{sec: definition BC}.}
            \State{Temporarily add the \glspl{bc} to $\mathcal{C}_{i,j,k=0}$.}
        \EndWhile
        \Comment {The set $\mathcal{S}_{i,j,k=0}$ enclosing $\mathcal{C}_{i,j,k=0}$ ($\mathcal{C}_{i,j,k=0} \subset \mathcal{S}_{i,j,k=0}$) is found} 
        \State{Using Eqs.~\eqref{eq: v2_energy_constraint} and \eqref{eq: injection_angle}, compute $(\dot{x},\dot{y},\dot{z})\; \forall\ (x, y) \in (\mathcal{G} \cap \mathcal{S}_{i,j,k=0})$, obtaining $\mathbf{x}_0$.}
        \State{Verify that $\mathbf{x}_0$ meets condition $\dot{\varepsilon_2}(\tau_0)<0$ with \cref{eq: check derivative eps2}, otherwise discard the current $\mathbf{x}_0$.}
        \State{Using \cref{eq: equations of motion}, propagate each initial condition $\mathbf{x}_0$ backward for a time $\tau_B$ and possibly forward for a time $\tau_{s}$.} 
        \State{Identify \glspl{bc} and add them to $\mathcal{C}_{i,j,k=0}$ according to the definition in \cref{sec: definition BC}.}
    \EndWhile
\EndWhile
\end{algorithmic}
\end{algorithm}

\begin{algorithm}[tbp]
\caption{\gls{etd} Ballistic Capture Set Computation: $\zeta$-sections}
\label{alg: search algorithm zdot}
\begin{algorithmic}[1]
\State {Build a grid $\mathcal{G}$ around $M_2$ with width $x_{\mathcal{G}}$, height $y_{\mathcal{G}}$, and  stepsize $h$. Initialize $\Gamma\gets0$, $i\gets0$, $z\gets0$, $j\gets0$.}
\While {$\Gamma < \Gamma_{max}$} 
    \State {$i \leftarrow i+2$}
    \State {$\Gamma_i \leftarrow \Gamma_{i-2}+2\Delta \Gamma$} \Comment {Step increased to $2\cdot\Delta \Gamma$ to limit the computational time required}
    \While {$\mathcal{C}_{i,j,k=0} \neq \emptyset$} \Comment {There are still \glspl{bc} for the current value of $z_j$}
        \State {$j \leftarrow j+2$}
        \State {$z_j \leftarrow z_{j-1}+2\Delta z$} \Comment {Step increased to $2\cdot\Delta z$ to limit the computational time required}
        \State {Set $k\gets0$, and load the capture set for the current $j-$th $z-$section $\mathcal{C}_{i,j,k} = \mathcal{C}(\Gamma_i,z_j, \zeta=0)$.}
        \While {$\mathcal{C}_{i,j,k} \neq \emptyset$} \Comment {There are still \glspl{bc} for the current value of $\zeta_k$}
            \State {$k \leftarrow k+1$}
            \State {$\zeta_k \leftarrow \zeta_{k-1}+\Delta \zeta$}
            \State {An additional step is included here to consider also negative values of $\zeta_k$ (i.e. negative $\dot{z}_k$).} \label{linealg: symm}
            \State {$\mathcal{S}_{i,j,k} \gets \mathcal{C}_{i,j,k-1}$} \Comment{The search space $\mathcal{S}_{i,j,k}$ is initialized as $\mathcal{C}_{i,j,k-1}$}
            \While {$ \mathcal{C}_{i,j,k} \cap \partial \mathcal{S}_{i,j,k} \ne \emptyset$}   \Comment {Check for \glspl{bc} on the boundary of $\mathcal{S}_{i,j,k}$}
                \State {$\partial \mathcal{S}_{i,j,k} \leftarrow \partial \mathcal{S}_{i,j,k} +\textrm{offset}(d_O)$}     \Comment{The polygon $\partial \mathcal{S}_{i,j,k}$ is enlarged with an offset of size $d_O$}
                \State {Using Eqs.~\eqref{eq: v2_energy_constraint} and \eqref{eq: injection_angle}, compute $(\dot{x},\dot{y},\dot{z})$ on every vertex $(x, y) \in \partial \mathcal{S}_{i,j,k}$ returned by \textit{polybuffer}, obtaining $\mathbf{x}_0$.}
                \State {Propagate each initial condition $\mathbf{x}_0$ with  \cref{eq: equations of motion} for a time $\tau_{\text{sp}}$.}
                \State{Identify \glspl{bc} according to the definition in \cref{sec: definition BC}.}
                \State{Temporarily add the \glspl{bc} to $\mathcal{C}_{i,j,k}$.}
            \EndWhile
            \Comment {The set $\mathcal{S}_{i,j,k}$ enclosing $\mathcal{C}_{i,j,k}$ ($\mathcal{C}_{i,j,k} \subset \mathcal{S}_{i,j,k}$) is found} 
            \State{Using Eqs.~\eqref{eq: v2_energy_constraint} and \eqref{eq: injection_angle}, compute $(\dot{x},\dot{y},\dot{z})\; \forall\ (x, y) \in (\mathcal{G} \cap \mathcal{S}_{i,j,k})$, obtaining $\mathbf{x}_0$.}
            \State{Verify that $\mathbf{x}_0$ meets condition $\dot{\varepsilon_2}(\tau_0)<0$ with \cref{eq: check derivative eps2}, otherwise discard the current $\mathbf{x}_0$.}
            \State{Using \cref{eq: equations of motion}, propagate each initial condition $\mathbf{x}_0$ backward for a time $\tau_B$ and possibly forward for a time $\tau_{s}$.}
            \State{Identify \glspl{bc} and add them to $\mathcal{C}_{i,j,k}$ according to the definition in \cref{sec: definition BC}.}
        \EndWhile
    \EndWhile
\EndWhile
\end{algorithmic}
\end{algorithm}

\subsection{Symmetries of the problem} \label{sec: symmetries in results}
In the \gls{cr3bp}, the dynamics is completely symmetric with respect to the $x$–$y$ plane. This means that two trajectories with symmetric initial conditions $(x,y,z,\dot{x},\dot{y},\dot{z})$ and $(x,y,-z,\dot{x},\dot{y},-\dot{z})$ will generate symmetric trajectories when propagated.
Therefore, the symmetry reflects only the position and velocity along the third component $(z, \dot{z})$ (both of them at the same time), and affects the orbital elements describing the orientation of the orbital plane. As a consequence, also the orientation of the orbit itself in this plane is affected. In fact, these elements are the \gls{raan} $\Omega$ and the argument of periapsis $\omega$.
To obtain the orbital elements $\Omega$ and $\omega$ of a symmetric orbit, both of them must be rotated of $\pi$. All the other orbital elements $a$, $e$, $i$, $\nu$ remain unchanged.

Therefore, \cref{alg: search algorithm z} is implemented only for positive values of $z$. All the resulting \glspl{bc} are duplicated exploiting the symmetry to take into account initial conditions with negative $z$ values. Hence, the new parameters are given by $z \leftarrow-z$, $\Omega\leftarrow\Omega+\pi$, and $\omega\leftarrow\omega-\pi$. For what concerns \cref{alg: search algorithm zdot}, as stated in line~\ref{linealg: symm}, the search must be carried out also for negative values of $\dot{z}$ for each $z$-section. Then, the symmetry can be exploited as mentioned above. In this way, all the \glspl{bc} with $z<0$ and their parameters can be also obtained just by duplicating and adapting the \glspl{bc} previously computed. The new parameters are given by $z \leftarrow-z$, $\dot{z} \leftarrow-\dot{z}$, $\Omega\leftarrow\Omega+\pi$, and $\omega\leftarrow\omega-\pi$.

\subsection{Sample spatial capture sets} \label{sec: sample spatial capture sets}
Example of the results given by Algorithms~\ref{alg: search algorithm z} and \ref{alg: search algorithm zdot} are provided in this section for the Earth-Moon system. \cref{tab: param algs} summarizes the parameters used to obtain the following results. Note that, at this stage, a step $\Delta z = 10h$ is adopted to limit the computational time required. For the same reason, \cref{alg: search algorithm zdot} employs steps of $2\cdot\Delta \Gamma$ and $2\cdot\Delta z$. As a consequence, only some capture sets $\mathcal{C}_{i,j,k} = \mathcal{C}(\Gamma_i,z_j, \zeta=0)$ are used to initialize the search detailed in \cref{alg: search algorithm zdot}.

\begin{table}[tbp]
\caption{\label{tab: param algs} Parameters used in Algorithms~\ref{alg: search algorithm z} and \ref{alg: search algorithm zdot}}
\centering{}
\begin{tabular}{ccccccccc}
\hline
\noalign{\vskip\doublerulesep}
$x_{\mathcal{G}}$ [LU] & $y_{\mathcal{G}}$ [LU] & $h$ [LU] & $d_O$ [LU] & $\Delta \Gamma$ & $\tau_{s}$ [TU] & $\tau_{\text{sp}}=\tau_B$ [TU] & $\Delta z$ [LU] & $\Delta \dot{z}$ [LU] \tabularnewline[\doublerulesep]
\hline
\noalign{\vskip\doublerulesep}
\noalign{\vskip\doublerulesep}
$7 \, r_H$ & $9 \, r_H$ & $4\cdot10^{-4}$ & $ 2\cdot10^{-3}$ & $0.02$  & $10 \cdot \left(2\pi\right)$ & $2 \cdot \left(2\pi\right)$ & $4\cdot10^{-3}$ & 1 deg \tabularnewline[\doublerulesep]
\noalign{\vskip\doublerulesep}
\end{tabular}
\end{table}

\begin{figure}[tbp]
    \centering
    \includegraphics[width=\textwidth]{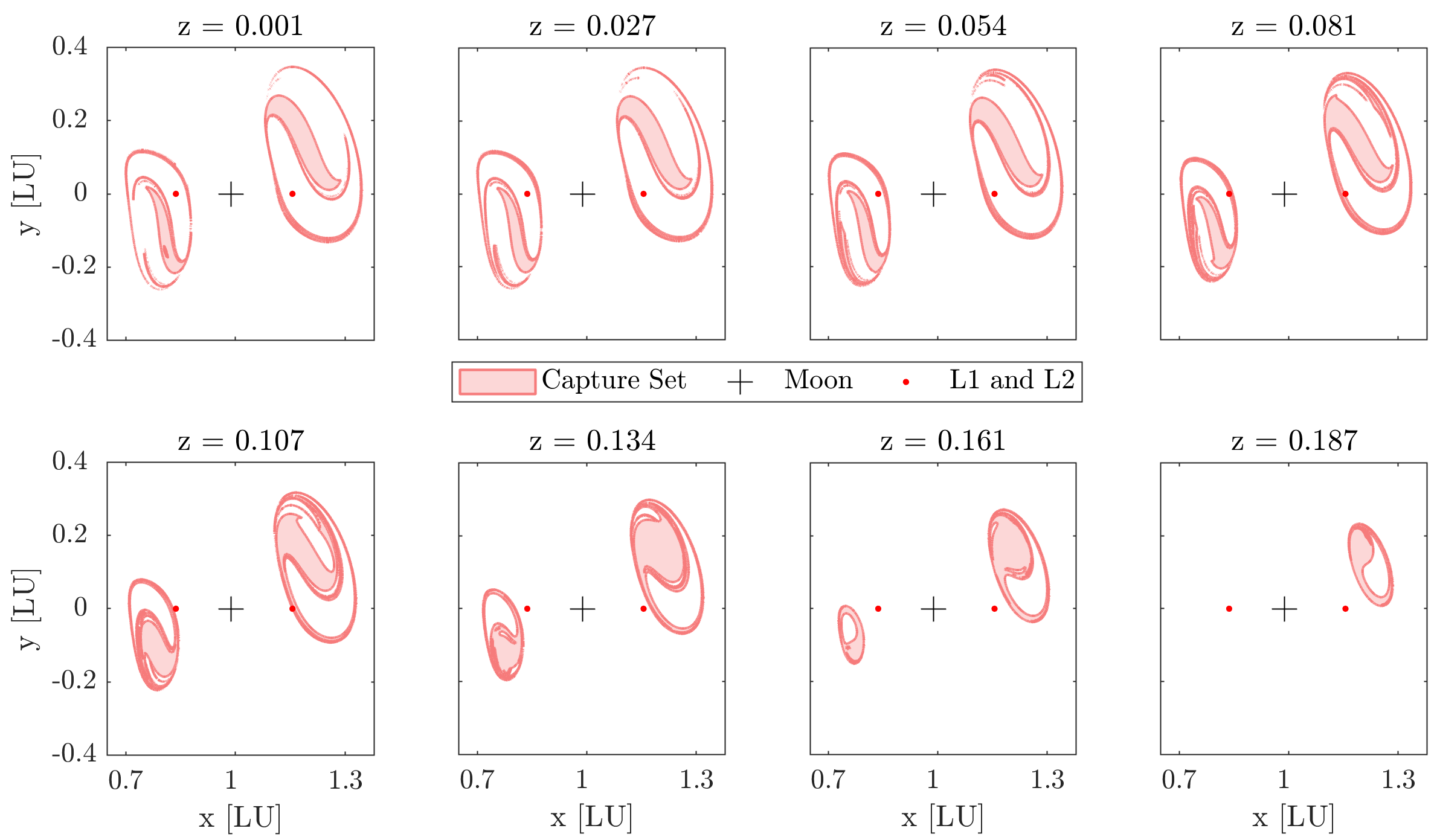}
    \caption{Results for \cref{alg: search algorithm z} when $\Gamma=0.52$, i.e. $\mathcal{C}(\Gamma=0.52, z, \zeta=0)$.}
    \label{fig: capture set z Gamma52}
\end{figure}
\begin{figure}[tbp]
    \centering
    \includegraphics[width=\textwidth]{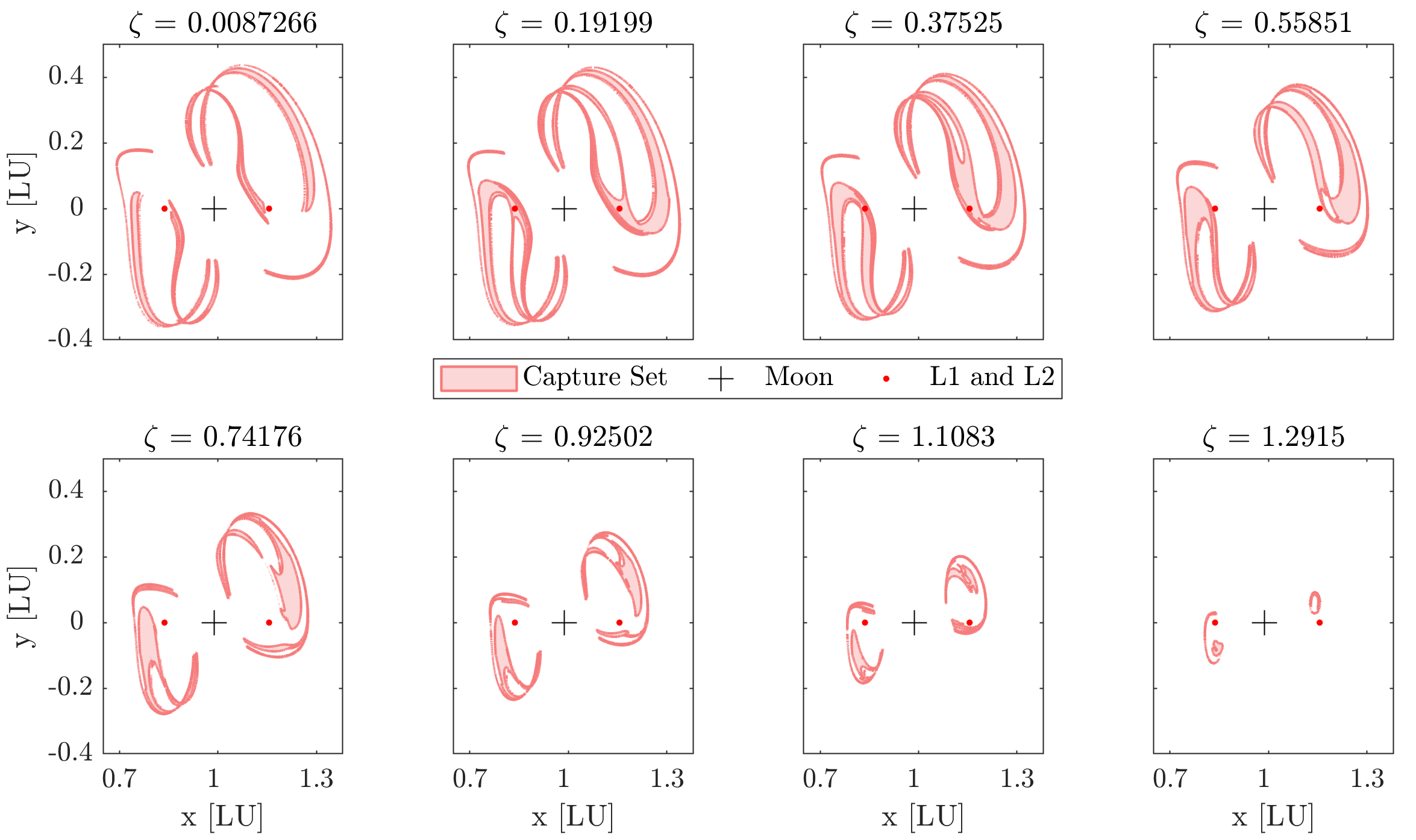}
    \caption{Results for \cref{alg: search algorithm zdot} when $\Gamma=0.88$ and $z=0$, i.e. $\mathcal{C}(\Gamma=0.88, z=0, \zeta)$.}
    \label{fig: capture set zdot Gamma88}
\end{figure}

In \cref{fig: capture set z Gamma52} a sample result of \cref{alg: search algorithm z} is presented, where only 8 $z-$sections of the capture sets are displayed for $\Gamma=0.52$, i.e. $\mathcal{C}(\Gamma=0.52, z, \zeta=0)$. 
Instead, \cref{fig: capture set zdot Gamma88} shows a sample result of \cref{alg: search algorithm zdot}, where only 8 $\zeta-$sections of the capture set are displayed for $\Gamma=0.88$ and $z=0$, i.e. $\mathcal{C}(\Gamma=0.88, z=0, \zeta)$.

\subsection{Database structure}

For each \gls{bc} belonging to a $\mathcal{C}_{i,j,k}$, a wide number of parameters defining its main features is computed and stored to form a database. A first group of parameters describes the initial conditions of the \glspl{bc} (which belong to the \gls{etd}). The position coordinates $(x,y,z)$, the three-body energy parameter $\Gamma$, and the velocity angles $\sigma$ and $\zeta$ are saved. Note that the injection angle $\sigma$ is only saved to distinguish between the two possible solutions of \cref{eq: injection_angle}.
A second group of parameters describes the robustness of the capture, with parameters such as the total number of revolutions, the number of prograde-only or retrograde-only revolutions (with sign), the time spent in the capture phase, the number of crossings of the $\varepsilon_2=0$ condition, and the time before a possible collision occurs. A third group of parameters defines the origin of a \gls{bc} as a two-body orbit around the Earth. In fact, for each \gls{bc}, the Cartesian coordinates at the escape time $\tau_e<\tau_0$ are transformed into orbital elements of an Earth-centered inertial frame ($a_T$, $e_T$, $i_T$, $\Omega_T$, $\omega_T$, $\nu_T$) and stored. The last and largest group - similarly to the previous one - stores the orbital elements of $M_3$ at selected periapsis around $M_2$. For the \gls{cr3bp} problem, the first periapsis and the closest periapsis (in terms of distance to $M_2$) orbital elements are both saved.
Note that the orbital elements for the \gls{bc} origin are extracted in a moment - the escape time $\tau_e$ - when the perturbation of $M_2$ is of minor importance, and the origin of the \gls{bc} can be well approximated by a two-body orbit around the Earth. Similarly, the orbital elements of all the periapsis are obtained when the gravitational pull of $M_2$ is largest and the perturbation of the Earth can be considered to be at its minimum.

\section{Transition to an ephemeris model} \label{sec: transition to ephemeris}

While the \gls{cr3bp} provides a simplified yet insightful model for three-body dynamics, higher-fidelity trajectory design requires accounting for additional gravitational perturbations. To this end, we refine promising trajectories identified within the \gls{cr3bp} by transforming their initial conditions to an ephemeris-based model that incorporates gravitational forces from multiple celestial bodies.
In this section, a method to transform the initial conditions from the \gls{cr3bp} to an ephemeris model is introduced. Since these conditions belong to the \gls{etd}, only the transformation of a single initial condition at a fixed epoch is required for each selected trajectory. Once the initial conditions are transformed, the trajectories are repropagated entirely within an ephemeris model (based on the Spice toolkit~\cite{spice}) which is also introduced in this section.
Finally, a distance metric is introduced to shrink the amount of \glspl{bc} that need to be transitioned to the ephemeris model.

\subsection{Transformation of initial conditions (ICs) into an ephemeris model} \label{sec: transformation of ICs}

Following previous work from the literature~\cite{DeiTos_ephemerisModel, ParkHowell_ephemerisModel}, a method for transforming initial conditions from the Earth-Moon \gls{cr3bp} into an ephemeris problem is presented here. For the current work, the synodic frame of the \gls{cr3bp} coincides with the rotopulsating frame introduced in~\cite{ParkHowell_ephemerisModel} (which is transitioned into ephemeris).
A key feature of this process is that only the initial conditions of selected trajectories (which belong to the \gls{etd}) are transformed at a fixed epoch, rather than the entire trajectory. The transformed initial conditions are then propagated within the ephemeris model and serve as starting points for refining the database into a higher-fidelity model, as detailed in the next section.

The first step is to identify the dimensionless time $\tau$ of the \gls{cr3bp} model and the dimensional time $T$ in seconds of the ephemeris model.
Different notations are adopted to distinguish between different time derivatives. The differentiation with respect to the dimensionless time is denoted using a dot, as for the quantities introduced in the \gls{cr3bp}: $\dot{\mathbf{r}}=\mathbf{v}$. Conversely, in an inertial frame centered on Earth, the position vector is denoted as $\mathbf{R}$, and its time derivative with respect to dimensional time is represented as $\mathbf{R}^{\prime}=d\mathbf{R}/dT$.

To transform an initial condition from the synodic frame of the \gls{cr3bp}, defined as $\mathbf{x} = (\mathbf{r}, \mathbf{v})$, into an initial condition in the Earth inertial frame, given by $(\mathbf{R}, \mathbf{R}^{\prime})$, the following relations are introduced~\cite{ParkHowell_ephemerisModel}:
\begin{equation}
\left\{ \begin{array}{l}
\mathbf{R}=\mathbf{B}+\ell\mathbf{C}\mathbf{r}
\\
\mathbf{R}^{\prime}= \mathbf{B}^{\prime}+(\ell^{\prime}\mathbf{C}+\ell\mathbf{C}^{\prime})\mathbf{r} +\ell\mathbf{C}\tau^{\prime}\mathbf{v}
\\
\end{array}\right. \; \; ,
\label{eq: rotopuls to inertial}
\end{equation}
where $\mathbf{B}$ is the position of the system barycentre in the new frame, $\ell$ is the dimensionalization factor which represents the pulsation of the Earth-Moon distance, and $\mathbf{C}$ is a direction cosine matrix. Geometrically speaking, $\mathbf{B}$ represents a translation, while the following terms represent a scaling and a rotation which aligns the synodic frame with the inertial frame. In addition, $\tau^{\prime}=d\tau/dT$ is an arbitrary function that describes the isochronous correspondence between the models and it is here defined in a non-uniform fashion as introduced by Park et al.~\cite{ParkHowell_ephemerisModel}: 
\begin{equation}
\tau^{\prime}=\sqrt{\cfrac{G(m_E+m_M)}{\ell^3}} \, .
\label{eq: t prime}
\end{equation}
Although approximate, this transformation ensures a close match between \gls{cr3bp} and ephemeris initial conditions, providing a smooth transition.

The algorithmic procedure begins by computing the Moon's position $\mathbf{r}_M$ and velocity $\mathbf{v}_M$ at a selected epoch from the ephemeris function. 
The norm $\ell=\left\lVert\mathbf{r}_M\right\lVert$ is the pulsation factor for the selected epoch, and its derivative over time can be computed as $\ell^{\prime}=\mathbf{r}_M \cdot \mathbf{v}_M/\ell$. As introduced by Dei Tos et al.~\cite{DeiTos_ephemerisModel}, the direction cosine matrix and its derivative can be defined as $\mathbf{C}=[\mathbf{e}_1, \mathbf{e}_2, \mathbf{e}_3]$ and $\mathbf{C}^{\prime}=[\mathbf{e}^{\prime}_1, \mathbf{e}^{\prime}_2, \mathbf{e}^{\prime}_3]$, respectively, where $\mathbf{e}$ are column vectors
\begin{equation}
\left\{ \begin{array}{l}
\mathbf{e}_1 = \mathbf{r}_M / \ell
\\
\mathbf{e}_2 = \mathbf{e}_3 \times \mathbf{e}_1
\\
\mathbf{e}_3 = \mathbf{r}_M \times \mathbf{v}_M/h
\end{array}\right. \, ,
\label{eq: C components}
\end{equation}
\begin{equation}
\left\{ \begin{array}{l}
\mathbf{e}^{\prime}_1 = \cfrac{\ell \mathbf{v}_M - \ell^{\prime} \mathbf{r}_M}{\ell^2}
\\
\mathbf{e}^{\prime}_2 = \mathbf{e}^{\prime}_3 \times \mathbf{e}_1 + \mathbf{e}_3 \times \mathbf{e}^{\prime}_1
\\
\mathbf{e}^{\prime}_3 = \cfrac{h \left(\mathbf{r}_M \times \mathbf{a}_M\right) - h^{\prime} \left(\mathbf{r}_M \times \mathbf{v}_M\right)}{h^2}
\end{array}\right. \, ,
\label{eq: C' components}
\end{equation}
and where the norm of the angular momentum of the Moon is
\begin{equation}
    h=\left\lVert\mathbf{r}_M\times\mathbf{v}_M\right\lVert \, .
\end{equation}

In this work, $\mathbf{a}_M$, $\mathbf{h}^{\prime}$ and its norm $h^{\prime}$ are obtained using a finite difference estimation from $\mathbf{v}_M$ and $\mathbf{h}=\mathbf{r}_M\times\mathbf{v}_M$. 

Given this, \cref{eq: rotopuls to inertial} can be used to compute the new initial conditions in the ephemeris model for the selected epoch.

\subsection{Ephemeris model} \label{sec: ephemeris model}
The ephemeris-based model implemented in this work accounts for the gravitational influences of Earth, Moon, and Sun using the Spice toolkit~\cite{spice} with the DE430 ephemerides. The positions of the celestial bodies are first retrieved in an Earth-centered inertial frame (e.g., EMO2000) and then rotated into an \gls{etd}-aligned Earth-inertial frame with coordinates $(x_1, y_1, z_1)$. This frame is also centered on Earth, with its $x_1$-axis is aligned with the Earth–Moon direction at a specific epoch, $T_{\text{ETD}}$, which in this work corresponds to the nominal insertion scenario described in the next section.
This choice mirrors the convention $\tau_0 = 0$ used in \cref{eq: coord wrt M2} and throughout \cref{sec: ETD} for the \gls{cr3bp} formulation, allowing a consistent alignment of frames at the epoch of interest. Additionally, the frame is defined such that the Moon’s osculating orbital plane at $T_{\text{ETD}}$ lies in the $x_1$–$y_1$ plane. Similarly, a Moon-centered inertial frame with coordinates $(x_2, y_2, z_2)$ is defined following the same convention. As a result, the synodic frame and both Earth- and Moon-centered inertial frames are all aligned at $T_{\text{ETD}}$, ensuring consistency across all dynamical models used.
Only these three frames are employed throughout the paper.

The model takes into account the positions of the satellite $\mathbf{R}$ and the ones of the aforementioned celestial bodies $\mathbf{R}_j$ to compute a dimensional acceleration $\mathbf{R}^{\prime\prime}$. Following Newton's law of universal gravitation~\cite{valladoNbody}, in an $N$-body framework, the equations of motion are given by:
\begin{equation}
    \mathbf{R}^{\prime\prime} = -\cfrac{G(m_1+m_{sat})}{R^3} \mathbf{R} - G\sum\limits_{j=2}^{N-1} m_j \left( \cfrac{\mathbf{R}-\mathbf{R}_{j}}{(R-R_{j})^3} + \cfrac{\mathbf{R}_{j}}{R_{j}^3} \right) \, .
    \label{eq: Newton law generic}
\end{equation}
where we set $N=4$, and subscripts $1$, $2$, and $3$ refer to Earth, Moon, and Sun, respectively. In addition, $G$ is the universal gravitational constant, and the contribution of $m_{sat}$ is negligible, consistently with \gls{cr3bp} assumptions.
Substituting the subscripts $1$, $2$, $3$ with $E$, $M$, $S$, the equations of motion become:
\begin{equation}
    \mathbf{R}^{\prime\prime} = -\cfrac{G m_E}{R^3} \mathbf{R} - G m_M \cfrac{\mathbf{R}-\mathbf{R}_M}{(R-R_M)^3} - G m_M \cfrac{\mathbf{R}_M}{R_M^3} - G m_S \cfrac{\mathbf{R}-\mathbf{R}_S}{(R-R_S)^3} - G m_S \cfrac{\mathbf{R}_S}{R_S^3} \, .
    \label{eq: Newton law EMS}
\end{equation}
Using the aforementioned Spice kernels, it is possible to retrieve the vectors $\mathbf{R}_M$ and $\mathbf{R}_S$ in the chosen Earth-inertial frame.

\subsection{Distance metric for the preliminary selection of BCs} \label{sec: distance metric}
To apply promising \glspl{bc} to an already known scenario, a method to select relevant trajectories is introduced.
The problem of identifying \glspl{bc} that are most compatible with a specific mission can be addressed using a distance metric that captures the similarity between two trajectories based on their orbital elements.

As previously mentioned, for each \gls{bc} belonging to the capture set, the orbital elements of the \gls{bc} origin (in the \gls{etd}-aligned Earth-inertial frame) are stored in a database. They are referred to as $a_T$, $e_T$, $i_T$, $\Omega_T$, and $\omega_T$, where the subscript stands for ``target'', as they will be the target orbital elements to obtain a ballistic insertion at the Moon. 
To compare the \glspl{bc} with a certain trajectory, the orbital elements of the latter ($a_I$, $e_I$, $i_I$, $\Omega_I$, and $\omega_I$) are extracted at the escape condition in the same \gls{etd}-aligned Earth-inertial frame. This is done coherently to what introduced in \cref{sec: definition BC}, i.e., when the dimensionless distance from the Moon is $0.9$ LU.

The distance metric $d_v$ provides a first-order estimate of the impulsive correction needed to transfer from a given trajectory to a \gls{bc}, following Wakker~\cite{Wakker_Gauss4manouvers}. It assumes a mono-impulsive maneuver and computes sub-costs for each orbital element individually while holding the others fixed at $a_I$, $e_I$, $i_I$, $\Omega_I$, and $\omega_I$, as summarized in \cref{tab: distance metric eq summary}. These sub-costs rely on simplifying assumptions that may differ between elements (e.g., maneuvers at perigee vs. apogee). The final metric is obtained as
$d_v = \sqrt{\Delta v_a^2+\Delta v_e^2+\Delta v_i^2+\Delta v_{\Omega}^2+\Delta v_{\omega}^2}$.
The true anomaly $\nu$ is not explicitly considered, and phasing costs are assumed to align with semimajor axis changes. While approximate, the metric effectively captures the main trends in transfer cost, as illustrated in the final results.

Please note that the equation used for $\Delta v_{\omega}$ is not the primary one recommended by Wakker~\cite{Wakker_Gauss4manouvers}, but rather an alternative expression from the same reference, selected because $i_I \xrightarrow{} 0$ in the case discussed in the next section.

Despite its lack of precision, this method allows to effectively obtain the most promising \glspl{bc} suitable for an applicative scenario. Validation cases will be provided in the result section, proving its applicability in this context. Using this distance metric, it can be guaranteed that all the \glspl{bc} with their origin compatible with a required orbit are picked just by selecting all the ones below an inflated threshold, $d_v^{\text{thr}}$. Using the criteria $d_v\leq d_v^{\text{thr}}$, a subset of the entire database $\mathcal{C}_{dv}(\Gamma,z,\zeta) \subset \mathcal{C}(\Gamma,z,\zeta)$ can be extracted.

\begin{table}[tbp]
\caption{Summary of the distance metric equations for the different orbital elements}
\centering{}\label{tab: distance metric eq summary}
\begin{tabular}{ccc}
\hline
\noalign{\vskip\doublerulesep}
Orbital element change & Distance metric equation & Assumptions (\cref{tab: assumptions}) \tabularnewline[\doublerulesep]
\hline
\noalign{\vskip\doublerulesep}
\noalign{\vskip\doublerulesep}
$\Delta a=a_T-a_I$ & $\Delta v_a=\dfrac{1}{2}\dfrac{\Delta a}{a}\sqrt{\dfrac{Gm_E(1-e)}{a(1+e)}}$ & $A \land B$ \tabularnewline[\doublerulesep]
\noalign{\vskip\doublerulesep}
\noalign{\vskip\doublerulesep}
$\Delta e=e_T-e_I$ & $\Delta v_e=\dfrac{1}{2}\Delta e\sqrt{\dfrac{Gm_E}{a(1-e^2)}}$ & $A \land (B \lor C)$ \tabularnewline[\doublerulesep]
\noalign{\vskip\doublerulesep}
\noalign{\vskip\doublerulesep}
$\Delta i=i_T-i_I$ & $\Delta v_i=\sqrt{\dfrac{Gm_E(1-e)}{a(1+e)}}\Delta i$ 
& $A \land C \land D$ \tabularnewline[\doublerulesep]
\noalign{\vskip\doublerulesep}
\noalign{\vskip\doublerulesep}
$\Delta \Omega=\Omega_T-\Omega_I$ & $\Delta v_{\Omega}=\sqrt{\dfrac{Gm_E(1-e)}{a(1+e)}} \sin{i} \, \Delta \Omega$ & $A \land C \land E$ \tabularnewline[\doublerulesep]
\noalign{\vskip\doublerulesep}
\noalign{\vskip\doublerulesep}
$\Delta \omega=\omega_T-\omega_I$ & $\Delta v_{\omega}=\dfrac{e}{2}\sqrt{\dfrac{Gm_E}{a(1-e^2)}} \Delta \omega$ & $A \land C \land E$ \tabularnewline[\doublerulesep]
\noalign{\vskip\doublerulesep}
\end{tabular}
\end{table}

\begin{table}[tbp]
\caption{Assumptions to obtain sub-costs equations in \cref{tab: distance metric eq summary}}
\centering{}\label{tab: assumptions}
\begin{tabular}{cc}
\hline
\noalign{\vskip\doublerulesep}
Assumption & Description \tabularnewline[\doublerulesep]
\hline
\noalign{\vskip\doublerulesep}
\noalign{\vskip\doublerulesep}
A & Small orbital elements change \tabularnewline[\doublerulesep]
\noalign{\vskip\doublerulesep}
\noalign{\vskip\doublerulesep}
B & Maneuver performed at perigee \tabularnewline[\doublerulesep]
\noalign{\vskip\doublerulesep}
\noalign{\vskip\doublerulesep}
C & Maneuver performed at apogee \tabularnewline[\doublerulesep]
\noalign{\vskip\doublerulesep}
\noalign{\vskip\doublerulesep}
D & \begin{tabular}{@{}c@{}} Maneuver performed at ascending \\ or descending node \end{tabular} \tabularnewline[\doublerulesep]
\noalign{\vskip\doublerulesep}
\noalign{\vskip\doublerulesep}
E & \begin{tabular}{@{}c@{}} Maneuver performed $90$ deg after \\ the ascending or descending node \end{tabular} \tabularnewline[\doublerulesep]
\noalign{\vskip\doublerulesep}
\end{tabular}
\end{table}

Given the symmetric properties introduced at the end of \cref{sec: planar2spatial4symm}, the metric must be applied also to symmetric conditions. As previously introduced, the \gls{raan} $\Omega_T$ and the argument of periapsis $\omega_T$ must be rotated of $\pi$. All the other orbital elements $a_T$, $e_T$, $i_T$, $\nu_T$ must remain untouched.

\section{Tailoring for mission-specific requirements} \label{sec: test case}
This section presents a case study based on a real space mission. The same methodology can be extended to other missions that utilize low-energy transfers to reach the Moon. However, this process is applicable only in advanced mission design phases, where an optimized nominal trajectory is already available. 

Key aspects of the \gls{ltb} mission and its trajectory are introduced. To preliminarily identify suitable \glspl{bc}, the distance metric with respect to \gls{ltb} reference mission trajectory - introduced in \cref{sec: distance metric} - is applied. In this was a significant reduction in the number of candidate \glspl{bc} requiring transition into an ephemeris model is obtained, thereby lowering computational costs. 

\subsection{Lunar Trailblazer}
\gls{ltb} was a Small Innovative Mission for Planetary Exploration (SIMPLEx) mission that studied the form, abundance and distribution of water on the Moon, as well as the lunar water cycle\footnote{\url{https://www.jpl.nasa.gov/missions/lunar-trailblazer/}}. To succeed in its objectives, the spacecraft was to be inserted into a scientific orbit around the Moon, which was polar and had an altitude of 100 km over the surface.
Moreover, being a low-budget mission, the spacecraft had very limited maneuver capabilities, most of all regarding thrust. This meant that impulsive burns greater than 60 m/s could not  be executed in a short time span near perilune. For this reason, a low-energy trajectory design is needed, making the insertion maneuver design especially critical.

\begin{figure}[tbp]
    \centering
    \begin{subfigure}[t!]{0.49\textwidth}
        \includegraphics[width=\textwidth]{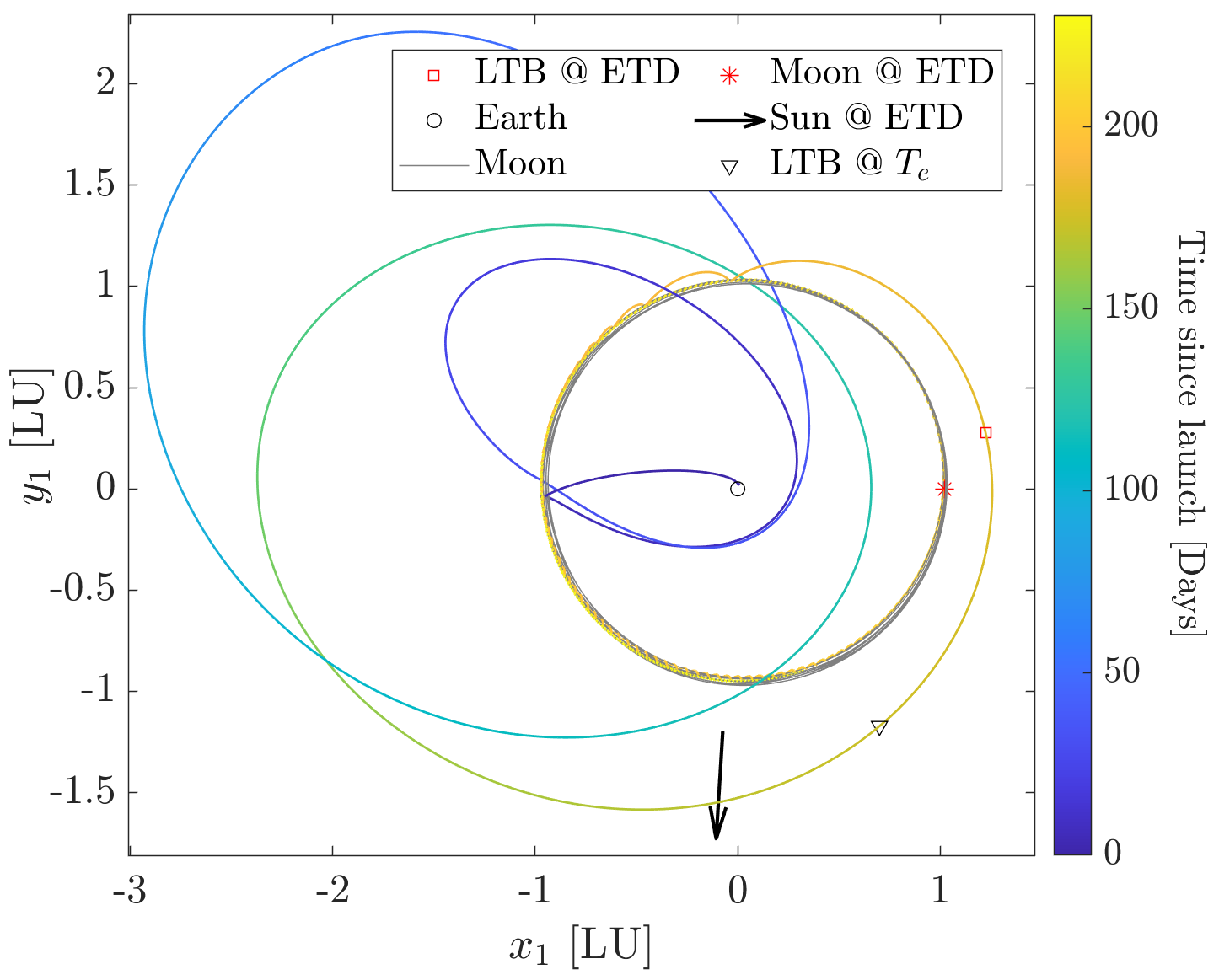}
    \end{subfigure}\\
    \begin{subfigure}[t!]{0.49\textwidth}
        \includegraphics[width=\textwidth]{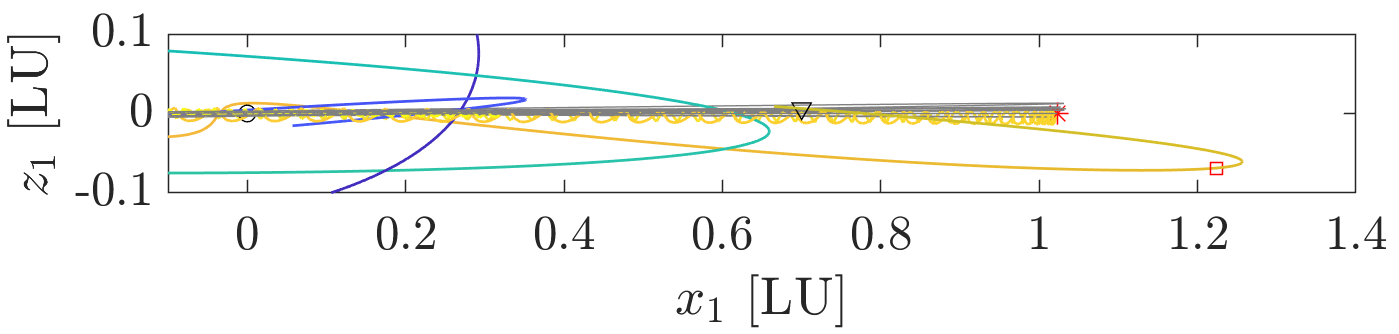}
    \end{subfigure}
    \caption{Nominal \gls{ltb} trajectory for departure on 4th Dec 2024 in the Earth-centered inertial frame. Colorbar indicates time after Day 0 (i.e. launch date).}
    \label{fig: LTB nominal planar and 3D}
\end{figure}

The nominal trajectory of \gls{ltb} used as a reference for the present work is represented in \cref{fig: LTB nominal planar and 3D} within the Earth-centered inertial frame introduced in \cref{sec: ephemeris model}. Dimensionless units are used, as defined in \cref{tab: scaling units}. Views along both the $x_1$–$y_1$ and $x_1$–$z_1$ planes are displayed to provide insight into the 3D dynamics. In this frame, the $x_1$-axis points toward the Moon at time $T_{\text{ETD}}$, the epoch when \gls{ltb} is at the \gls{etd} condition (Day 181, see Appendix~\ref{appendix: sample ICs} for the precise value). Time in days after the launch date (Day 0) is represented by the colorbar. The design considered here is based on an outdated trajectory with a launch date of December 4, 2024, which is referred to as Day 0 in the following. Although this was not the actual trajectory used for the mission, its overall structure remains consistent with the nominal trajectory planned for the February 26, 2025, launch\footnote{\url{https://www.jpl.nasa.gov/news/how-nasas-lunar-trailblazer-will-make-a-looping-voyage-to-the-moon/}}.
The design analyzed in this work presents a launch in an orbit directly pointing at the Moon. Due to this low-budget mission's limited propulsive capabilities, a direct \gls{loi} is not possible, as it would require an impulsive burn of nearly $200$ m/s. Therefore, the trajectory must be corrected to perform a first (almost polar) flyby of the Moon, followed by a second one after one revolution around the Earth. 
This second flyby increases the total energy of the two-body orbit around the Earth, as can be noticed by the high semimajor axis of the following two revolutions. During this time, the Sun perturbation kicks in by mainly decreasing the eccentricity value and guaranteeing a much lower excess velocity with respect to the Moon $v_\infty$ during the final approach. A correction maneuver of approximately 10 m/s is also implemented on Day 160 
to correct the spacecraft's trajectory and optimize the insertion at the Moon. 
The contribution of the Sun is so important that it gets to a $v_\infty<0$. Consequently, an insertion into a very short \gls{bc} is achieved.
In other words, the two-body energy with respect to the Moon $\varepsilon_2$ during the insertion becomes lower than zero (locally, it becomes an elliptic orbit around the Moon), and this allows for a very small nominal \gls{loi} (at the time $T_{LOI}$ on Day 189 
), in the order of 30 m/s. As a consequence, this type of trajectory needs to cross the \gls{etd}. This particular condition occurs at time $T_{\text{ETD}}$ on Day 181 
and it is represented with a red square in \cref{fig: LTB nominal planar and 3D}. The Moon's position at the same epoch is represented with a red star. Finally, the insertion in a (weakly) bounded orbit around the Moon is completed, and successive period-reduction maneuvers to reach science orbit are implemented.
It must be noted that the discussion above refers to \gls{ltb} nominal trajectory, which exists in a real ephemeris model. This leads to a significant discrepancy between the concept of \gls{etd} in the \gls{cr3bp} and the \gls{etd} condition of \gls{ltb} in the ephemeris model represented in \cref{fig: LTB nominal planar and 3D}. Therefore, in the latter case, the \gls{etd} condition is just the $\varepsilon_2=0$ condition, and the \gls{bc} behavior of \gls{ltb} cannot be found in the capture set $\mathcal{C}(\Gamma,z,\zeta)$.

\subsection{Application of the distance metric to the insertion phase of LTB}

\begin{figure}[tbp]
    \centering
    \includegraphics[width=0.49\textwidth]{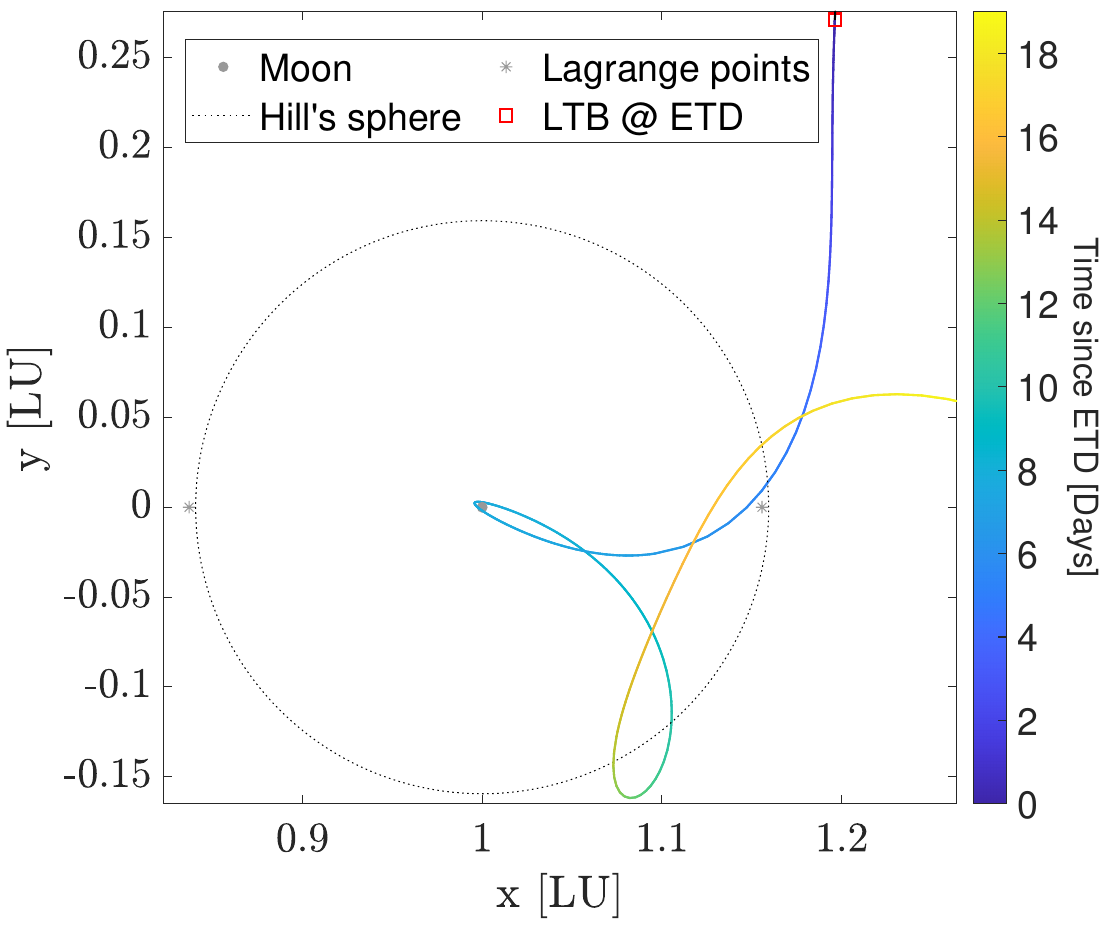}
    \caption{Close-up around the Moon in the synodic frame, without an applied burn at \gls{loi}. Colorbar indicates time after \gls{etd} state (red square).}
    \label{fig: LTB no LOI}
\end{figure}

In this work, the focus is only on the insertion phase of \gls{ltb}, which takes place approximately between Day 170 and Day 200. 
While this trajectory is designed for optimal transfer, exploring ballistic capture alternatives can provide additional flexibility, particularly in contingency scenarios.
The fundamental idea is to identify alternative insertion orbits that exploit \gls{bc} and can be reached from the \gls{ltb} nominal trajectory with minimal correction maneuvers. 
This way, an insertion can be designed and customized with the required mission constraints. The primary goal is to ensure the insertion remains robust against potential thrust execution errors during the \gls{loi}. This is not the case for the \gls{ltb} nominal trajectory (see \cref{fig: LTB no LOI} in the synodic frame), where a thrust failure scenario leads to an immediate escape, resulting in the loss of the spacecraft. At least another suitable perilune could ensure the robustness for a second insertion chance to follow soon after the first.

The starting condition of the \gls{ltb} insertion phase is fixed at the escape condition (see \cref{sec: definition BC}), which occurs on Day 174 at time $T_e$ and it is represented with a black triangle in \cref{fig: LTB nominal planar and 3D}. 
At this epoch, its osculating orbital elements around Earth $a_I$, $e_I$, $i_I$, $\Omega_I$, and $\omega_I$ are given in \cref{tab: oe LTB}. These elements are fundamental for computing the distance metric (\cref{sec: distance metric}).

\begin{figure}[tbp]
    \centering
    \begin{subfigure}[t!]{0.48\textwidth}
        \includegraphics[width=\textwidth]{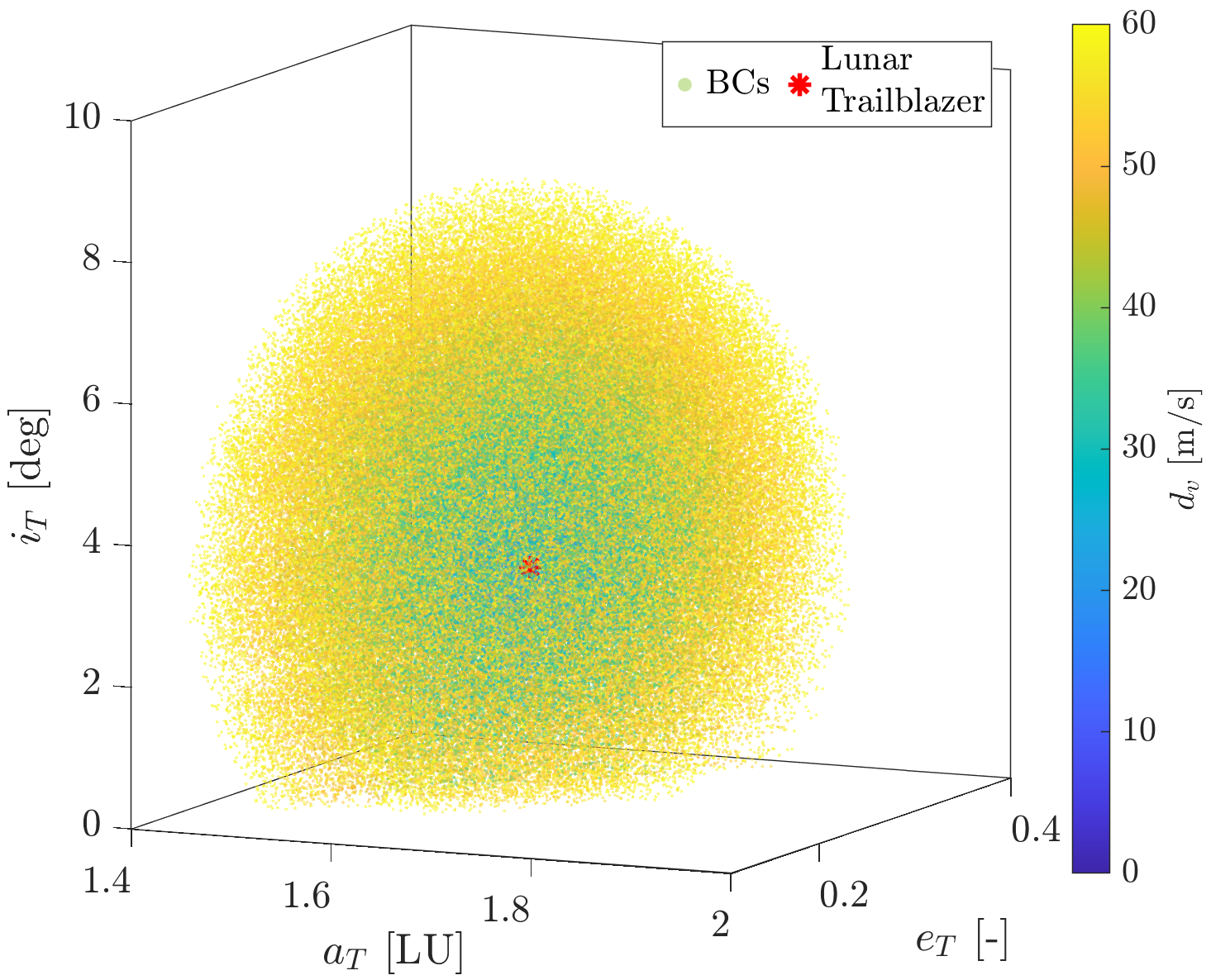}
        \subcaption{Distance metric $\tilde{d}_v$ in the orbital elements $a_T$, $e_T$, $i_T$}
        \label{fig: distance metric CR3BP aei}
    \end{subfigure} \hfill
    \begin{subfigure}[t!]{0.48\textwidth}
        \includegraphics[width=\textwidth]{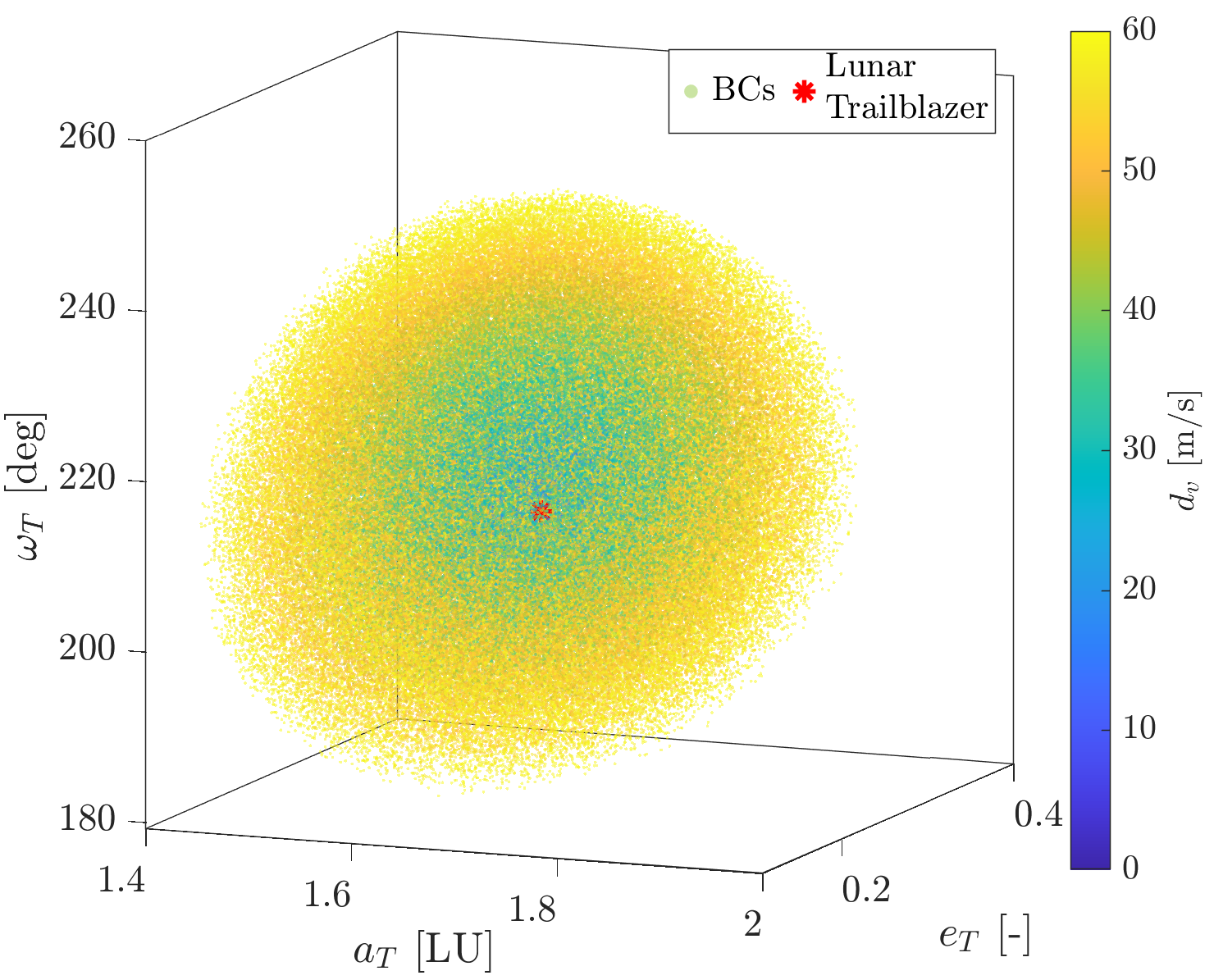}
        \subcaption{Distance metric $\tilde{d}_v$ in the orbital elements $a_T$, $e_T$, $\omega_T$}
        \label{fig: distance metric CR3BP aeomega}
    \end{subfigure}
    \caption{Distance metric in \gls{cr3bp}: $\tilde{d}_v$.}
    \label{fig: distance metric CR3BP}
\end{figure}

To enable a direct comparison with the \glspl{bc}, \gls{ltb} is backward-propagated within the \gls{cr3bp} framework, yielding the corresponding orbital elements at escape: $\tilde{a}_I$, $\tilde{e}_I$, $\tilde{i}_I$, $\tilde{\Omega}_I$, and $\tilde{\omega}_I$, as listed in \cref{tab: oe LTB CR3BP}. 
The primary difference between these elements and those obtained in the full-ephemeris model (\cref{tab: oe LTB}) lies in the semimajor axis and eccentricity, as expected. This deviation is attributed to the influence of the Sun’s gravitational perturbations, which predominantly act in the orbital plane and strongly affect these parameters.
For consistency, the orbital elements computed in the \gls{cr3bp} are first used to compare \gls{ltb} with the \glspl{bc}, since the latter are obtained within this dynamical framework. 
At a later stage, the orbital elements computed in the ephemeris model ($a_I$, $e_I$, $i_I$, $\Omega_I$, and $\omega_I$) will serve as key parameters for selecting the most suitable full-ephemeris \glspl{bc} for \gls{ltb}. The distance metric $\tilde{d}_v$ shown in \cref{fig: distance metric CR3BP} is computed using the \gls{cr3bp} orbital elements of \gls{ltb} $\tilde{a}_I$, $\tilde{e}_I$, $\tilde{i}_I$, $\tilde{\Omega}_I$, and $\tilde{\omega}_I$. Based on this metric, \glspl{bc} are selected using the threshold $\tilde{d}_v<60$ m/s, yielding a capture sub-set $\mathcal{C}_{\tilde{d}v}(\Gamma,z,\zeta) \subset \mathcal{C}(\Gamma,z,\zeta)$. Only $0.1\%$ of the approximately 200 million \glspl{bc}~$\in \mathcal{C}_{\tilde{d}v}(\Gamma,z,\zeta)$ are represented here. Two projections are presented: one in the $(a_T, e_T, i_T)$ space (\cref{fig: distance metric CR3BP aei}) and another in the $(a_T, e_T, \omega_T)$ space (\cref{fig: distance metric CR3BP aeomega}). As expected, inclination $i_T$ plays the most significant role in the metric. In contrast, \gls{raan} $\Omega_T$ has a negligible impact due to the small inclination value ($i_T \approx 3^\circ$), making its representation unnecessary.

\begin{table}[tbp]
\caption{Orbital elements of \gls{ltb} at its escape time $T_e$ in the full-ephemeris model}
\centering{}\label{tab: oe LTB}
\begin{tabular}{cccccc}
\hline
\noalign{\vskip\doublerulesep}
$a_{I}$ [LU] & $e_{I}$ [-] & $i_{I}$ [deg] & $\Omega_{I}$ [deg] & $\omega_{I}$ [deg] & $\nu_{I}$ [deg] \tabularnewline[\doublerulesep]
\hline
\noalign{\vskip\doublerulesep}
\noalign{\vskip\doublerulesep}
$1.8137$ & $0.2915$ & $3.4401$ & $125.2589$ & $217.7110$  & $317.7567$ \tabularnewline[\doublerulesep]
\noalign{\vskip\doublerulesep}
\end{tabular}
\end{table}

\begin{table}[tbp]
\caption{Orbital elements of \gls{ltb} at its escape in the \gls{cr3bp} model}
\centering{}\label{tab: oe LTB CR3BP}
\begin{tabular}{cccccc}
\hline
\noalign{\vskip\doublerulesep}
$\tilde{a}_I$ [LU] & $\tilde{e}_I$ [-] & $\tilde{i}_I$ [deg] & $\tilde{\Omega}_I$ [deg] & $\tilde{\omega}_I$ [deg] & $\tilde{\nu}_I$ [deg] \tabularnewline[\doublerulesep]
\hline
\noalign{\vskip\doublerulesep}
\noalign{\vskip\doublerulesep}
$1.6839$ & $0.2282$ & $3.434$ & $124.9858$ & $213.7120$  & $313.9816$ \tabularnewline[\doublerulesep]
\noalign{\vskip\doublerulesep}
\end{tabular}
\end{table}

The uniform distribution of \glspl{bc} around \gls{ltb} (red star) in both \cref{fig: distance metric CR3BP aei} and \cref{fig: distance metric CR3BP aeomega} confirms that the proposed distance metric effectively identifies trajectories that are dynamically similar to \gls{ltb}.

\subsection{Transition of BCs into ephemeris model} \label{sec: transition to ephemeris algorithm}

Having selected the sub-database $\mathcal{C}_{\tilde{d}v}(\Gamma,z,\zeta) \subset \mathcal{C}(\Gamma,z,\zeta)$ using the distance metric $\tilde{d}_v$, the transition to an ephemeris model can be achieved. First, the transformation of the initial conditions is obtained following the procedure of \cref{sec: transformation of ICs} for the epoch $T_{\text{ETD}}$ (Day 181). 
The precise value used is presented in Appendix~\ref{appendix: sample ICs}. Then, a further adaptation of the polygonal-based algorithm is proposed to follow the model transition, and avoid loss of information in the process.

\cref{alg: update algorithm ephemeris} summarizes the method used to compute portions of the capture set $\mathcal{C}^{\text{eph}}(\Gamma,z,\zeta)$ in the ephemeris model that follow the metric $d_v<60$ m/s, such that a set of \glspl{bc} $\mathcal{C}_{dv}^{\text{eph}}(\Gamma,z,\zeta)$ tailored on \gls{ltb} can be build.
This process is very similar to the previously presented polygonal-based methods for the generation of the capture sets in the spatial \gls{cr3bp}. However, some main differences are here underlined.
\begin{itemize}
    \item The first part of the process focuses on extracting all the capture sets $\mathcal{C}_{\tilde{d}v}(\Gamma_i, z_j, \zeta_k)$ that belong to $\mathcal{C}_{\tilde{d}v}(\Gamma,z,\zeta)$, and on finding the corresponding triplets of parameters $(i,j,k)_{\tilde{d}v}$.
    \item For each triplet $(i,j,k)_{\tilde{d}v}$, the corresponding capture set $\mathcal{C}(\Gamma_i, z_j, \zeta_k)$ — originally computed in the \gls{cr3bp} — is retrieved. 
    \item A flag is implemented within the \textit{polybuffer} while loop (lines~\ref{linealg: start while polybuff}-\ref{linealg: end while polybuff} of \cref{alg: update algorithm ephemeris}) to ensure that at least a partial transition to the ephemeris capture set is achieved. Statistical analysis indicates that approximately $98\%$ of transitions are successfully completed. The remaining $2\%$ of cases correspond to instances where dynamical differences in the ephemeris model alter the structure of the capture sets, causing some \glspl{bc} to disappear or not to be retrieved by this method.
    \item An additional step (lines~\ref{linealg: start coarse grid}-\ref{linealg: end coarse grid} of \cref{alg: update algorithm ephemeris}) is introduced, and executed only within a coarser grid $\mathcal{G}_{dv}$ with stepsize $h_{dv}>h$. This grid enables a more computationally efficient search by performing a backward propagation to determine the orbital elements $a_T$, $e_T$, $i_T$, $\Omega_T$, and $\omega_T$, from which the distance metric $d_v$ is then computed. In this case, $a_I$, $e_I$, $i_I$, $\Omega_I$, and $\omega_I$ are used as the reference orbital elements of \gls{ltb}. Following~\cite{BC-journal}, it is established that the orbital elements of \glspl{bc} origins around Earth exhibit smooth behavior rather than chaotic variation. Therefore, the coarser grid is sufficient to approximate the region where \glspl{bc} close to \gls{ltb} are located, applying a slightly relaxed threshold of $d_v<70$ m/s to ensure completeness in this intermediate selection process. This refinement reduces the search space, $(\mathcal{S}_{i, j, k})_{dv} \subset \mathcal{S}_{i, j, k}$, optimizing the selection of candidate trajectories before full propagation. As previously stated, the final capture set $\mathcal{C}_{dv}^{\text{eph}}(\Gamma,z,\zeta)$ is constructed with a stricter threshold of $d_v<60$ m/s.
    \item Given that the scope of the current work is to provide robust ballistic insertions, only \glspl{bc} presenting two or more revolutions are considered ad included in the capture set. 
\end{itemize}

\begin{algorithm}[tbp]
\caption{\gls{etd} Ballistic Capture Set transition into an ephemeris model}
\label{alg: update algorithm ephemeris}
\begin{algorithmic}[1]
\State {Build a grid $\mathcal{G}$ around $M_2$ with width $x_{\mathcal{G}}$, height $y_{\mathcal{G}}$, and  stepsize $h$. Build another coarser grid $\mathcal{G}_{dv}$ with $h_{dv} = 10\cdot h$. Initialize a vector of $P$ triplets to be addressed $(i,j,k)_{\tilde{d}v}$, with overall dimensions $P \times 3$.} 
\For{$p \gets 1$ to $P$}
    \State {$(i,j,k) \gets (i_p, j_p, k_p)_{\tilde{d}v}$.}
    \State {Compute corresponding values $(\Gamma_i, z_j, \zeta_k)$ and load the section of capture set $\mathcal{C}_{i, j, k} = \mathcal{C}(\Gamma_i, z_j, \zeta_k)$.}
    \State {$\mathcal{S}_{i, j, k} \gets \mathcal{C}_{i, j, k}$; $\mathcal{C}_{i, j, k}^{\text{eph}} \gets \mathcal{C}_{i, j, k}$} \Comment{The search space $\mathcal{S}_{i, j, k}$ and the candidate $\mathcal{C}_{i, j, k}^{\text{eph}}$ are initialized as $\mathcal{C}_{i, j, k}$}
    \State{$\text{flag} \gets \text{false}$; \; $\text{iter} \gets 0$}
    \While {$\mathcal{C}_{i, j, k}^{\text{eph}} \cap \partial \mathcal{S}_{i, j, k} \ne \emptyset \; \lor \text{flag} = \text{false}$}   \Comment {Check for \glspl{bc} on the boundary of $\mathcal{S}_{i, j, k}$}  \label{linealg: start while polybuff}
        \State {$\partial \mathcal{S}_{i, j, k} \leftarrow \partial \mathcal{S}_{i, j, k} +\textrm{offset}(d_O)$}     \Comment{The polygon $\partial \mathcal{S}_{i, j, k}$ is enlarged with an offset of size $d_O$}
        \State {Using Eqs.~\eqref{eq: v2_energy_constraint} and \eqref{eq: injection_angle}, compute $(\dot{x},\dot{y},\dot{z})$ on every vertex $(x, y) \in \partial \mathcal{S}_{i,j,k}$ returned by \textit{polybuffer}, obtaining $\mathbf{x}_0$.}
        \State{Use process presented in \cref{sec: transformation of ICs} to transform $\mathbf{x}_0$ into its ephemeris correspondent $\mathbf{x}_{\text{eph},0}$.}
        \State {Propagate each initial condition $\mathbf{x}_{\text{eph},0}$ with \cref{eq: Newton law EMS} for a time  $T_{\text{sp}}\approx54.6$ days.}
        \State{Identify \glspl{bc} according to the definition in \cref{sec: definition BC}.}
        \State{Temporarily add the \glspl{bc} to $\mathcal{C}_{i, j, k}^{\text{eph}}$. If \glspl{bc} are found, set $\text{flag} \gets \text{true}$ (successful transition)}
        \If{$\text{iter}=20$} \Comment {No \glspl{bc} found in any $\partial \mathcal{S}_{i, j, k}$: the transition process may have failed for this set}
            \State {Exit while loop, setting $\mathcal{S}_{i, j, k} \gets \mathcal{C}_{i, j, k}$} \Comment {As a final attempt, search for \glspl{bc} within $\mathcal{C}_{i, j, k}$} 
        \EndIf
        \State{$\text{iter} \gets \text{iter}+1$}
    \EndWhile \label{linealg: end while polybuff}
    \Comment {The set $\mathcal{S}_{i, j, k}$ enclosing $\mathcal{C}_{i, j, k}^{\text{eph}}$ ($\mathcal{C}_{i, j, k}^{\text{eph}} \subset \mathcal{S}_{i, j, k}$) is found} 
    \State{Using Eqs.~\eqref{eq: v2_energy_constraint} and \eqref{eq: injection_angle}, compute $(\dot{x},\dot{y},\dot{z})\; \forall\ (x, y) \in (\mathcal{G}_{dv} \cap \mathcal{S}_{i,j,k})$, obtaining $\mathbf{x}_0$.} \label{linealg: start coarse grid}
    \State{Use process presented in \cref{sec: transformation of ICs} to transform $\mathbf{x}_0$ into its ephemeris correspondent $\mathbf{x}_{\text{eph},0}$.}
    \State{Using \cref{eq: Newton law EMS}, propagate each initial condition $\mathbf{x}_{\text{eph},0}$ backward for a time  $T_B\approx54.6$ days.}
    \State{Evaluate the distance metric for each final state and extract the subset $(\mathcal{S}_{i, j, k})_{dv} \subset \mathcal{S}_{i, j, k}$ satisfying $d_v<70$ m/s.} \label{linealg: end coarse grid}
    
    \State{Using Eqs.~\eqref{eq: v2_energy_constraint} and \eqref{eq: injection_angle}, compute $(\dot{x},\dot{y},\dot{z})\; \forall\ (x, y) \in [\mathcal{G} \cap (\mathcal{S}_{i, j, k})_{dv}]$, obtaining $\mathbf{x}_0$.}
    \State{Use process presented in \cref{sec: transformation of ICs} to transform $\mathbf{x}_0$ into its ephemeris correspondent $\mathbf{x}_{\text{eph},0}$.}
    \State{Verify that $\mathbf{x}_{\text{eph},0}$ meets condition $\dot{\varepsilon_2}(\tau_0)<0$ with \cref{eq: check derivative eps2}, otherwise discard it.}
    \State{Using \cref{eq: Newton law EMS}, propagate each initial condition $\mathbf{x}_{\text{eph},0}$ backward for a time  $T_B\approx54.6$ days and possibly forward for a time $T_{s}\approx273$ days.}
    \State{Identify only \glspl{bc} with 2+ revolutions and add them to $\mathcal{C}_{i, j, k}^{\text{eph}}$ according to the definition in \cref{sec: definition BC}.}
\EndFor
\State{After iterating through all $(i,j,k)_{\tilde{d}v}$ triplets, the final capture set $\mathcal{C}_{dv}^{\text{eph}}(\Gamma,z,\zeta)$ consists of \glspl{bc} that have been successfully transitioned into the ephemeris model.}
\end{algorithmic}
\end{algorithm}


The parameters used in \cref{alg: update algorithm ephemeris} are identical to the ones introduced in \cref{sec: sample spatial capture sets} for the previous algorithms. 

\section{Results} \label{sec: results}

This section provides an overall analysis of the database built upon the ephemeris-model capture set, $\mathcal{C}_{dv}^{\text{eph}}(\Gamma,z,\zeta)$, and selects representative trajectories to demonstrate the effectiveness of the proposed approach.

As introduced in \cref{sec: sample spatial capture sets}, the subscripted notation $a_T$, $e_T$, $i_T$, $\Omega_T$, $\omega_T$, and $\nu_T$ refers to the orbital elements of the \gls{bc} origin with respect to Earth (in this section, they are calculated in the ephemeris model). In contrast, un-subscripted symbols (e.g., $i$) will be used throughout the following sections to denote the Moon-centered osculating orbital elements at selected perilunes.

\subsection{Database analysis} \label{sec: capture set analysis}
\begin{figure}[tbp]
    \centering
    \begin{subfigure}[t]{\widthof{\includegraphics[height=5cm]{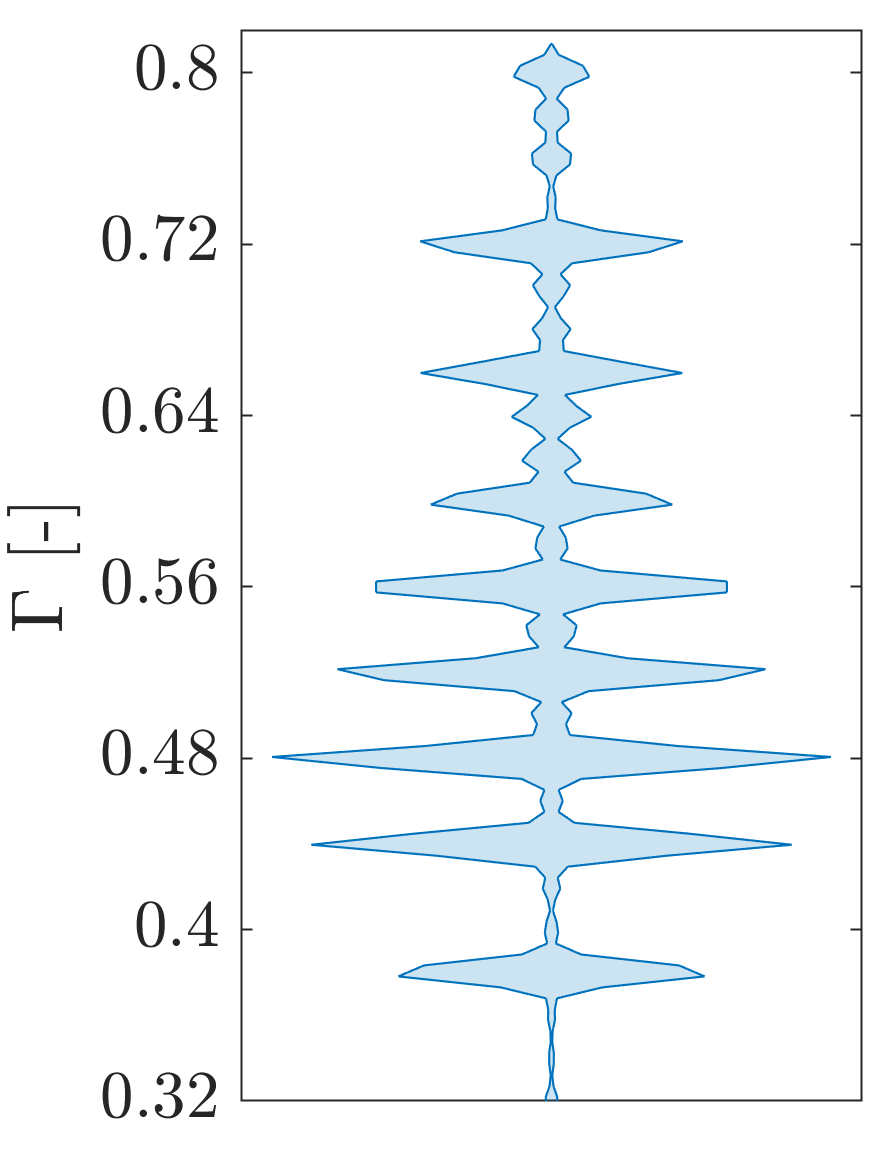}}}
        \includegraphics[height=5cm]{Figures/GlobalRes/GlobalRes_Gamma.png}
        \subcaption{Distribution of $\Gamma$}
        \label{fig: GlobalRes Gamma}
    \end{subfigure} \hspace{0.5cm}
    \begin{subfigure}[t]{\widthof{\includegraphics[height=5cm]{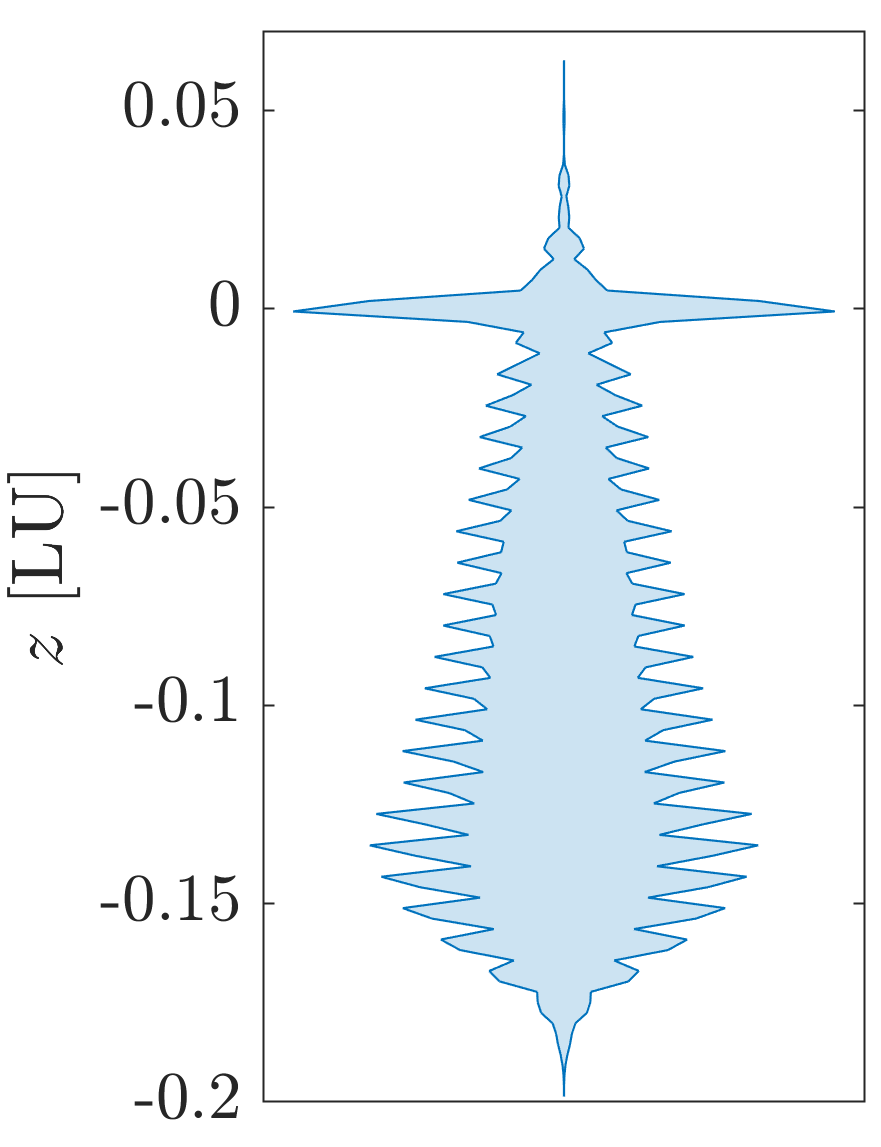}}}
        \includegraphics[height=5cm]{Figures/GlobalRes/GlobalRes_z.png}
        \subcaption{Distribution of $z$}
        \label{fig: GlobalRes z}
    \end{subfigure} \hspace{0.5cm}
    \begin{subfigure}[t]{\widthof{\includegraphics[height=5cm]{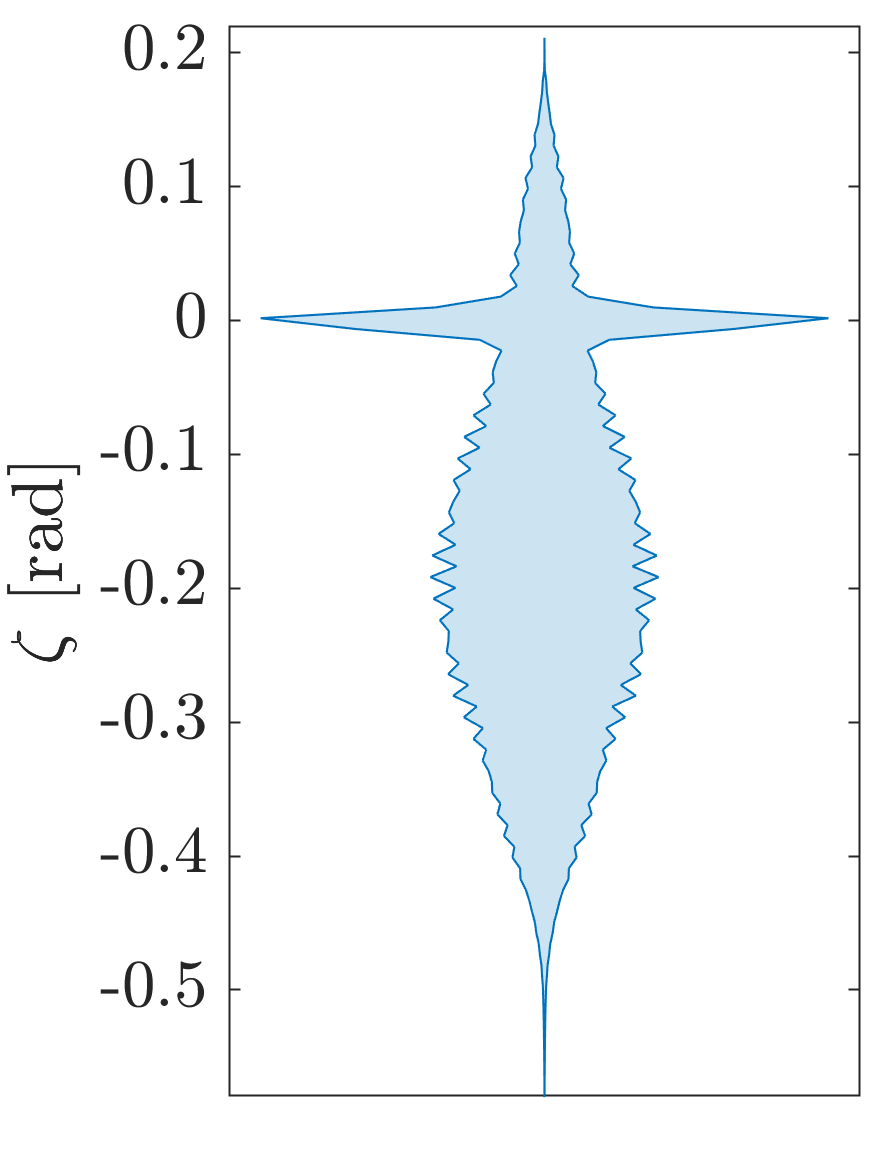}}}
        \includegraphics[height=5cm]{Figures/GlobalRes/GlobalRes_zeta.png}
        \subcaption{Distribution of $\zeta$}
        \label{fig: GlobalRes zeta}
    \end{subfigure}
    \caption{Distribution of the \gls{etd} state parameters in the dataset $\mathcal{C}^{\text{eph}}_{dv}(\Gamma,z,\zeta)$.}
    \label{fig: global res params}
\end{figure}

\begin{figure}[tbp]
    \centering
    \begin{subfigure}[t]{\widthof{\includegraphics[height=5cm]{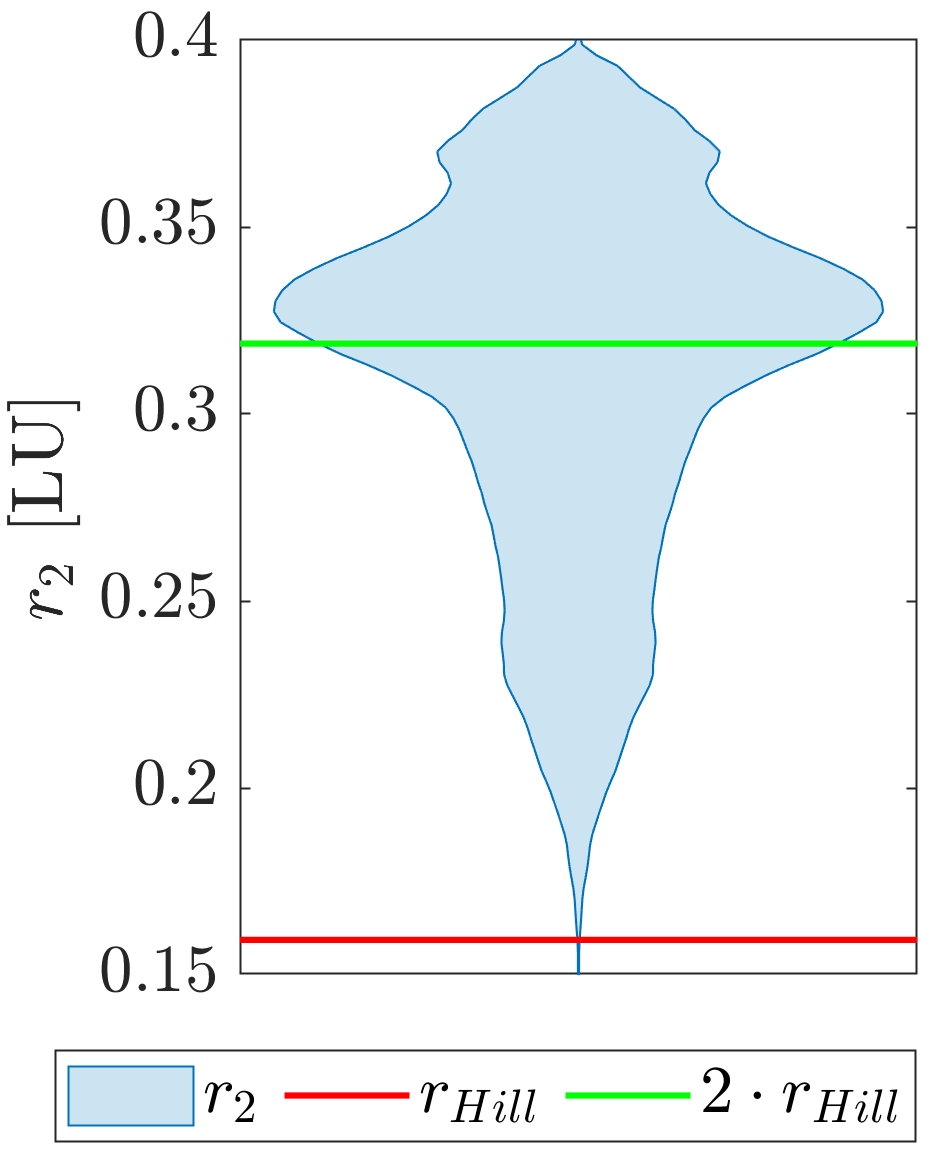}}}
        \includegraphics[height=5cm]{Figures/GlobalRes/GlobalRes_r2.png}
        \subcaption{Distribution of $r_2$ at \gls{etd}}
        \label{fig: GlobalRes r2}
    \end{subfigure} \hspace{0.5cm}
    \begin{subfigure}[t]{\widthof{\includegraphics[height=5cm]{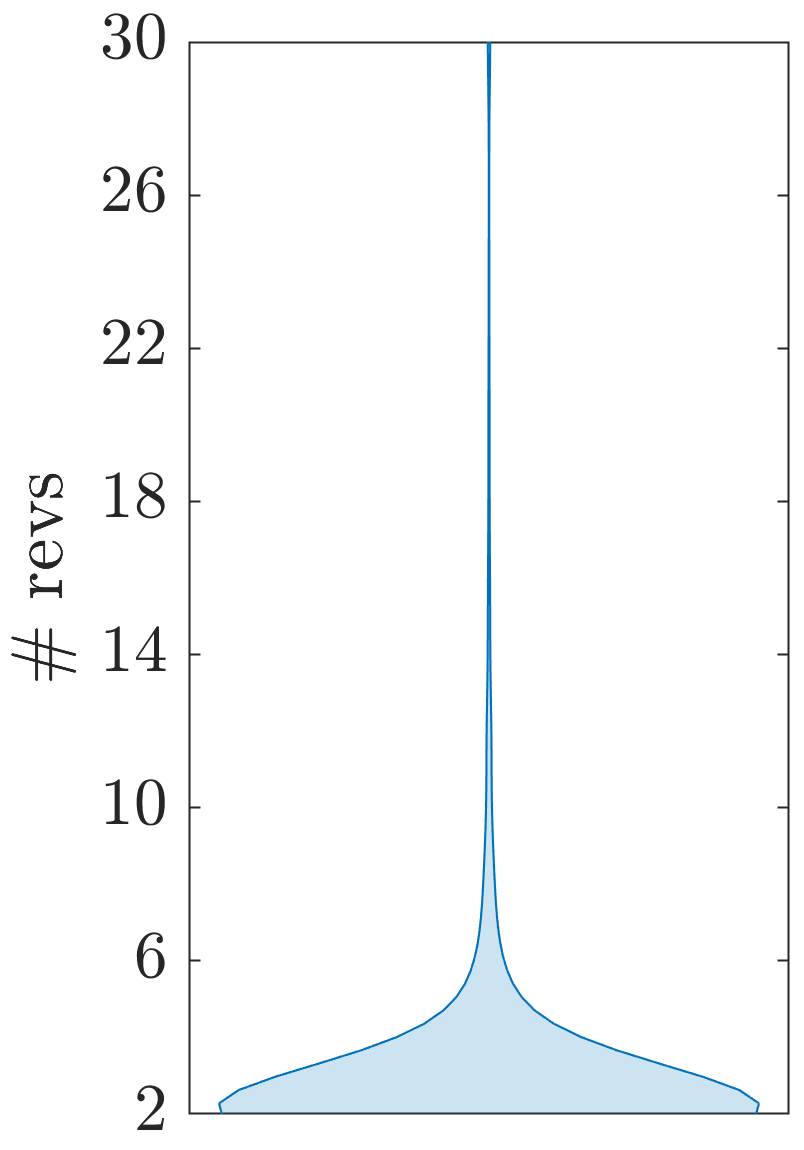}}}
        \includegraphics[height=5cm]{Figures/GlobalRes/GlobalRes_NRevALL.png}
        \subcaption{Distribution of total number of revolutions}
        \label{fig: GlobalRes NRevALL}
    \end{subfigure} \hspace{0.5cm}
    \begin{subfigure}[t]{\widthof{\includegraphics[height=5cm]{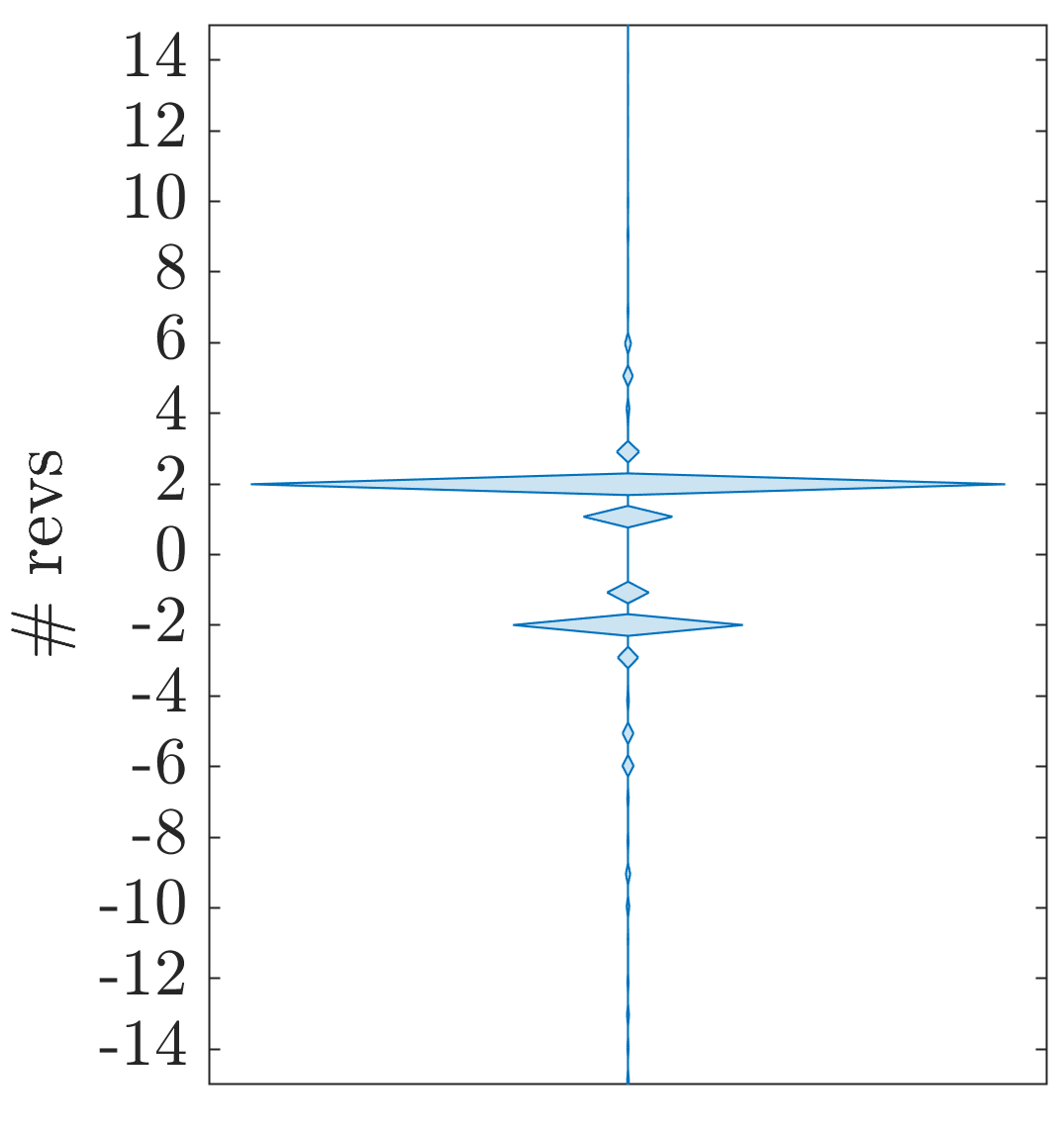}}}
        \includegraphics[height=5cm]{Figures/GlobalRes/GlobalRes_NRev.png}
        \subcaption{Distribution of revolutions with sign}
        \label{fig: GlobalRes NRev}
    \end{subfigure}
    \caption{Distribution of the \gls{etd} state parameter $r_2$ and number of revolutions in the dataset $\mathcal{C}^{\text{eph}}_{dv}(\Gamma,z,\zeta)$.}
    \label{fig: global res r2 revs}
\end{figure}

The dataset $\mathcal{C}_{dv}^{\text{eph}}(\Gamma,z,\zeta)$, containing approximately 20 million \glspl{bc}, is analyzed to highlight key features.
First, the parameter distribution $(\Gamma,z,\zeta)$ is summarized in \cref{fig: global res params}, showing the discrete steps from \cref{tab: param algs} and the increased step sizes in $\Delta z$ and $\Delta \Gamma$ used in \cref{alg: search algorithm zdot}.
The $(x, y)$ component distribution is presented in terms of distance from the Moon, $r_2$, in \cref{fig: GlobalRes r2}. The longevity of \glspl{bc} is assessed through the total number of completed revolutions, shown in \cref{fig: GlobalRes NRevALL}, with individual and cumulative revolution counts detailed in \cref{tab: percentages revs all}. As shown in \cref{fig: GlobalRes NRev}, prograde motion (positive sign) is more prevalent, while retrograde orbits (negative sign), known for their higher stability, dominate at higher revolution counts. \cref{tab: percentages revs} further distinguishes between prograde and retrograde cases. Notably, \glspl{bc} with a single revolution can appear because a transition between prograde and retrograde motion (and vice versa) may occur.

\begin{figure}[tbp]
    \centering
    \begin{subfigure}[t]{\widthof{\includegraphics[height=5cm]{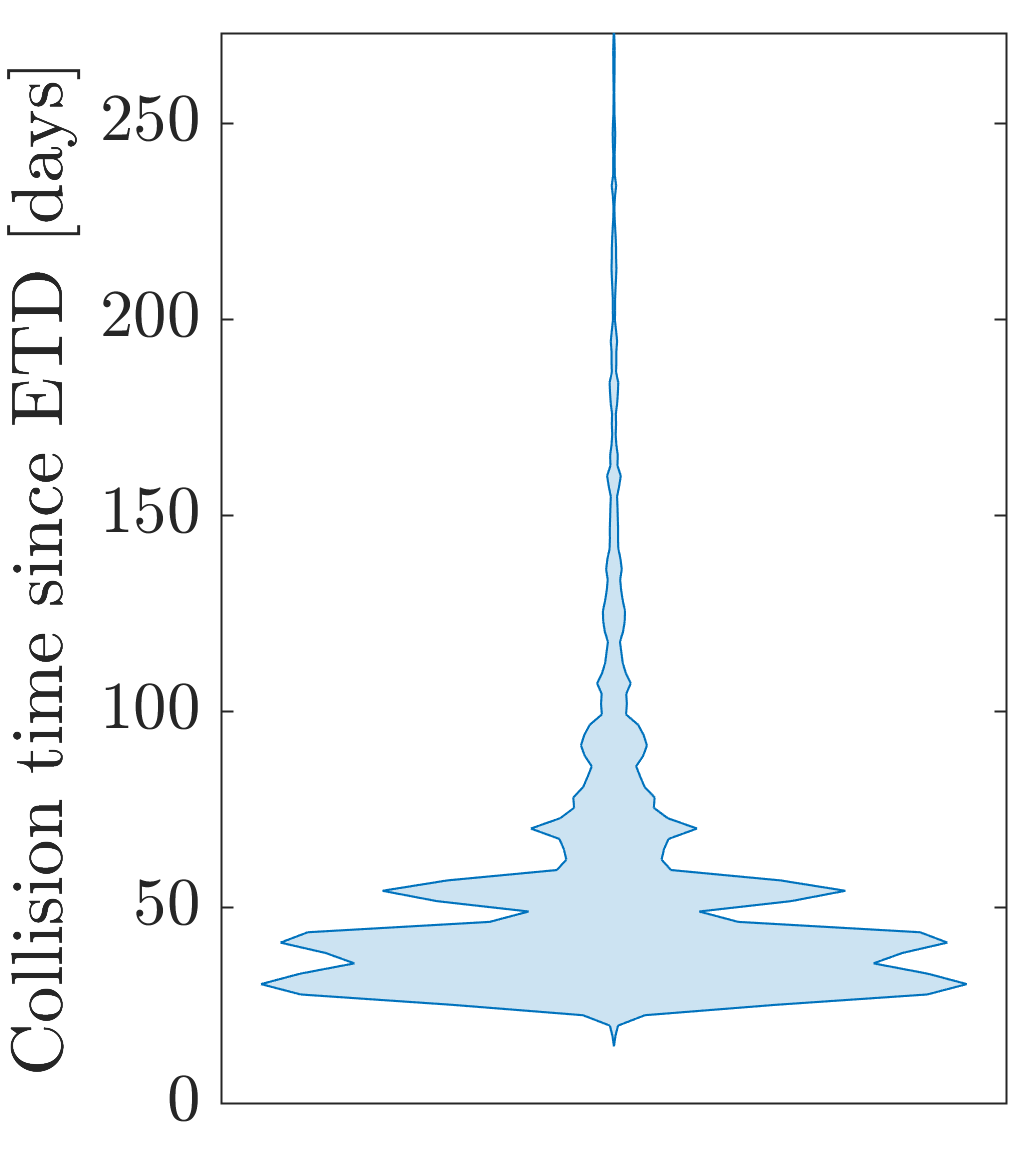}}}
        \includegraphics[height=5cm]{Figures/GlobalRes/GlobalRes_Coll.png}
        \subcaption{Distribution of collision times}
        \label{fig: GlobalRes coll}
    \end{subfigure} \hspace{0.5cm}
    \begin{subfigure}[t]{\widthof{\includegraphics[height=5cm]{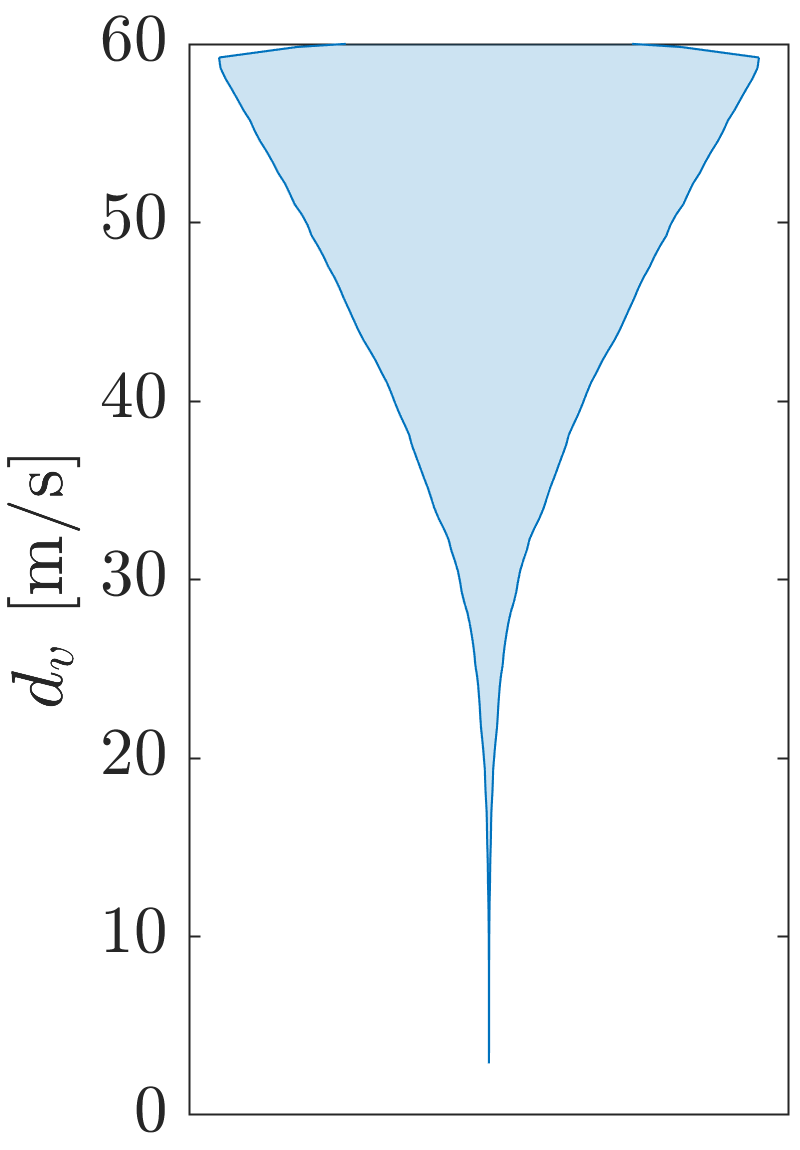}}}
        \includegraphics[height=5cm]{Figures/GlobalRes/GlobalRes_dv.png}
        \subcaption{Distribution of $d_v$}
        \label{fig: GlobalRes dv}
    \end{subfigure}
    \caption{Distribution of collision times and $d_v$ in $\mathcal{C}^{\text{eph}}_{dv}(\Gamma,z,\zeta)$.}
    \label{fig: global res coll dv}
\end{figure}

\begin{figure}[tbp]
    \centering
    \begin{subfigure}[t!]{0.47\textwidth}
        \includegraphics[width=\textwidth]{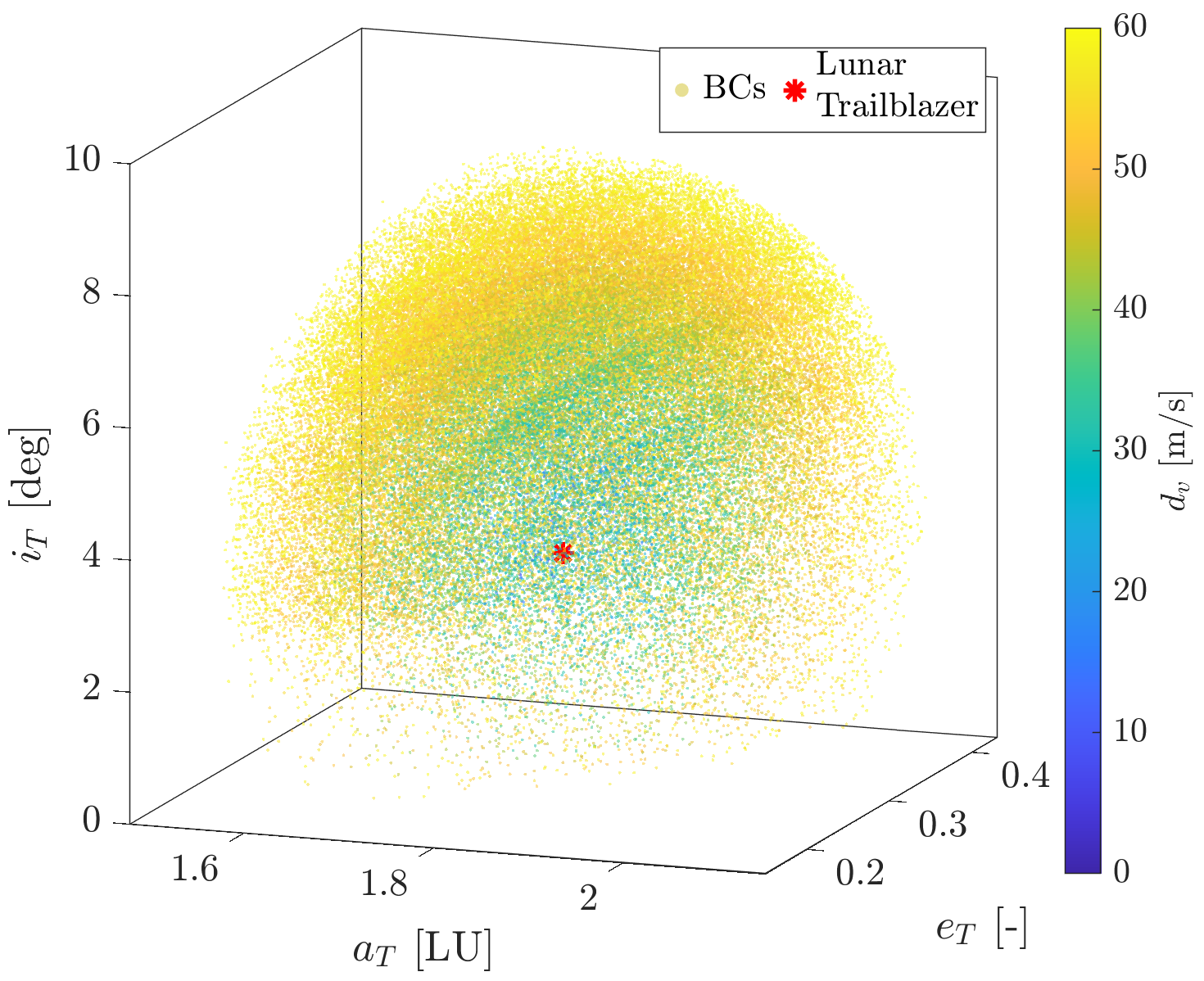}
        \subcaption{Distance metric $d_v$ in the orbital elements $a_T$, $e_T$, $i_T$}
        \label{fig: distance metric ephem aei}
    \end{subfigure}
    \begin{subfigure}[t!]{0.50\textwidth}
        \includegraphics[width=\textwidth]{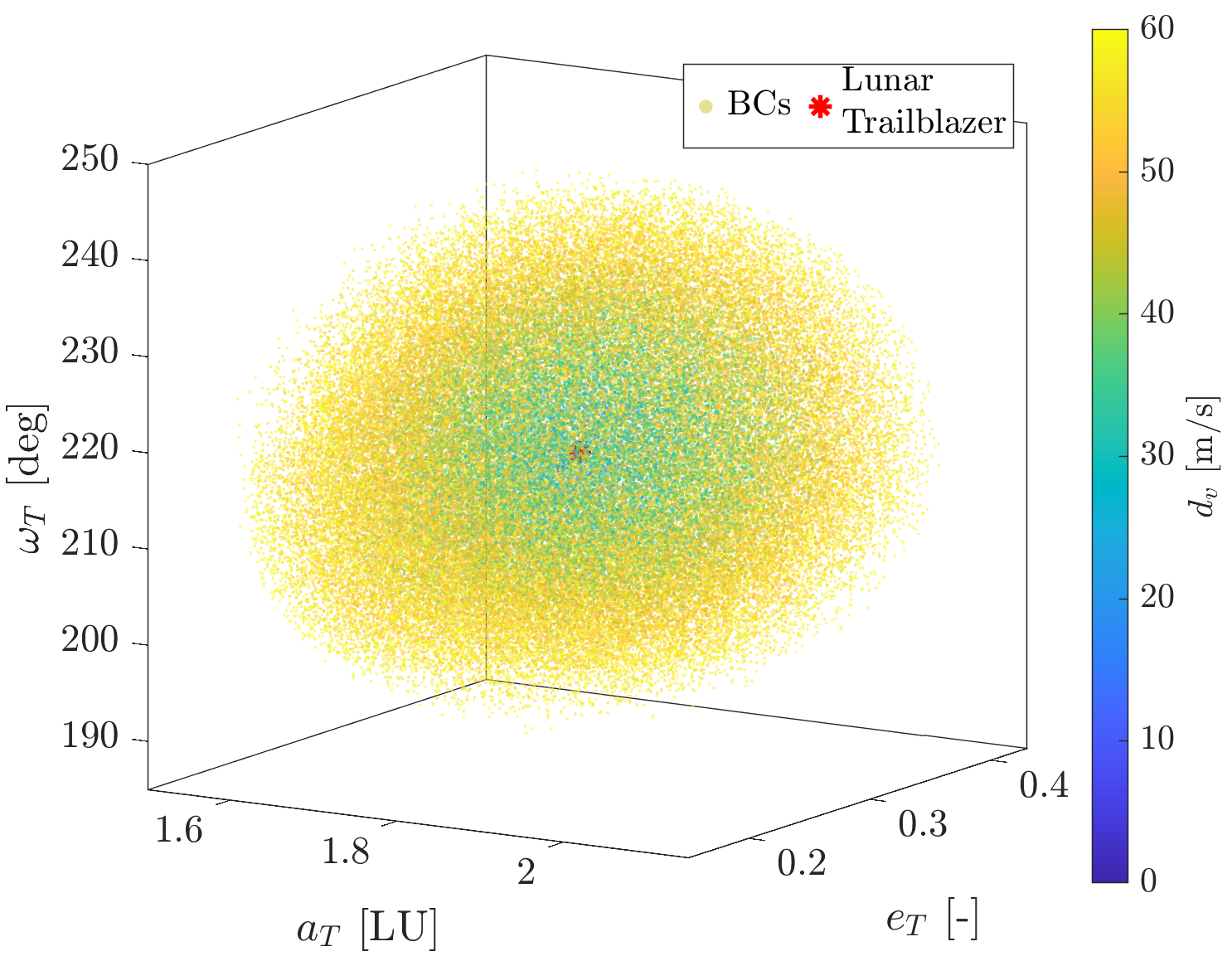}
        \subcaption{Distance metric $d_v$ in the orbital elements $a_T$, $e_T$, $\omega_T$}
        \label{fig: distance metric ephem aeomega}
    \end{subfigure}
    \caption{Distance metric $d_v$.}
    \label{fig: distance metric ephem}
\end{figure}

\cref{fig: global res coll dv} presents the distribution of collision times since \gls{etd} and the distance metric $d_v$, computed from the full-ephemeris orbital elements of \gls{ltb} ($a_I$, $e_I$, $i_I$, $\Omega_I$, and $\omega_I$). The latter shows an almost exponential increase in the number of \glspl{bc} as the threshold $d_v^{\text{thr}}$ increases.
\cref{fig: distance metric ephem} further illustrates $d_v$ as a function of the initial orbital elements. The contribution of $\Omega$ is negligible and is thus omitted. Only $0.5\%$ of the dataset $\mathcal{C}_{dv}^{\text{eph}}(\Gamma,z,\zeta)$ is represented here. The nominal \gls{ltb} point lies within the cloud of \gls{bc} points. From \cref{fig: distance metric ephem aei}, it is evident that \glspl{bc} are less frequent at very low inclinations. 

\begin{figure}[tbp]
    \begin{minipage}[c]{0.49\linewidth}
        \centering
        \includegraphics[width=\linewidth]{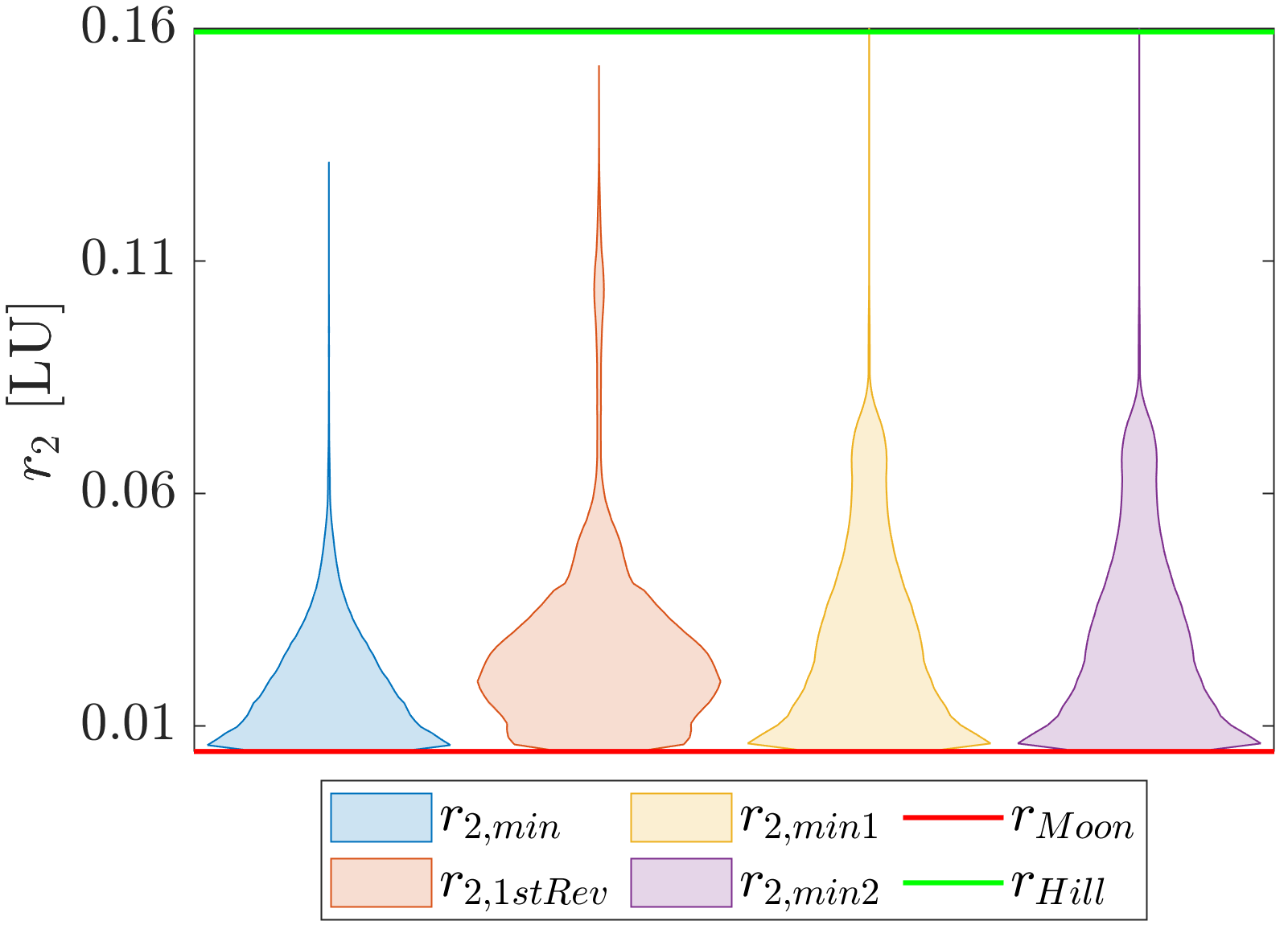}
        \caption{Distributions of the perilune distances in the dataset $\mathcal{C}^{\text{eph}}_{dv}(\Gamma,z,\zeta)$.}
        \label{fig: global res r2mins}
    \end{minipage}
    \hfill
    \begin{minipage}[c]{0.49\linewidth}
        \centering
        \includegraphics[width=\linewidth]{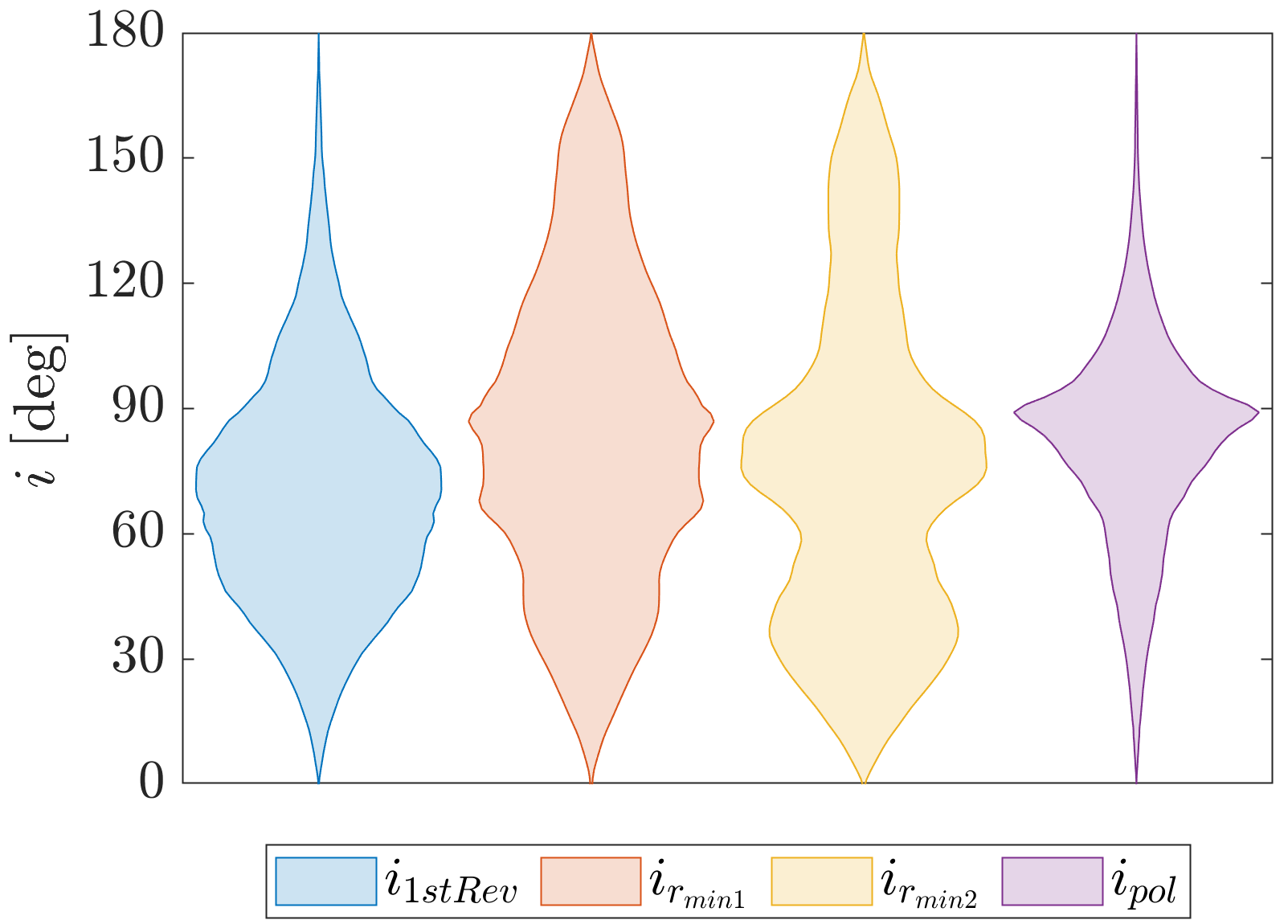}
        \caption{Distributions of inclinations at perilune in the dataset $\mathcal{C}^{\text{eph}}_{dv}(\Gamma,z,\zeta)$.}
        \label{fig: global res incls}
    \end{minipage}
\end{figure}

Finally, insights into the \glspl{bc} perilunes are provided. \cref{fig: global res r2mins} shows the distributions of the most relevant minimum distances to the Moon. The first perilune distance is denoted as $r_{2,1stRev}$, while $r_{2,min1}$ and $r_{2,min2}$ correspond to the two closest approaches occurring after the first perilune (with $r_{2,min2} > r_{2,min1}$), and can serve to identify repeated close encounters. For reference, the overall minimum perilune distance, $r_{2,min}$, is also shown—it may coincide with either $r_{2,1stRev}$ or $r_{2,min1}$, depending on the specific trajectory. These distributions highlight the availability of low-to-moderate altitude perilunes, offering valuable opportunities for mission design.
Moreover, \cref{fig: global res incls} presents the inclinations of the osculating orbits around the Moon at key perilunes. Following the same nomenclature, $i_{1stRev}$ denotes the inclination at the first perilune, while $i_{r_{min1}}$ and $i_{r_{min2}}$ correspond to the inclinations at the two closest perilunes following the first one (with $r_{2,min2} > r_{2,min1}$). Additionally, the most polar perilune inclination is included for statistical purposes, showing a significant number of polar or near-polar perilunes. This is particularly relevant for \gls{ltb} mission design. Sample trajectories using a combination of these parameters will be investigated in the following.

\begin{table}[tbp]
\caption{Distribution of total number of revolutions in $\mathcal{C}^{\text{eph}}_{dv}(\Gamma,z,\zeta)$}
\centering{}\label{tab: percentages revs all}
\begin{tabular}{cccccccc}
\hline
\noalign{\vskip\doublerulesep}
 & 2 revs & 3 revs & 4 revs & 5 revs & 6 revs & 7 revs & 8 revs \tabularnewline[\doublerulesep]
\hline
\noalign{\vskip\doublerulesep}
\noalign{\vskip\doublerulesep}
Individual $\#$ revs & $72.7 \%$ & $14.2 \%$ & $5.1 \%$ & $2.5 \%$  & $1.2 \%$ & $0.8 \%$  & $0.7 \%$ \tabularnewline[\doublerulesep]
\noalign{\vskip\doublerulesep}
\noalign{\vskip\doublerulesep}
Cumulative $\#$ revs & $72.7 \%$ & $86.9 \%$ & $92.0 \%$ & $94.5 \%$ & $95.7 \%$  & $96.5 \%$ & $97.2 \%$ \tabularnewline[\doublerulesep]
\noalign{\vskip\doublerulesep}
\end{tabular}
\end{table}

\begin{table}[tbp]
\caption{Distribution of the prograde/retrograde number of revolutions in $\mathcal{C}^{\text{eph}}_{dv}(\Gamma,z,\zeta)$}
\centering{}\label{tab: percentages revs}
\begin{tabular}{cccccccccc}
\hline
\noalign{\vskip\doublerulesep}
 & 1 rev & 2 revs & 3 revs & 4 revs & 5 revs & 6 revs & 7 revs & 8 revs & >9 revs \tabularnewline[\doublerulesep]
\hline
\noalign{\vskip\doublerulesep}
\noalign{\vskip\doublerulesep}
Prograde $\#$ revs & $10.5 \%$ & $40.8 \%$ & $3.8 \%$ & $2.5 \%$  & $0.9 \%$ & $0.4 \%$  & $0.2 \%$ & $0.1 \%$ & $0.1 \%$ \tabularnewline[\doublerulesep]
\noalign{\vskip\doublerulesep}
\noalign{\vskip\doublerulesep}
Retrograde $\#$ revs & $4.9 \%$ & $12.4 \%$ & $3.4 \%$ & $1.5 \%$ & $1.0 \%$  & $0.7 \%$ & $0.4 \%$ & $0.5 \%$ & $2.4 \%$ \tabularnewline[\doublerulesep]
\noalign{\vskip\doublerulesep}
\end{tabular}
\end{table}

\subsection{Sample BCs} \label{sec: sample BCs for LTB}
A simple selection process is implemented to identify \glspl{bc} with specific features within the database. The sample trajectories are visualized in multiple reference frames within the ephemeris model. A synodic frame rotating with the Moon is employed, with nondimensional coordinates obtained through the inverse procedure from \cref{sec: transformation of ICs}. Additionally, two inertial frames centered on Earth and Moon are considered, both with their $x_1$- and $x_2$-axes aligned with the Earth–Moon direction at time $T_{\text{ETD}}$, as previously introduced in \cref{sec: ephemeris model}. All coordinates are expressed in nondimensional form, using the mean Earth–Moon distance defined in \cref{tab: scaling units}. Initial conditions for each \gls{bc} are provided in Appendix~\ref{appendix: sample ICs}.

\subsubsection{Longest BCs} \label{sec: longest BCs}
\begin{figure}[tbp]
    \centering
    \begin{subfigure}[t!]{0.49\textwidth}
        \includegraphics[width=\textwidth]{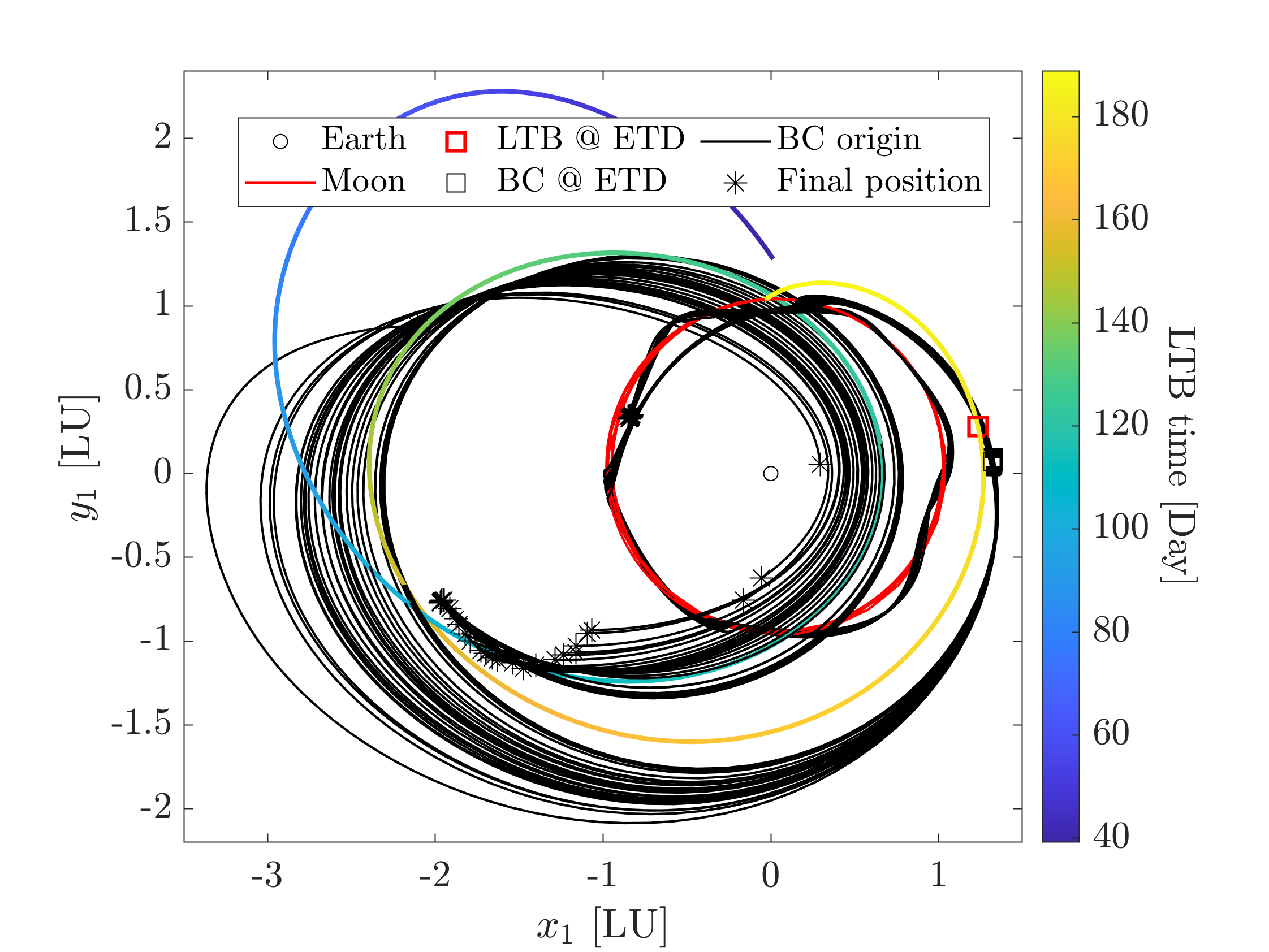}
    \end{subfigure}
    \\
    \begin{subfigure}[t!]{0.49\textwidth}
        \includegraphics[width=\textwidth]{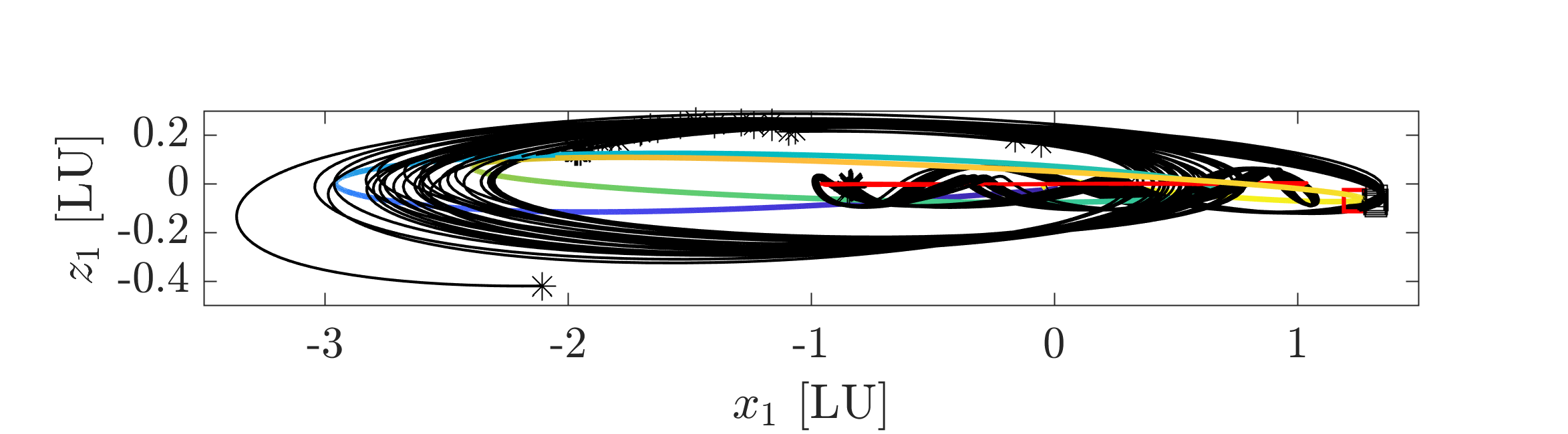}
    \end{subfigure}
    \caption{Longest \glspl{bc} corridor and \gls{ltb} in Earth inertial frame.}
    \label{fig: BC corridor inert1}
\end{figure}

\begin{figure}[tbp]
    \centering
    \begin{subfigure}[t!]{0.49\textwidth}
        \includegraphics[width=\textwidth]{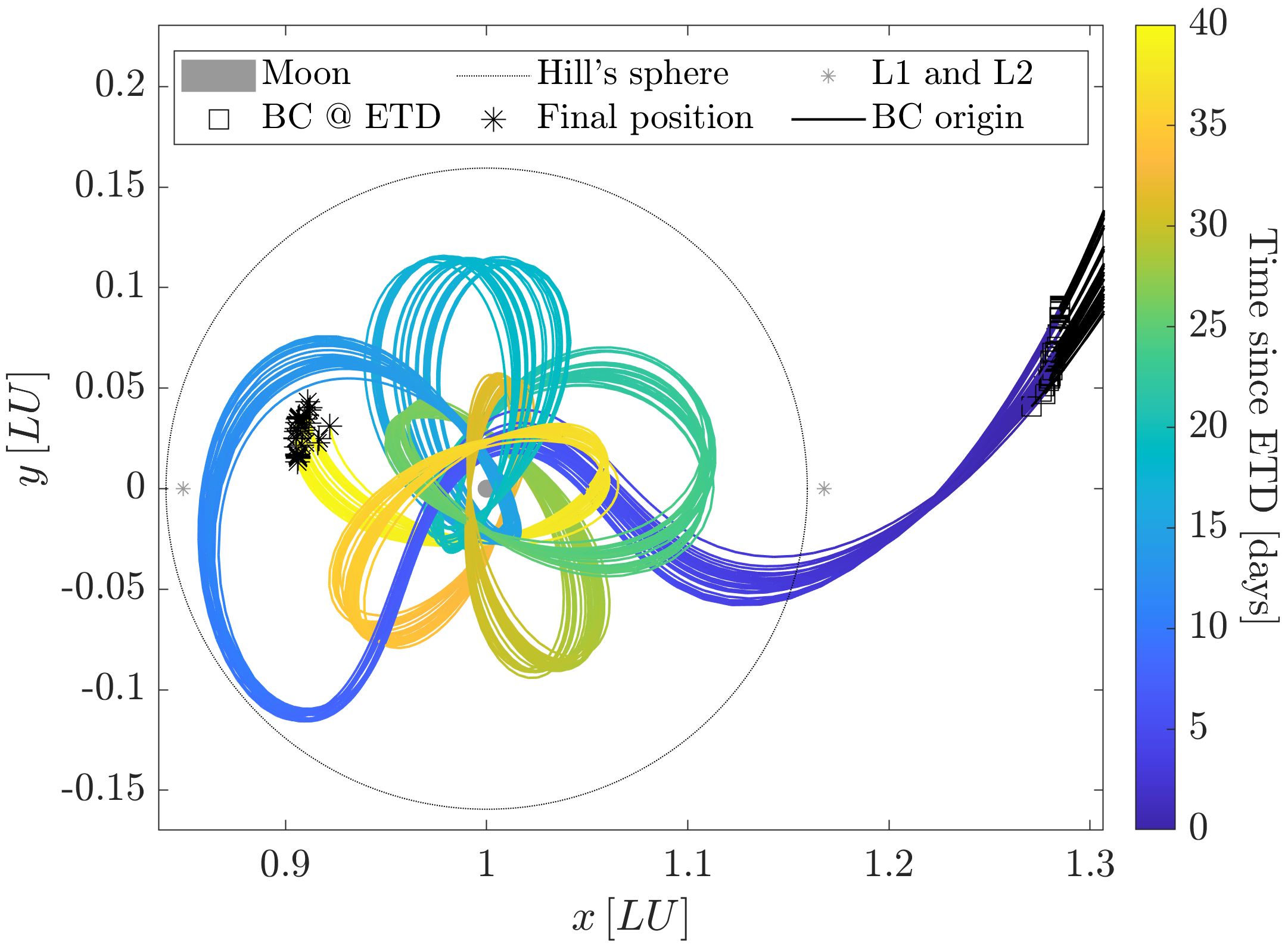}
    \end{subfigure}
    \\
    \begin{subfigure}[t!]{0.49\textwidth}
        \includegraphics[width=\textwidth]{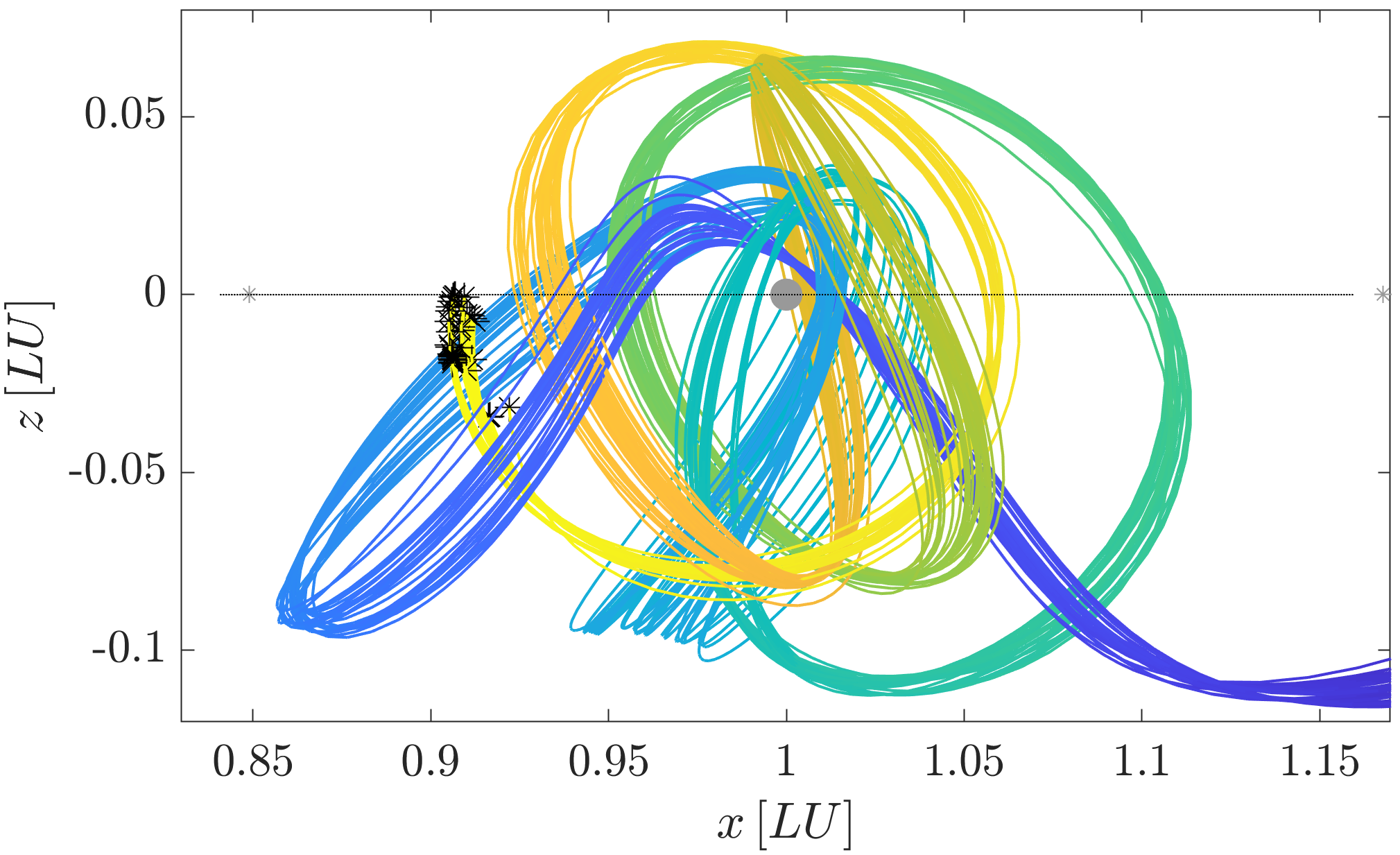}
    \end{subfigure}
    \caption{Longest \glspl{bc} corridor in the synodic frame.}
    \label{fig: BC corridor syn}
\end{figure}

\begin{figure}[tbp]
    \centering
    \includegraphics[width=0.63\textwidth]{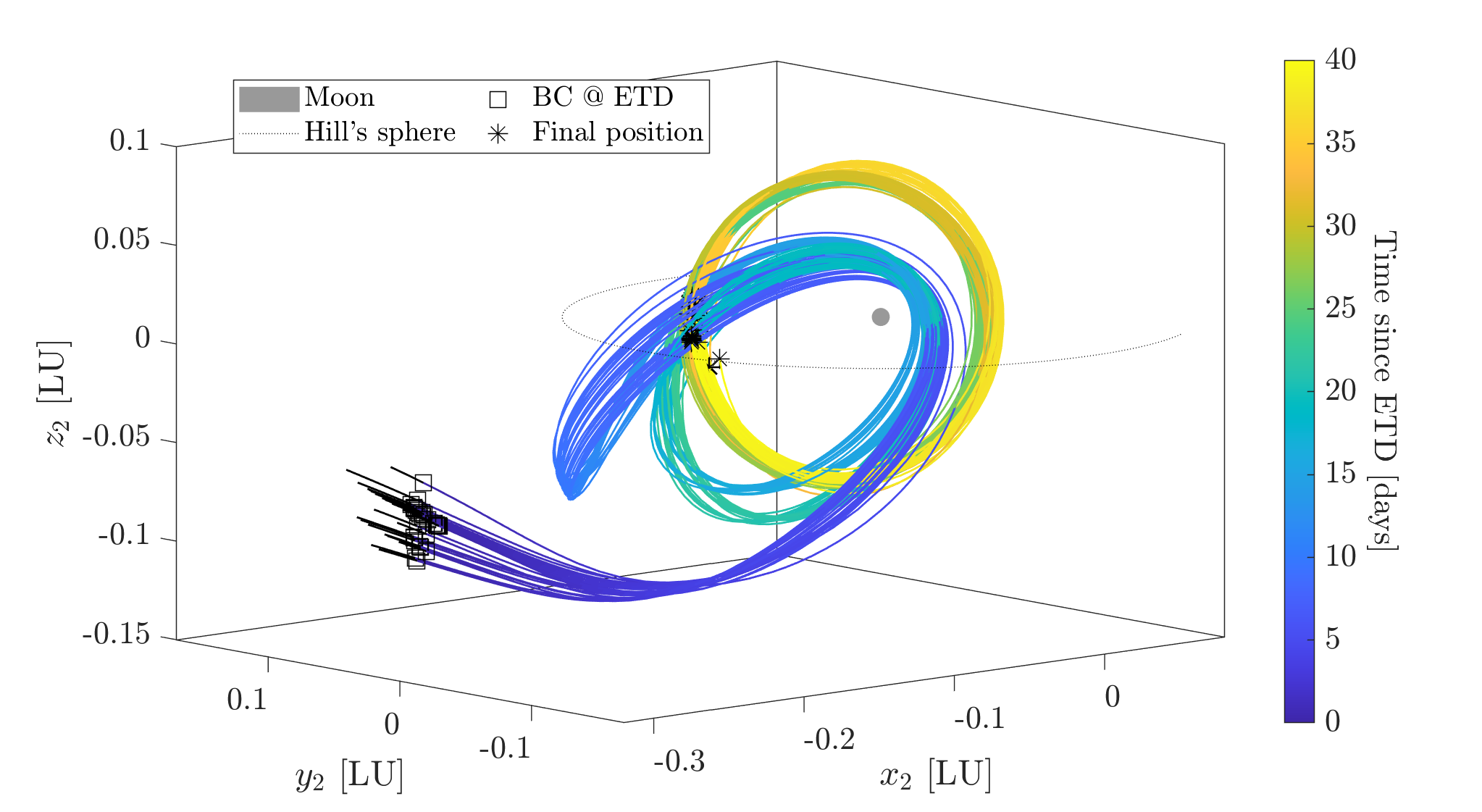}
    \caption{Longest \glspl{bc} corridor in Moon inertial frame.}
    \label{fig: BC corridor inert2}
\end{figure}

First, a representative \gls{bc} corridor is shown in Figs.~\ref{fig: BC corridor inert1}, \ref{fig: BC corridor syn}, and \ref{fig: BC corridor inert2}. In this context, a corridor refers to a subset of the longest available \glspl{bc} (i.e. completing at least 45 revolutions) within the capture set $\mathcal{C}_{dv}^{\text{eph}}(\Gamma,z,\zeta)$. These \glspl{bc} share similar characteristics - such as comparable orbital elements at their origin - and collectively exhibit a recognizable pattern. For instance, the convergence of trajectories as they approach the \gls{etd} in \cref{fig: BC corridor inert1} illustrates this behavior. Similarly, the distance metric can be visualized as the difference between \glspl{bc} origin (black line) and \gls{ltb} nominal trajectory (colored line), with an average of $d_v \sim 45$~m/s. This is mostly due to the difference in inclination (see $x_1$–$z_1$ view), which is slightly higher for the \glspl{bc}.
While the term corridor typically denotes a subset with consistent features, the entire capture set $\mathcal{C}_{dv}^{\text{eph}}(\Gamma,z,\zeta)$ can be viewed as a broader corridor composed of trajectories originating near the nominal \gls{ltb} trajectory.

\Cref{fig: BC corridor syn} provides insight into the underlying reason for the unusual stability of these trajectories. A highly inclined prograde approach quickly transitions into nearly polar, retrograde revolutions. Over time, the average distance from the Moon decreases, and the \glspl{bc} remain well within the Hill sphere. Such behavior does not occur in the \gls{cr3bp}, indicating that these trajectories fully leverage the perturbations present in the ephemeris model — primarily the Moon's orbital ellipticity around Earth and the gravitational influence of the Sun — to sustain stable motion around the Moon. This effect is further illustrated in \cref{fig: BC corridor inert2}, where a temporary stabilization into an almost circular orbit is observed, persisting for at least two full revolutions between days 25 and 40.

\begin{figure}[tbp]
    \centering
    \begin{subfigure}[t!]{0.49\textwidth}
        \includegraphics[width=\textwidth]{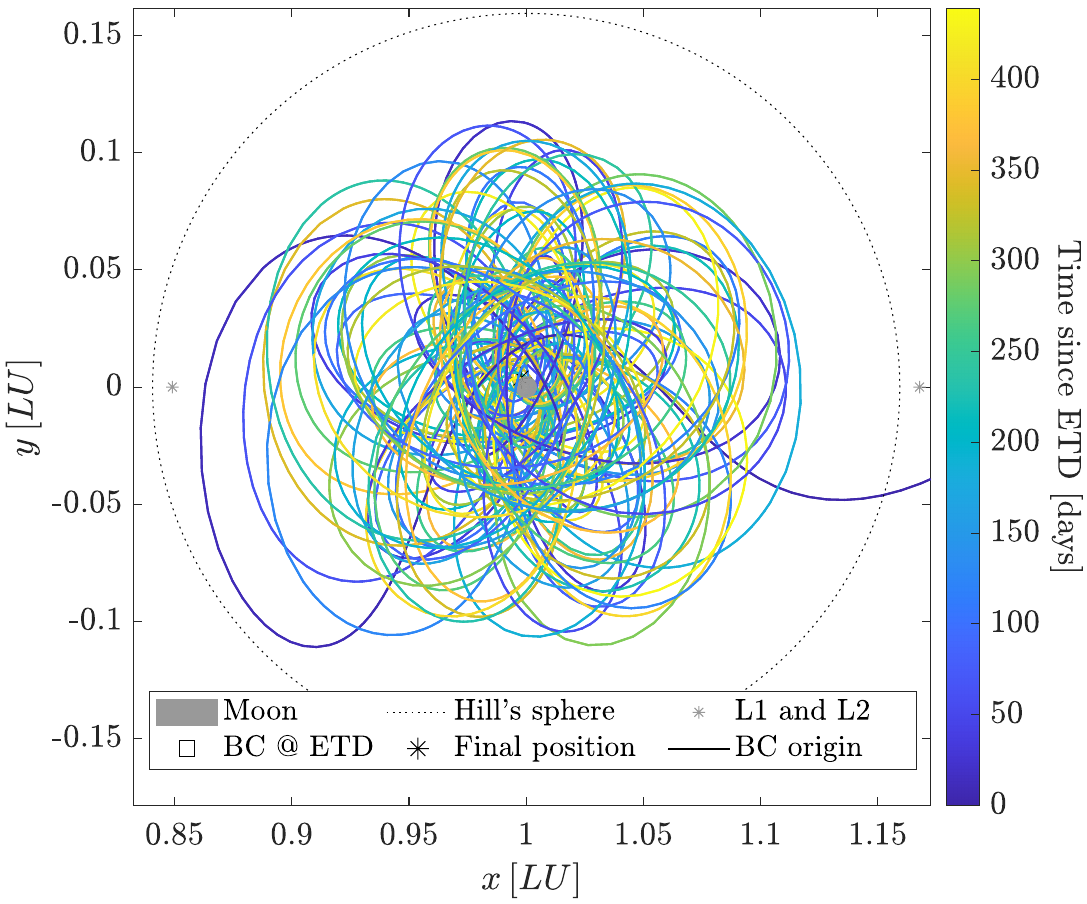}
    \end{subfigure}
    \\
    \begin{subfigure}[t!]{0.49\textwidth}
        \includegraphics[width=\textwidth]{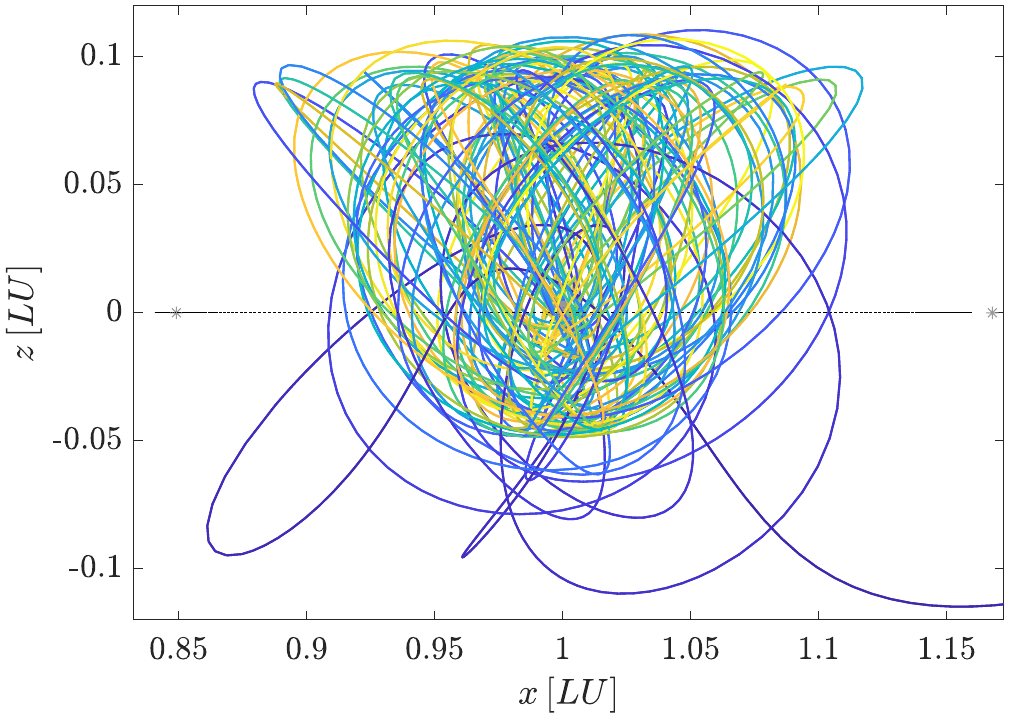}
    \end{subfigure}
    \caption{Longest \gls{bc} in the synodic frame.}
    \label{fig: longest BC syn}
\end{figure}

\begin{figure}[tbp]
    \centering
    \includegraphics[width=0.63\textwidth]{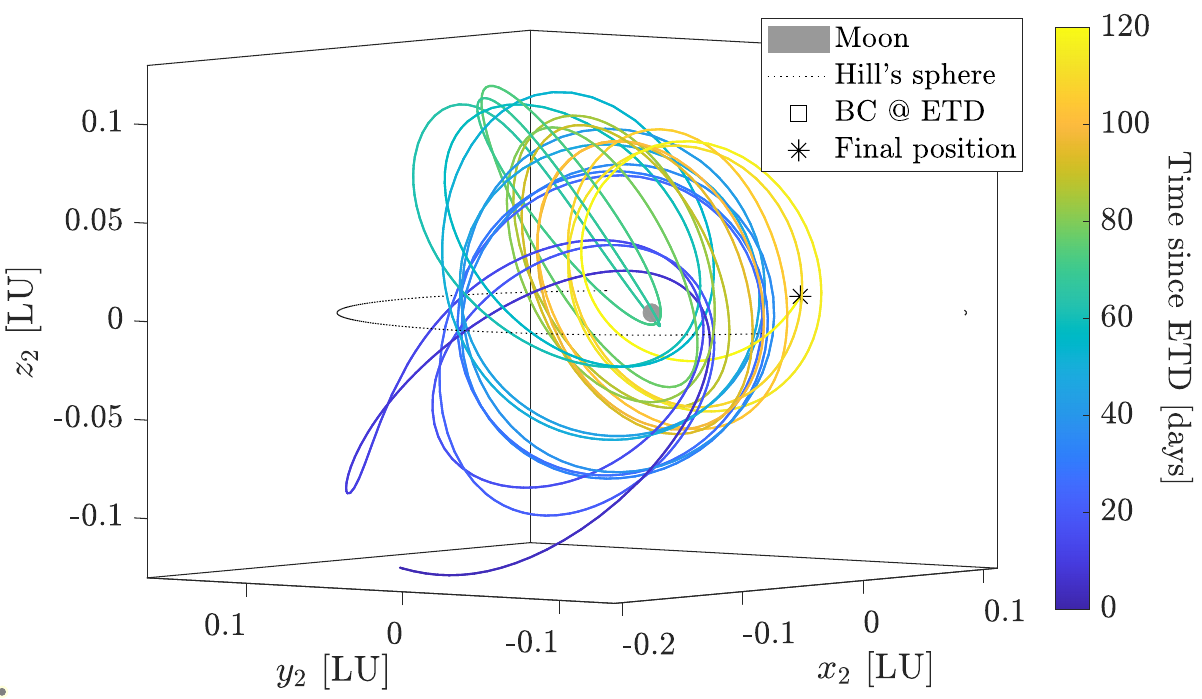}
    \caption{Longest \gls{bc} in Moon inertial frame.}
    \label{fig: longest BC inert2}
\end{figure}

A \gls{bc} belonging to the corridor is extracted and shown in Figs.~\ref{fig: longest BC syn} and \ref{fig: longest BC inert2}. This trajectory lasts 439 days (approximately 14.6 months) before ultimately impacting the Moon. While retaining the previously discussed features, its extended propagation time allows for observing additional dynamical characteristics. The post-transition dynamics appear to follow a repeating pattern in the synodic frame (\cref{fig: longest BC syn}). This suggests the presence of an underlying dynamical structure associated with a periodic orbit in the \gls{cr3bp}. In particular, strong similarities can be observed with (Period-$n$) Distant Retrograde Orbits (DROs) and butterfly periodic orbits. This \gls{bc} exhibits features of both: the $x$–$y$ view indicates a periodic retrograde motion reminiscent of DROs, while the $x$–$z$ projection bears resemblance to a northern butterfly orbit around $C_J = 3.09$ (see Appendix\ref{appendix: sample ICs}).
The \gls{bc} resembles a bounded, though highly perturbed, lunar orbit. As shown in \cref{fig: longest BC inert2}, all orbital elements exhibit fluctuations. Their amplitudes varying significantly over time. Nevertheless, chaotic phases are interspersed with more stable ones. For example, two relatively stable intervals can be identified between days 25–60 and 80–120, characterized primarily by the precession of the \gls{raan} $\Omega$ and the argument of perilune $\omega$.

\begin{figure}[tbp]
    \centering
    \begin{subfigure}[t!]{0.49\textwidth}
        \includegraphics[width=\textwidth]{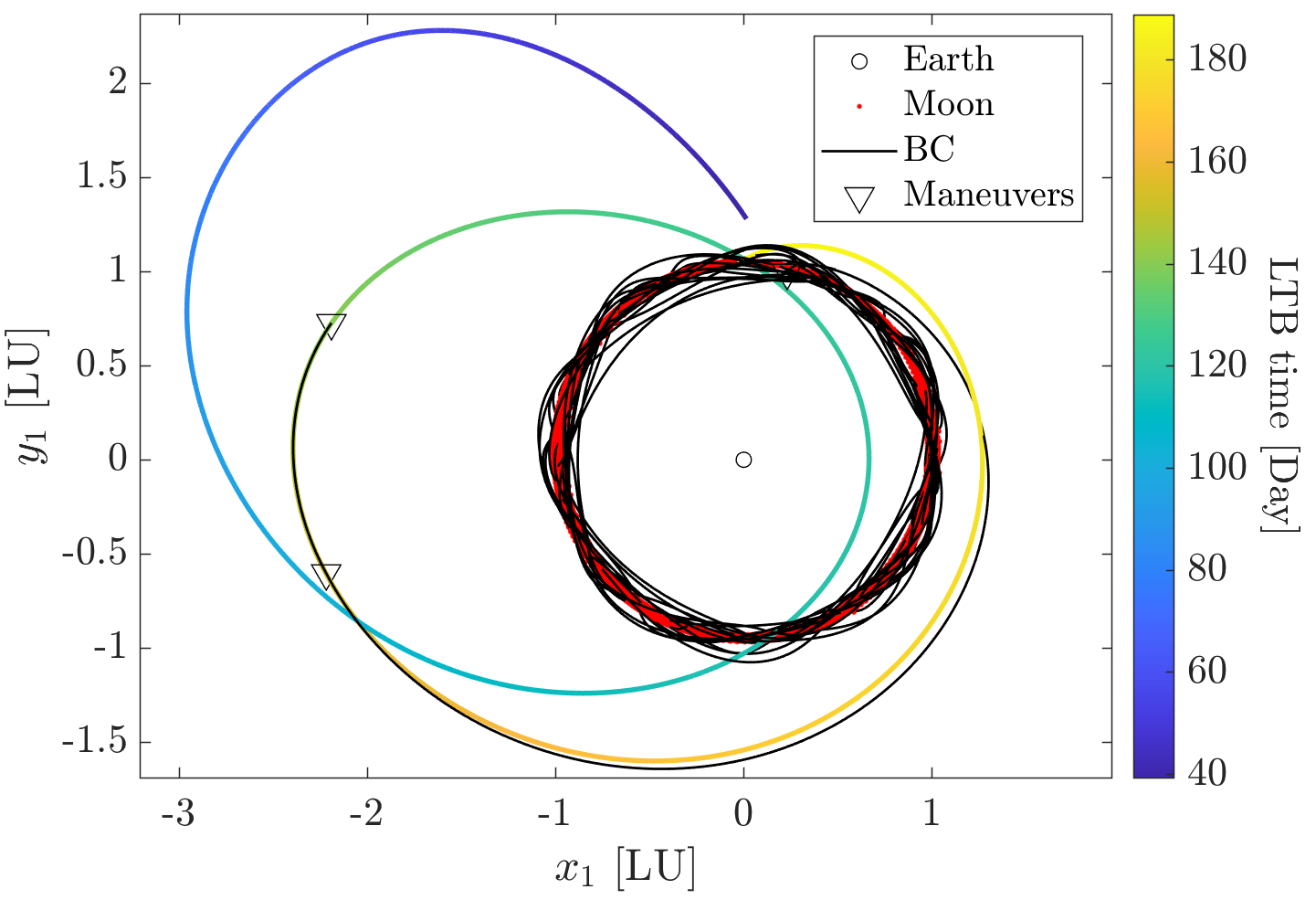}
    \end{subfigure}
    \\
    \begin{subfigure}[t!]{0.49\textwidth}
        \includegraphics[width=\textwidth]{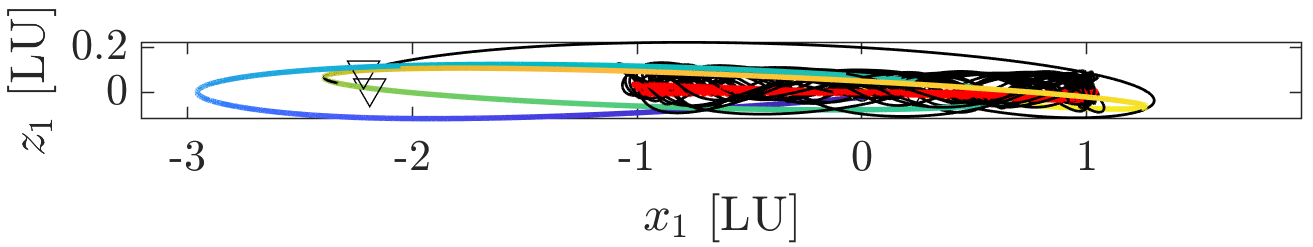}
    \end{subfigure}
    \caption{Distance metric assessment: \gls{ltb} to longest \gls{bc} three-impulses transfer in Earth inertial frame.}
    \label{fig: longest BC inert1 cost}
\end{figure}

The estimated mono-impulsive transfer cost, as given by the distance metric, is $d_v = 44.2$~m/s. To assess its accuracy, a simple optimization tool was implemented using MATLAB\copyright's \textit{fmincon} to compute a three-impulse transfer from \gls{ltb} to this \gls{bc}, minimizing the total $\Delta v$. The resulting transfer, shown in \cref{fig: longest BC inert1 cost}, requires $\Delta v = 66.1$~m/s. The transfer is nearly bi-impulsive, with a negligible first correction. The majority of the correction — $59.8$~m/s — is applied at the intermediate maneuver, while the remaining $6.3$~m/s are used during the approach to the Moon.

Even longer \glspl{bc}, lasting up to 4 years, can be found when extending the search to $d_v < 85$~m/s. 

\subsubsection{Stable polar BC} \label{sec: stable polar BC}
\begin{figure}[tbp]
    \centering
    \begin{subfigure}[t!]{0.49\textwidth}
        \includegraphics[width=\textwidth]{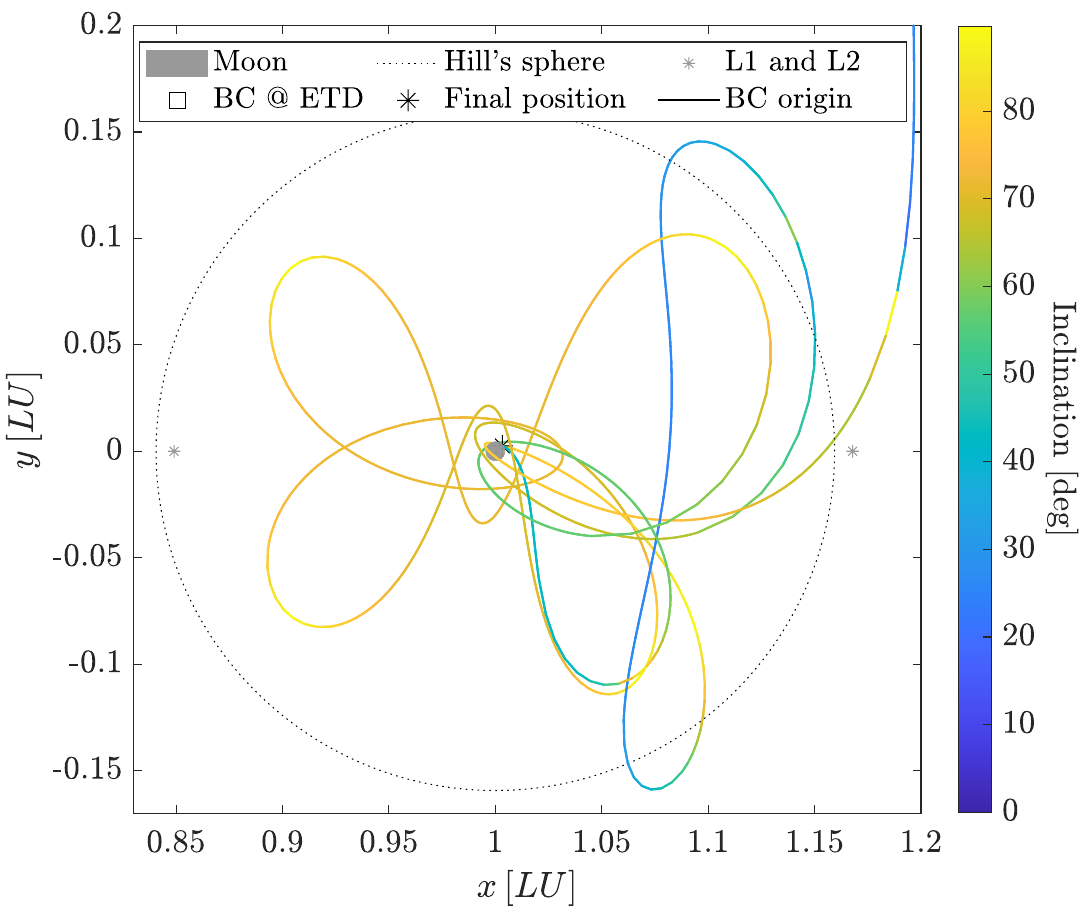}
    \end{subfigure}
    \\
    \begin{subfigure}[t!]{0.49\textwidth}
        \includegraphics[width=\textwidth]{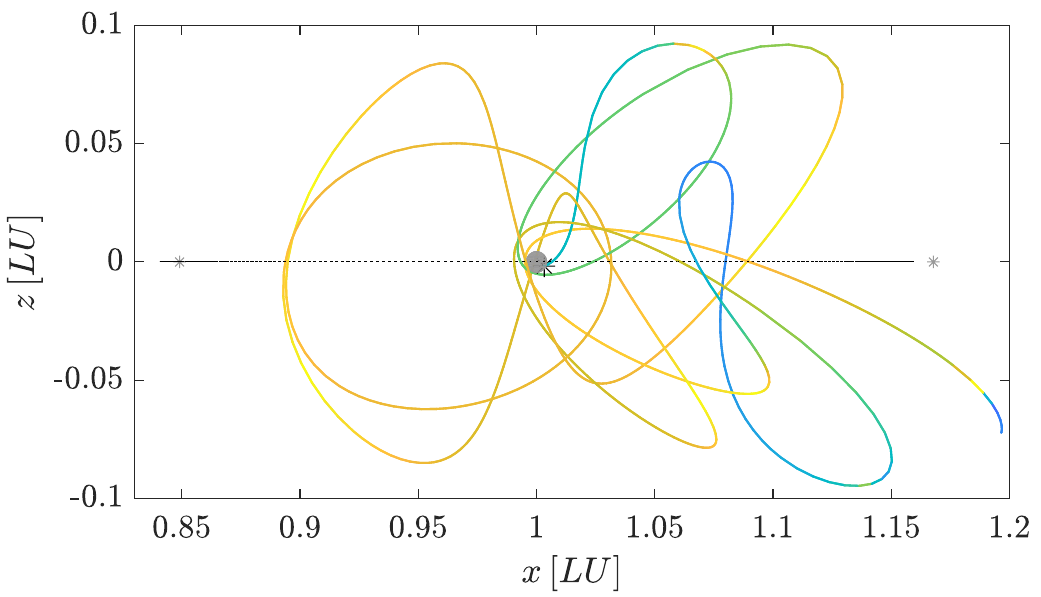}
    \end{subfigure}
    \caption{Stable polar \gls{bc} in the synodic frame.}
    \label{fig: stable polar BC syn}
\end{figure}

\begin{figure}[tbp]
    \centering
    \includegraphics[width=0.63\textwidth]{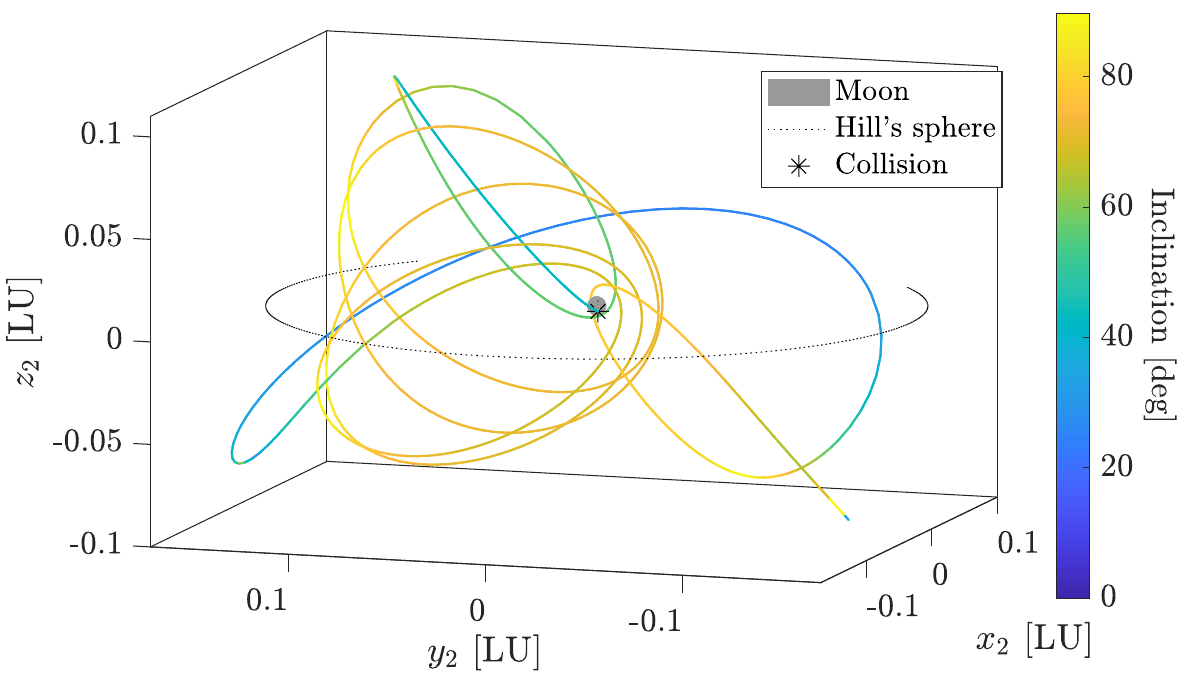}
    \caption{Stable polar \gls{bc} in Moon inertial frame.}
    \label{fig: stable polar BC inert2}
\end{figure}

Another \gls{bc} is extracted from the database and shown in Figs.~\ref{fig: stable polar BC syn} and \ref{fig: stable polar BC inert2}. This trajectory is selected by minimizing the value of $d_v$ among all \glspl{bc} that complete at least four revolutions and exhibit a high inclination at the first perilune, specifically those satisfying $\left|i_{1stRev}-90^{\circ}\right|<20^{\circ}$.
The capture begins with an initial leg of approximately 8 days, culminating in a highly inclined retrograde perilune passage (see \cref{fig: stable polar BC syn}). This is followed by a 17-day leg, during which the spacecraft travels far from the Moon along a trajectory resembling a spatial periodic orbit of Hénon's \textit{g'} family~\cite{henon1969V}. This type of motion, observed in several other \glspl{bc} within the dataset, appears to facilitate a transition into a more stable capture phase.
Between days 25 and 45, the trajectory exhibits behavior similar to that seen in the previous example: a relatively stable orbital phase dominated by slow precession of the \gls{raan} $\Omega$ and the argument of perilune $\omega$ (see \cref{fig: stable polar BC inert2}). Eventually, a new transition occurs, leading to a lunar impact around day 55.


\subsubsection{Multiple polar insertion opportunities} \label{sec: multiple polar insertions}
\begin{figure}[tbp]
    \centering
    \begin{subfigure}[t!]{0.49\textwidth}
        \includegraphics[width=\textwidth]{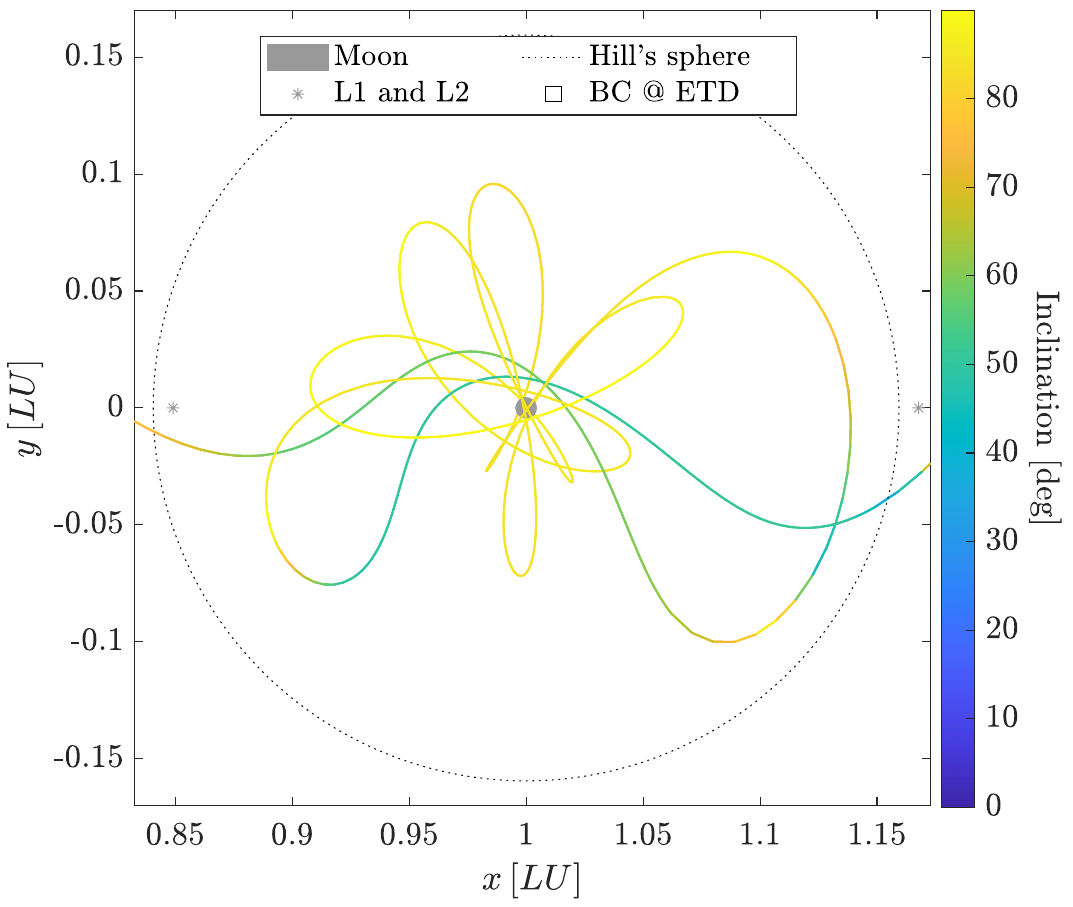}
    \end{subfigure}
    \\
    \begin{subfigure}[t!]{0.49\textwidth}
        \includegraphics[width=\textwidth]{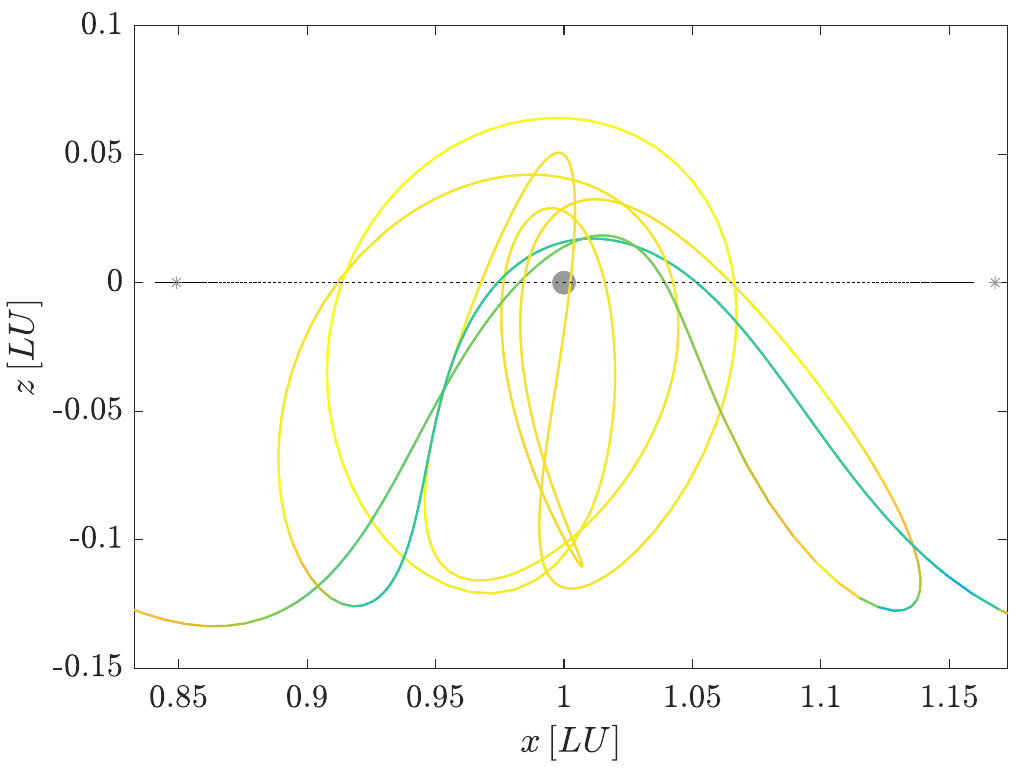}
    \end{subfigure}
    \caption{Multiple polar insertion opportunities \gls{bc} in the synodic frame.}
    \label{fig: similar polar periselenia BC syn}
\end{figure}

\begin{figure}[tbp]
    \centering
    \includegraphics[width=0.63\textwidth]{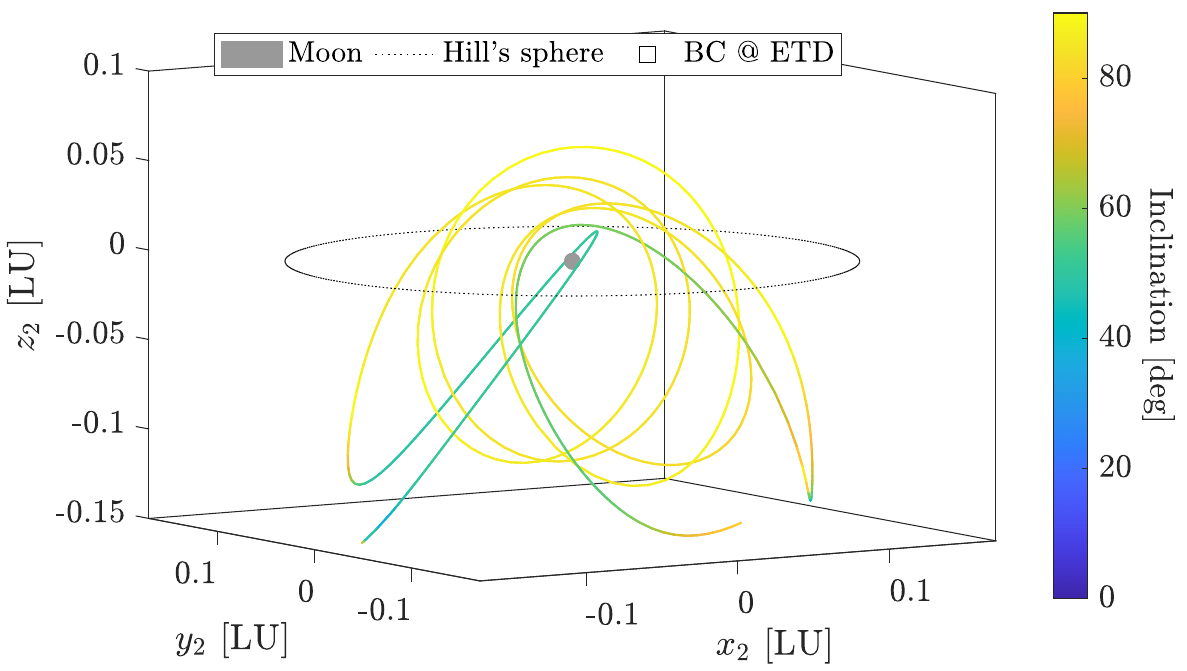}
    \caption{Multiple polar insertion opportunities \gls{bc} in Moon inertial frame.}
    \label{fig: similar polar periselenia BC inert2}
\end{figure}

A third \gls{bc} is extracted from the database and illustrated in Figs.~\ref{fig: similar polar periselenia BC syn} and \ref{fig: similar polar periselenia BC inert2}. This trajectory is selected by minimizing a newly defined distance metric, $d_{v,2}$, which compares the orbital elements at three stored perilunes along each trajectory: the first perilune $\mathbf{r}_{2,1stRev}$ with inclination $i_{1stRev}$, the second $\mathbf{r}_{2,min1}$ with $i_{r_{min1}}$, and the third $\mathbf{r}_{2,min2}$ with $i_{r_{min2}}$.
The metric $d_{v,2}$ follows the same formulation as the previously introduced $d_v$, but it is computed with the Moon as the central body. It is evaluated for all pairwise combinations of the three perilunes, and the pair yielding the lowest $d_{v,2}$ value is selected for further analysis.
As part of the selection process, an additional constraint is applied: the inclination $i$ at the two perilunes must satisfy $\left|i-90^{\circ}\right|<6^{\circ}$. This ensures that the selected pair corresponds to nearly polar perilunes. As a result, the selected \gls{bc} features two closely matched perilunes, both suitable for a polar \gls{loi}.


\subsubsection{Multiple grazing and polar insertion opportunities} \label{sec: multiple grazing polar}
\begin{figure}[tbp]
    \centering
    \begin{subfigure}[t!]{0.49\textwidth}
        \includegraphics[width=\textwidth]{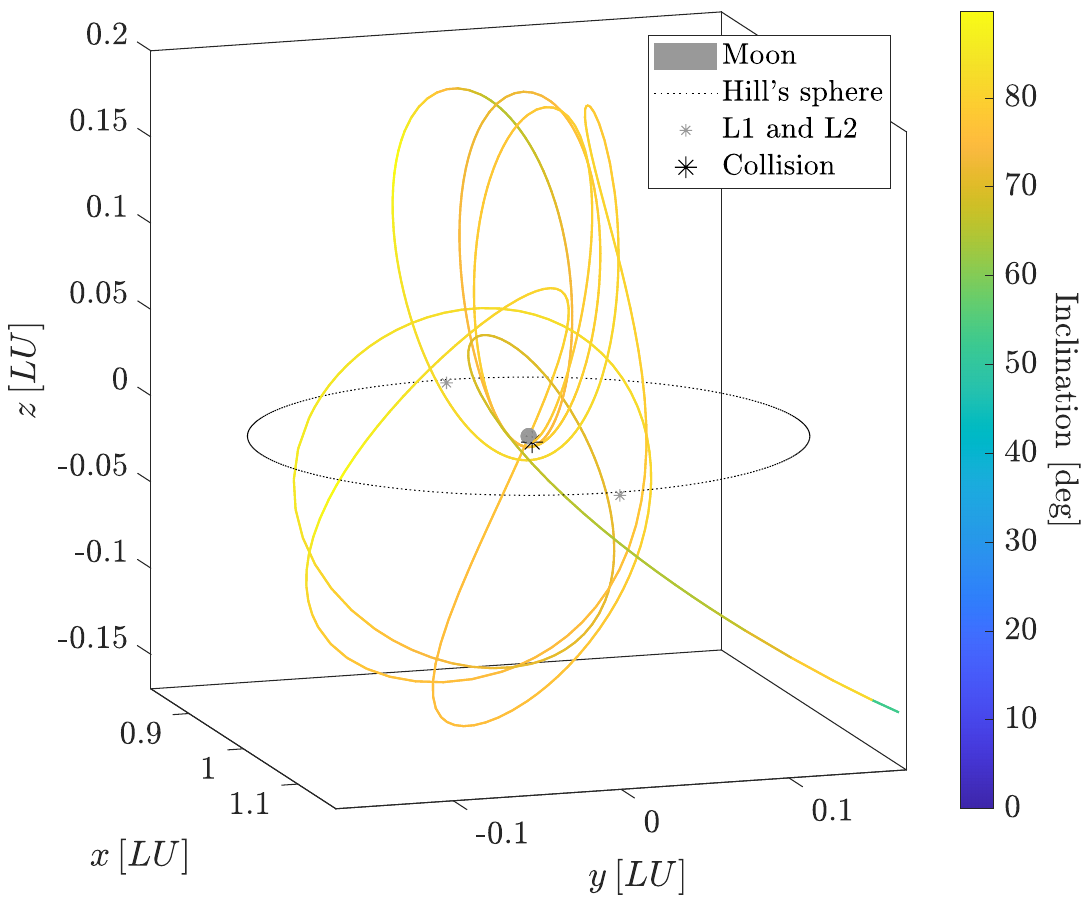}
    \end{subfigure}
    \begin{subfigure}[t!]{0.49\textwidth}
        \includegraphics[width=\textwidth]{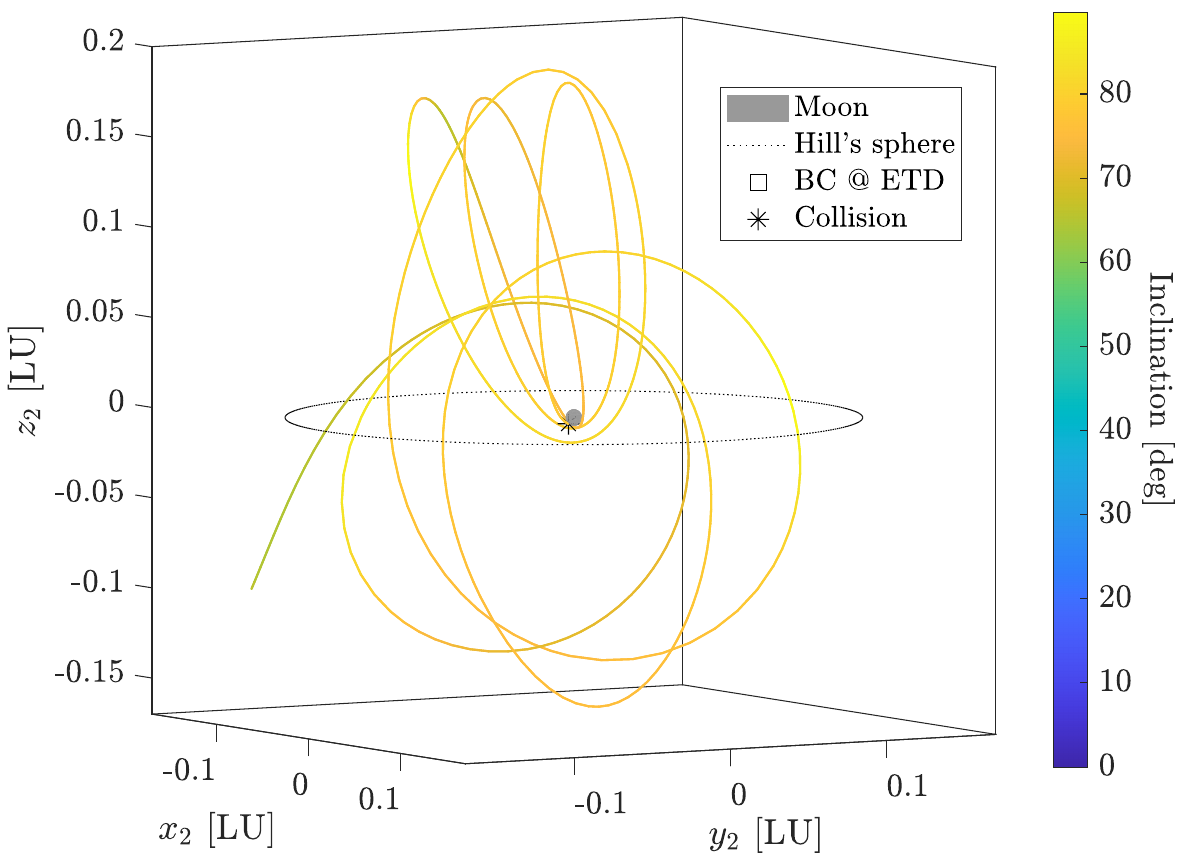}
    \end{subfigure}
    \caption{Multiple grazing polar insertion opportunities \gls{bc} in the synodic (left) and Moon inertial (right) frames.}
    \label{fig: similar grazing polar periselenia BC syn and inert2}
\end{figure}

A fourth \gls{bc} is extracted and displayed in \cref{fig: similar grazing polar periselenia BC syn and inert2}.
As in the previous case, the metric $d_{v,2}$ is used, with a constraint on high inclination: $\left|i-90^{\circ}\right|<12^{\circ}$. Additionally, the distance from the Moon at both selected perilunes must satisfy $r_2 < r_M + 600$~km.
This yields a \gls{bc} with two similar insertion opportunities, both polar and approximately $520$~km above the Moon. This portion of the \gls{bc} closely resembles a Near-Rectilinear Halo Orbit (NRHO), making this approach particularly appealing.

\begin{figure}[tbp]
    \centering
    \begin{subfigure}[t!]{0.49\textwidth}
        \includegraphics[width=\textwidth]{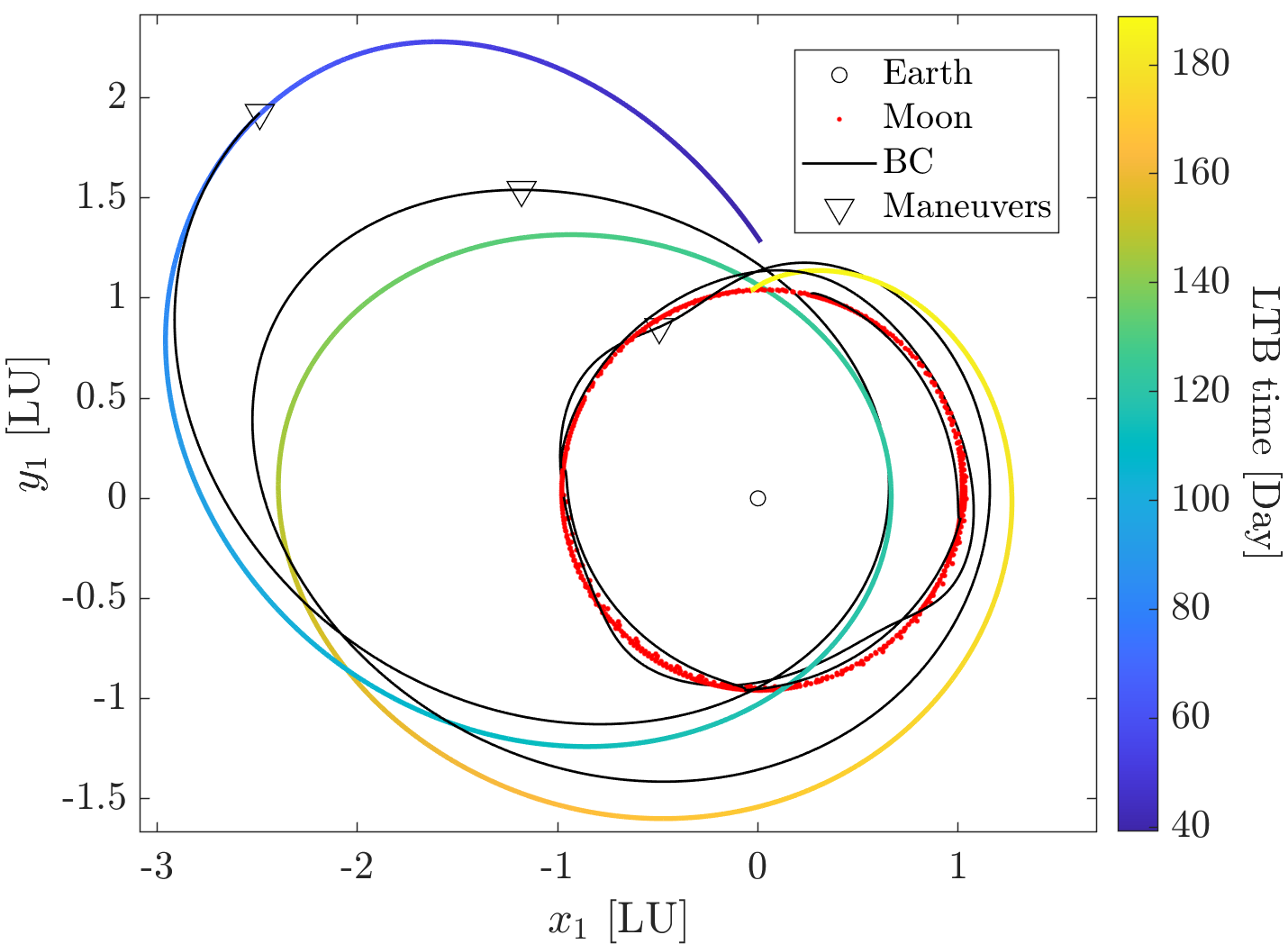}
    \end{subfigure}
    \\
    \begin{subfigure}[t!]{0.49\textwidth}
        \includegraphics[width=\textwidth]{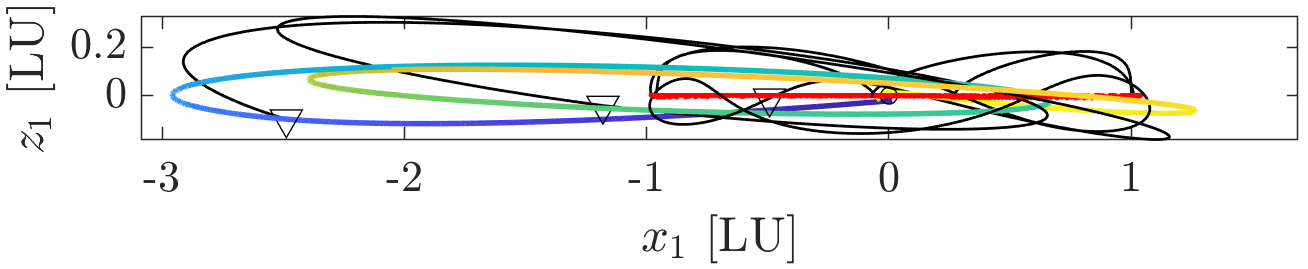}
    \end{subfigure}
    \caption{Distance metric assessment: \gls{ltb} to grazing, polar \gls{bc} three-impulses transfer in Earth inertial frame..}
    \label{fig: similar grazing polar periselenia BC inert1 cost}
\end{figure}

The estimated mono-impulsive transfer cost given by the distance metric is $d_v = 57.4$~m/s. To assess the accuracy of this estimate, a three-impulse transfer from \gls{ltb} to the selected \gls{bc} is optimized. The resulting trajectory is shown in \cref{fig: similar grazing polar periselenia BC inert1 cost} and requires a total $\Delta v$ of $65.8$~m/s. Notably, the majority of the cost stems from the inclination change (see $x_1$–$z_1$ view), with the first maneuver contributing approximately $46.4$~m/s. The final correction, performed during the Moon approach, accounts for just $1.9$~m/s.

\subsubsection{Stabilization through small successive maneuvers} \label{sec: successive braking man}
\begin{figure}[tbp]
    \centering
    \begin{subfigure}[t!]{0.50\textwidth}
        \includegraphics[width=\textwidth]{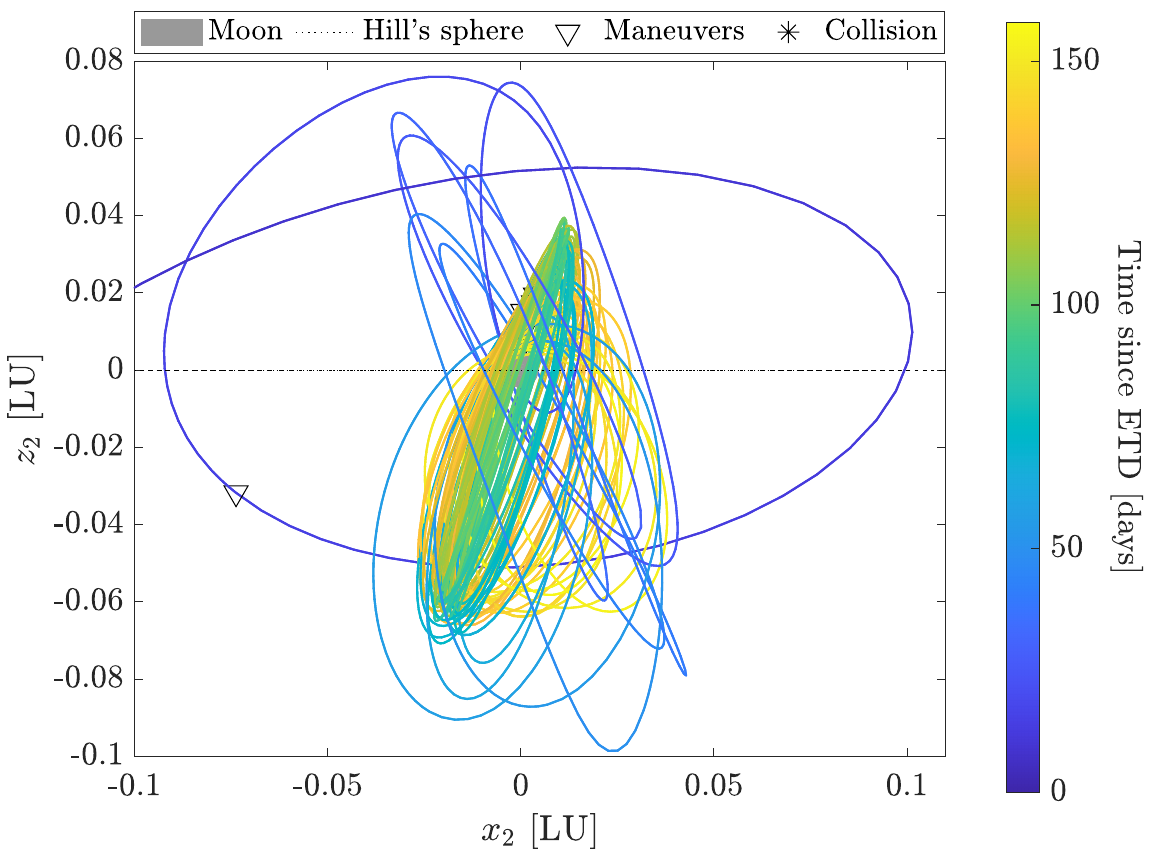}
    \end{subfigure}
    \hfill
    \begin{subfigure}[t!]{0.48\textwidth}
        \includegraphics[width=\textwidth]{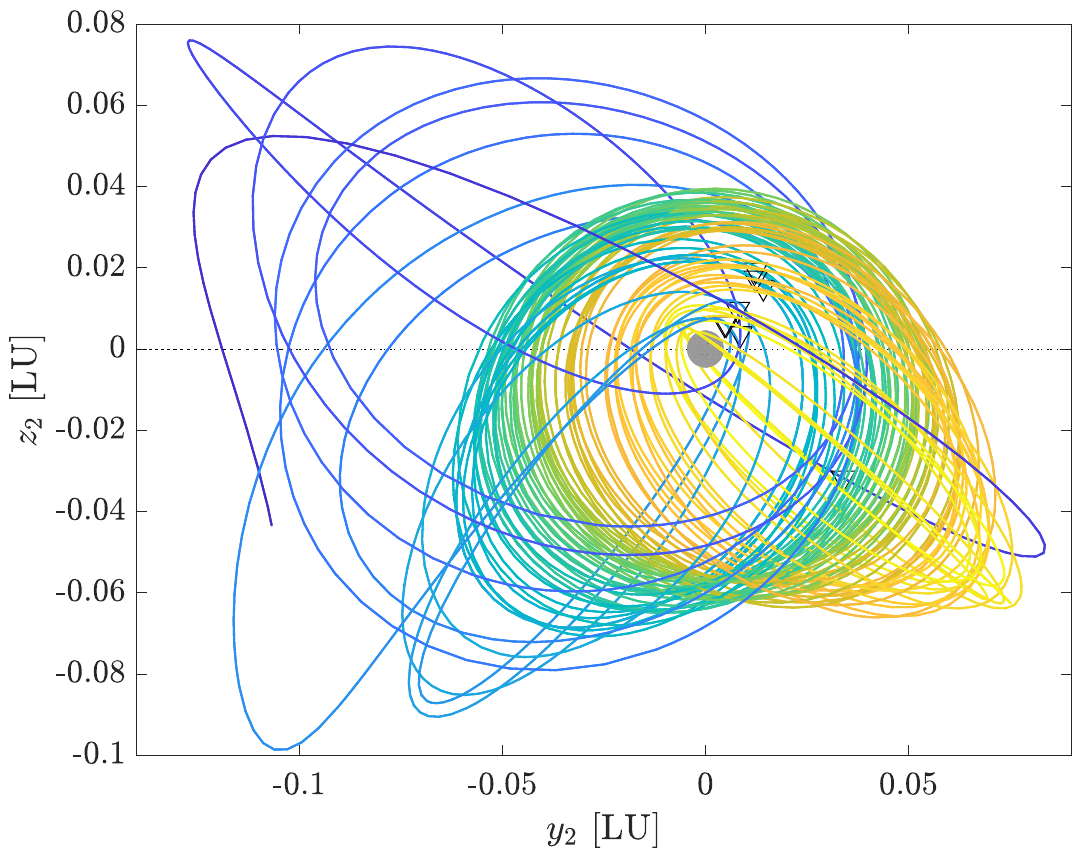}
    \end{subfigure}
    \caption{Successive braking maneuvers in Moon inertial frame.}
    \label{fig: braking man}
\end{figure}

A sample polar \gls{bc} is selected to illustrate the feasibility of stabilizing this class of chaotic trajectories. A sequence of small braking maneuvers, each imparting just $5$~m/s opposite to the velocity vector, is applied at successive lunar perilunes. The objective is to progressively extend its orbital lifetime, increasingly stabilizing the motion within the lunar environment.
An illustrative case is shown in \cref{fig: braking man}, where seven maneuvers are performed for a total $\Delta v$ of $35$~m/s. After an initial transition phase lasting nearly 40 days, the trajectory settles into a lower-energy orbit with reduced semimajor axis, maintained between days 40 and 140. This stabilized phase exhibits only mild precession in the argument of periapsis $\omega$ and the longitude of the ascending node $\Omega$. The final maneuver is executed around day 70 after the \gls{etd}, followed by approximately 80 days of continued lunar orbiting before eventual impact with the Moon.
Although this example is not optimized, it clearly demonstrates the viability of such an approach. A more refined stabilization strategy targeting maneuvers closer to the Moon could increase efficiency while preventing future collisions.

\section{Conclusions}
This work presents a comprehensive framework for the generation and analysis of Ballistic Capture (BC) trajectories in the spatial Circular Restricted Three-Body Problem (CR3BP). In the first part, a method is developed to compute a complete database of BC trajectories, referred to as the capture set. This objective is achieved by extending the previously introduced Energy Transition Domain (ETD) concept—originally defined in the planar case—to three dimensions. By using initial conditions from the ETD, billions of candidate BCs are retrieved for the Earth–Moon system, and their key features are stored for subsequent analysis. Although demonstrated here for the Earth–Moon system, the methodology is directly applicable to other low-energy lunar transfers and planetary systems.

Secondly, a mission-specific distance metric is introduced to refine the selection of promising BC candidates. This filtering step enables the isolation of a smaller subset of trajectories that are then transitioned into a full ephemeris model, making the process computationally feasible. The methodology is illustrated through application to NASA’s Lunar Trailblazer mission, showing how low-energy insertion trajectories can be systematically identified and refined for specific mission requirements.

Finally, this work focuses on analyzing the refined subset of ephemeris trajectories tailored to Lunar Trailblazer. After providing an overall view of the database, representative examples are presented, including trajectories with multiple, polar, and repeated close approaches to the Moon. Particular attention is given to BCs offering multiple successive insertion opportunities or allowing the insertion maneuver to be distributed across several perilunes—characteristics that enhance robustness during the critical insertion phase. Candidate trajectories display minimal deviation from the nominal Trailblazer trajectory, as indicated by low values of the distance metric, suggesting their potential as backup options. An a posteriori evaluation confirms the usefulness of the distance metric for pre-selecting viable trajectories.

The examples shown, selected based on simple constraints, provide an initial demonstration of how the capture set can be exploited to identify alternative insertion opportunities. However, fully leveraging the potential of the capture set requires more comprehensive investigations. Future work could involve advanced clustering or machine learning techniques~\cite{Miceli_clusteringLLFO, TemporaryCaptureTaxonomy_MachineLearning} to process the large dataset and extract meaningful patterns, or the development of optimization-based methods using tailored cost functions. Such approaches would further refine the BC candidates and enhance their suitability for precise mission planning.

Ultimately, this framework lays a scalable foundation for bridging the gap between dynamical systems theory and real-world trajectory design, facilitating the integration of ballistic capture strategies into future low-energy missions.

\section*{Funding Sources}
The research was carried out at the Jet Propulsion Laboratory, California Institute of Technology, under a contract with the National Aeronautics and Space Administration (80NM0018D0004).

Lorenzo Anoè's visit at NASA Jet Propulsion Laboratory (JPL) through the JPL Visiting Student Research Program (JVSRP) has been funded by the Royal Society Te Aparangi Catalyst Seeding general grant: Advanced Cislunar Space Mission Design.

\section*{Acknowledgments}
The authors wish to acknowledge the Centre for eResearch at the University of Auckland for their assistance in facilitating this research. http://www.eresearch.auckland.ac.nz

\bibliography{sample}

\pagebreak
\begin{appendices}
\section{Sample trajectories initial conditions} \label{appendix: sample ICs}
\begin{figure}[tbp]
    \centering
    \includegraphics[width=0.63\textwidth]{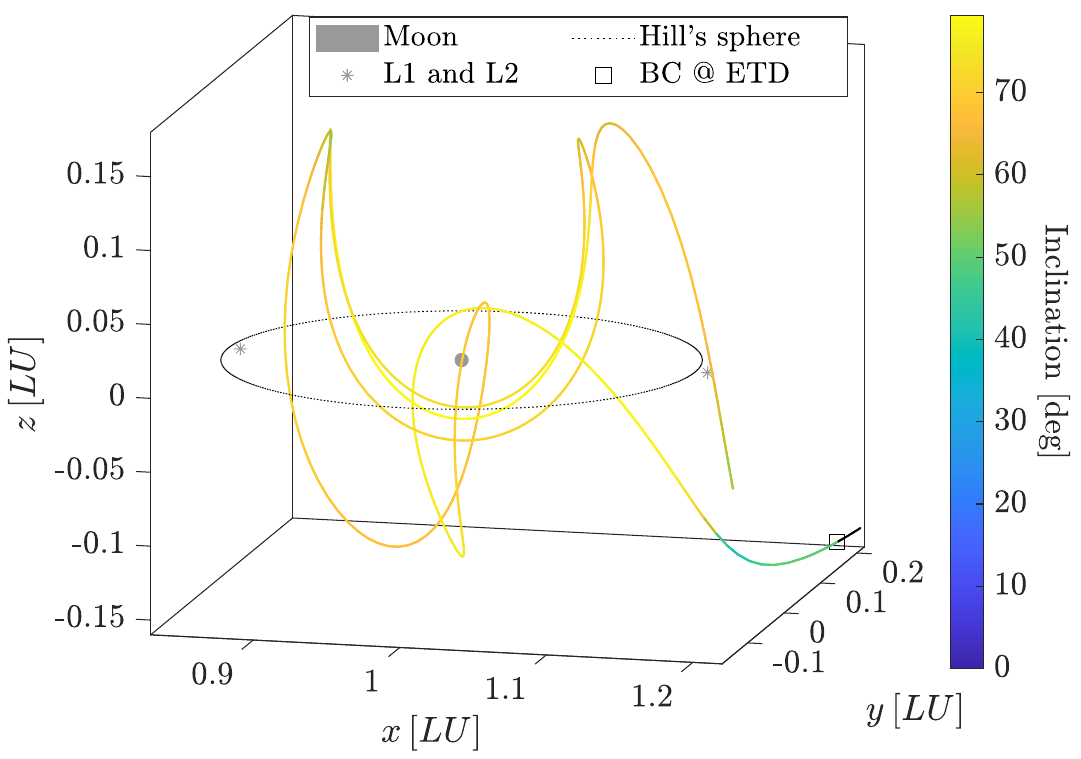}
    \caption{Butterfly-like \gls{bc} in the synodic frame.}
    \label{fig: butterfly BC syn}
\end{figure}

In the following tables, the initial conditions for the sample \gls{bc} trajectories given in this work are given.
Firstly, \cref{tab: parameters systems} contains masses and primaries distance used to set up the \gls{cr3bp} and ephemeris models for the Earth-Moon system.
Then, Tables \ref{tab: ICs ephemeris 1} and \ref{tab: ICs ephemeris 2} express the dimensional initial conditions for the sample \glspl{bc} of \cref{sec: sample BCs for LTB} in the EMO2000 frame at time $T_{ETD}$. An additional initial condition is provided for a \gls{bc} resembling a northern butterfly periodic orbit for $C_J \sim 3.09$, which is shown in \cref{fig: butterfly BC syn}.

\begin{table}[htbp]
\caption{Bodies masses $m_S$, $m_E$ and $m_M$, Earth-Moon relative distance $r_{EM}$, Moon physical radius $r_M$, and epoch (in both ephemeris seconds past J2000 $T_{ETD}$ and calendar format $\text{date}_{ETD}$) for ephemeris models implementation}
\centering{}\label{tab: parameters systems}
\begin{tabular}{ccccccc}
\hline
\noalign{\vskip\doublerulesep}
$m_S$ [kg] & $m_E$ [kg] & $m_M$ [kg] & $r_{EM}$ [km] & $r_M$ [km] & $T_{ETD}$ [s] & $\text{date}_{ETD}$ \tabularnewline[\doublerulesep]
\hline
\noalign{\vskip\doublerulesep}
\noalign{\vskip\doublerulesep}
$1.9885 \cdot 10^{30}$  & $5.9724 \cdot 10^{24}$ & $7.3461 \cdot 10^{22}$ & $384399$ & $1737.4$ & $802221652.5$ & 2025 JUN 03 11:19:43.3 \tabularnewline[\doublerulesep]
\noalign{\vskip\doublerulesep}
\end{tabular}
\end{table}


\begin{table}[htbp]
\caption{Initial conditions for sample orbits in the ephemeris model at epoch $T_{ETD}$: part 1}
\centering{}\label{tab: ICs ephemeris 1}
\begin{tabular}{cccc}
\hline
\noalign{\vskip\doublerulesep}
Variable & \cref{sec: longest BCs} \gls{bc} & \cref{sec: stable polar BC} \gls{bc} & \cref{sec: multiple polar insertions} \gls{bc} \tabularnewline[\doublerulesep]
\hline
\noalign{\vskip\doublerulesep}
\noalign{\vskip\doublerulesep}
$x_0$ [km] & $-500754.648873973$ & $-485952.557622184$ & $-502151.104316433$ \tabularnewline[\doublerulesep]
\noalign{\vskip\doublerulesep}
\noalign{\vskip\doublerulesep}
$y_0$ [km] & $96930.0726983651$ & $12484.7053447739$ & $67890.4520791561$ \tabularnewline[\doublerulesep]
\noalign{\vskip\doublerulesep}
\noalign{\vskip\doublerulesep}
$z_0$ [km] & $-28315.7473944248$ & $-32398.9385774915$ & $-50012.4217631616$ \tabularnewline[\doublerulesep]
\noalign{\vskip\doublerulesep}
\noalign{\vskip\doublerulesep}
$\dot{x}_0$ [km/s] & $-0.0262302667192358$ & $-0.0290637180948451$  & $-0.0279032774486339$ \tabularnewline[\doublerulesep]
\noalign{\vskip\doublerulesep}
\noalign{\vskip\doublerulesep}
$\dot{y}_0$ [km/s] & $-0.966278607253216$ & $-0.972684625927066$ & $-0.942356982288852$ \tabularnewline[\doublerulesep]
\noalign{\vskip\doublerulesep}
\noalign{\vskip\doublerulesep}
$\dot{z}_0$ [km/s] & $-0.182778054922072$ & $-0.0988095375176495$ & $-0.141908744346009$ \tabularnewline[\doublerulesep]
\noalign{\vskip\doublerulesep}
\end{tabular}
\end{table}

\begin{table}[htbp]
\caption{Initial conditions for sample orbits in the ephemeris model at epoch $T_{ETD}$: part 2}
\centering{}\label{tab: ICs ephemeris 2}
\begin{tabular}{cccc}
\hline
\noalign{\vskip\doublerulesep}
Variable & \cref{sec: multiple grazing polar} \gls{bc} & \cref{sec: successive braking man} \gls{bc} & \cref{fig: butterfly BC syn} butterfly-like \gls{bc} \tabularnewline[\doublerulesep]
\hline
\noalign{\vskip\doublerulesep}
\noalign{\vskip\doublerulesep}
$x_0$ [km] & $-456081.713990439$ & $-509026.731598873$ & $-483653.619368984$ \tabularnewline[\doublerulesep]
\noalign{\vskip\doublerulesep}
\noalign{\vskip\doublerulesep}
$y_0$ [km] & $-2451.70963369324$ & $46561.0023634826$ & $43594.8128371648$ \tabularnewline[\doublerulesep]
\noalign{\vskip\doublerulesep}
\noalign{\vskip\doublerulesep}
$z_0$ [km] & $-76462.9670875461$ & $-42387.9887255228$ & $-61529.9871265759$ \tabularnewline[\doublerulesep]
\noalign{\vskip\doublerulesep}
\noalign{\vskip\doublerulesep}
$\dot{x}_0$ [km/s] & $-0.0366896820620522$ & $-0.0351690791481468$  & $-0.0272893185001766$ \tabularnewline[\doublerulesep]
\noalign{\vskip\doublerulesep}
\noalign{\vskip\doublerulesep}
$\dot{y}_0$ [km/s] & $-0.983876646961336$ & $-0.961755240640447$ & $-0.996138672414697$ \tabularnewline[\doublerulesep]
\noalign{\vskip\doublerulesep}
\noalign{\vskip\doublerulesep}
$\dot{z}_0$ [km/s] & $-0.0958194874764763$ & $-0.0900457348096386$ & $-0.129261358637946$ \tabularnewline[\doublerulesep]
\noalign{\vskip\doublerulesep}
\end{tabular}
\end{table}

\section{Additional spatial characteristics of the ETD} \label{appendix: etd characteristics}
Building on the discussion of \cref{sec: ETD}, we can derive additional analytical properties of the \gls{etd}.
First of all, the \gls{etd} is never defined exactly in the $M_2$ position $(x, y, z)=(1-\mu,0,0)$, because \cref{eq: v2_energy_constraint} is not defined in that point (where $r_2=0$). In addition, a singularity always occurs at $(x, y, z)=(1-\mu,0,z)$, located directly above or below $M_2$. Here, the \gls{etd} is never defined apart from a value of
\begin{equation}
    C_J=f(z)=1-2\mu+\mu^2+2(1-\mu)/\sqrt{1+z^2} \, ,
    \label{eq: CJ above M2}
\end{equation}
which has a minimum for $z=0$ with the aforementioned value $C_J=3-4\mu+\mu^2$. This entails that for each value of $z$ above or below $M_2$ the \gls{etd} is always defined, but for a specific value of $C_J$ only. This value for $C_J$ decreases (and hence $\Gamma$ increases) for increasing $z$, and this feature can be observed in all the figures with $\Gamma\geq1.2$ (\cref{fig: ETD dominio all 2}).
In the geometric representation using spheres, the condition $(x,y,z)=(1-\mu,0,z)$ yields $C_1=\{0,0,0\}$, and it can be represented by two concentric spheres. As a consequence, the two radii $r_J$ and $r_\varepsilon$ must coincide to provide a solution for the \gls{etd}. When this occurs, the solution space does not consist of a one-dimensional circumference, like in the usual case. Instead, it is a two-dimensional surface defined by the two overlapping spheres.
In other words, for the points with coordinates $(x,y,z)=(1-\mu,0,z)$, the two constraints imposed (i.e. two-body energy $\varepsilon_2=0$ and three-body energy $C_J=\text{const}$) are always incompatible except when the value is the on in \cref{eq: CJ above M2}. In the latter case, when $C_J$ takes exactly that value, the two constraints coincide, and the dimension of the solution space is two. This means that all the possible combinations of angles $(\eta,\zeta)$ are solutions, instead of the correlation $\eta=f(\zeta)$ introduced in \cref{sec: ics in ETD}.

As previously highlighted in~\cite{BC-journal}, a region where $\varepsilon_2$ remains strictly negative surrounds $M_2$ (i.e., the Moon) for $C_J>3-4\mu+\mu^2$. However, for lower values of the Jacobi constant (i.e., higher three-body energies), such negative-only regions cease to exist.
\end{appendices}

\end{document}